\documentclass{amsart}

\usepackage{amsmath,amssymb,graphics}
\usepackage[headings]{fullpage}

\makeatletter
\def\ps@headings{\ps@empty
  \def\@evenhead{%
    \setTrue{runhead}%
    \normalfont\scriptsize
    \rlap{\thepage}\hfil
    \def\thanks{\protect\thanks@warning}%
    \leftmark{}{}\hfil}%
  \def\@oddhead{%
    \setTrue{runhead}%
    \normalfont\scriptsize \hfil
    \def\thanks{\protect\thanks@warning}%
    \rightmark{}{}\hfil \llap{\thepage}}%
  \let\@mkboth\markboth
}
\makeatother

\usepackage[all]{xy}

\usepackage{fullpage}

\theoremstyle{plain}
\newtheorem{theorem}{Theorem}[section]
\newtheorem{proposition}[theorem]{Proposition}

\newtheorem{fact}[theorem]{Fact}
\newtheorem{facts}[theorem]{Facts}

\theoremstyle{definition}
\newtheorem{definition}[theorem]{Definition}
\newtheorem{example}[theorem]{Example}
\newtheorem{remark}[theorem]{Remark}

\newtheorem{observation}[theorem]{Observation}

\numberwithin{equation}{section}

\def\gg {\mathfrak{g}}

\def\uu {\mathfrak{u}}

\def\cal{\mathcal}

\renewcommand{\(}{\begin{equation}}
\renewcommand{\)}{\end{equation}}
\newcommand{\bea}{\begin{eqnarray}}
\newcommand{\eea}{\end{eqnarray}}
\newcommand{\R}{{\mathbb R}}

\newcommand{\Z}{{\mathbb Z}}
\newcommand{\Q}{{\mathbb Q}}

\def\proof {{Proof.}\hspace{7pt}}
\def\endofproof {\hfill{$\Box$}\\}

\def\pull{\rfloor}

\def\pull{\mathbf{\rfloor}}
\def\hpull{\mathbf{\rfloor\!\!\rfloor}}

\begin{document}

\title{Twisted differential String and Fivebrane structures}

\author{Hisham Sati}
\address{Department of Mathematics, Yale University, New Haven, CT 06511
\newline
   Current address: Department of Mathematics,
University of Maryland, College Park, MD 20742}
\email{hsati@math.umd.edu}

\author{Urs Schreiber}
\address{Fachbereich Mathematik\\Universit\"at
Hamburg\\Bundesstra\ss e 55\\D--20146 Hamburg, Germany
\newline
Current address: 
Department of Mathematics,
Utrecht University,
3508 TA Utrecht,
The Netherlands}
\email{schreiber@math.uni-hamburg.de}

\author{Jim Stasheff}
\address{Department of Mathematics, University of Pennsylvania,
David Rittenhouse Lab,
Philadelphia, PA 19104-6395}
\email{jds@math.upenn.edu}

\begin{abstract}

In the background effective field theory of heterotic string theory, 
the Green-Schwarz anomaly cancellation mechanism plays a key role.
Here we reinterpret it and its magnetic dual version in terms of,
\emph{differential twisted String}- and \emph{differential twisted Fivebrane-structures}
that generalize the notion of $\mathrm{Spin}$-structures and $\mathrm{Spin}$-lifting gerbes 
and their differential refinement to smooth $\mathrm{Spin}$-connections.
We show that all these structures can be encoded in terms of \emph{nonabelian}
cohomology, \emph{twisted} nonabelian cohomology, and
\emph{differential} twisted nonabelian cohomology, extending the
differential generalized abelian cohomology as developed by Hopkins and Singer and
shown by Freed to formalize the global description of anomaly cancellation
problems in higher gauge theories arising in string theory.
We demonstrate that
the Green-Schwarz mechanism for the $H_3$-field,
as well as its magnetic dual version for the $H_7$-field define cocycles in 
{\it differential twisted nonabelian cohomology} that may be called, respectively, 
differential twisted
$\mathrm{Spin}(n)$-, $\mathrm{String}(n)$- and $\mathrm{Fivebrane}(n)$-
structures on target space, where the twist in each case is provided by the obstruction 
to lifting the classifying map of the gauge bundle through a higher connected cover of $U(n)$ or $O(n)$. We show that the twisted Bianchi identities in string theory can be captured by 
the (nonabelian) 
$L_\infty$-algebra valued differential form data provided by the
differential refinements of these twisted cocycles. 
\end{abstract}

\maketitle



\tableofcontents

\section{Introduction}

String theory and M-theory involve various higher gauge-fields, 
which are locally given by differential form fields of higher degree
and which are globally modeled by higher bundles with connection (higher
gerbes with connection, higher differential characters)
\cite{Freed}\cite{nactwist}. 
Some of these entities arise in terms of
lifts through various connected covers of orthogonal or unitary groups. 
For example, an orientation of a Riemannian manifold $M$ can be 
given by a lifting of the classifying map
${\it riem} :M\to B\mathrm{O}$ 
for the tangent or frame bundle of $M$ to a map ${\it or}:M\to B \mathrm{SO}$.
In turn, a Spin structure on $M$ can be given by a further lifting 
${\it sp}:M\to B \mathrm{Spin}$.
The existence of a Spin structure is an anomaly cancellation condition for
fermionic particles propagating on $M$.
The spaces $B\mathrm{O}$, $B \mathrm{SO}$ and $B \mathrm{Spin}$ are the first 
steps
 in the 
Whitehead tower of $B O$. The next step 
above $B \mathrm{Spin}$ is known as $B \mathrm{String}$, with $\mathrm{String}$
the topological group known as the \emph{String group}.

\vspace{3mm}
Originally Killingback \cite{Ki} defined a \emph{String structure} on $M$ as 
a lift of the transgressed map $L M \to B L \mathrm{Spin}$ on loop space 
through the Kac-Moody central extension $B {\hat L}\mathrm{Spin}(n)$.
The existence of such a lift cancels an anomaly of the 
heterotic superstring on $M$ in the case that the gauge bundle is trivial.
Later it was realized that this is captured down on $M$ by a lift of 
$sp : M \to B \mathrm{Spin}(n)$ to $str : M \to B \mathrm{String}(n)$ 
\cite{ST1}.
A further lift $fiv : M \to B \mathrm{Fivebrane}(n)$ 
through the next step in the Whitehead tower of $B \mathrm{O}(n)$ is 
similarly related to anomaly cancellation for the physical fivebrane on $M$
and accordingly the corresponding space is called $B \mathrm{Fivebrane}(n)$
\cite{SSS2}.

\begin{figure}[h]
$$
  \xymatrix{
    && B \mathrm{Fivebrane}(n)
    \ar[d]
    && \mbox{Fivebrane structure}
    \\
    && B \mathrm{String}(n)
    \ar[d]
    && \mbox{String structure}
    \\
    && B \mathrm{Spin}(n)
    \ar[d]
    && \mbox{Spin structure}
    \\
    && B \mathrm{SO}(n)
    \ar[d]
    && \mbox{Orientation}
    \\
    X \ar[rr]
    \ar@{-->}[urr]
    \ar@{-->}[uurr]
    \ar@{-->}[uuurr]
    \ar@{-->}[uuuurr]
    &&
    B \mathrm{O}(n)
    &&
    \mbox{Riemannian
    }
    \\
    \mbox{target space}
    &&
    \mbox{Whitehead tower of}\; B\mathrm{O}(n)
  }
$$
\caption{
  Topological structures generalizing $\mathrm{Spin}(n)$ structure.
  In application to effective background field theories appearing 
  in string theory, these bare structures are
  \emph{twisted} and moreover refined to \emph{differential} 
  structures.
}
\end{figure}
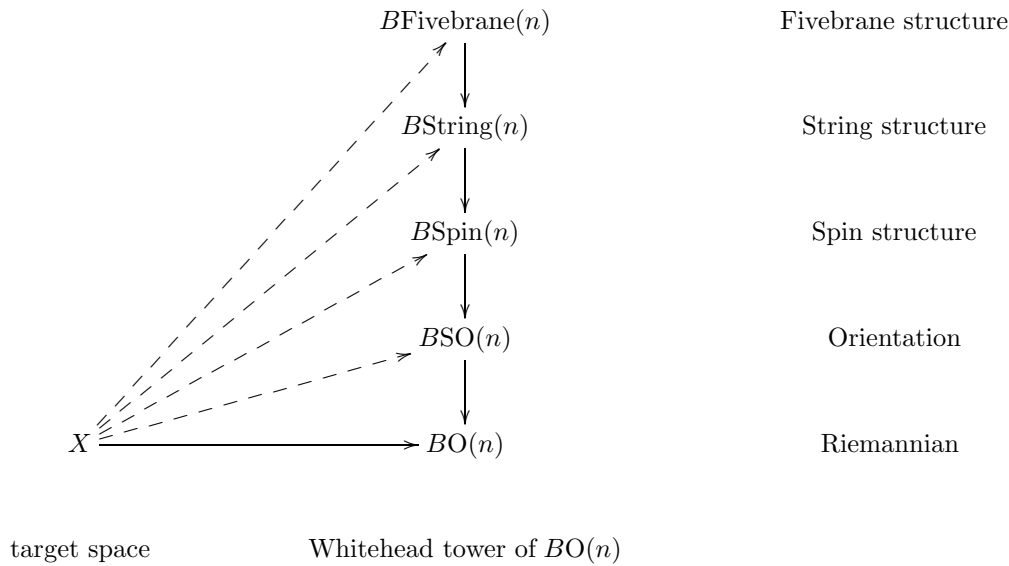

\vspace{3mm}
While anomalies are canceled by these lifts of maps of topological
spaces, the \emph{dynamics} of these systems is controlled by smooth refinements
of such maps. This is well understood for the first steps: the topological groups $O(n)$,
$SO(n)$ and $\mathrm{Spin}(n)$ naturally carry Lie group structures
and the differential refinement of $X \to B \mathrm{Spin}(n)$ is well
known to be given by a differential nonabelian $\mathrm{Spin}(n)$-cocycle,
namely a smooth $\mathrm{Spin}(n)$-principal bundle with connection. 

\vspace{3mm}
However, the higher connected 
topological groups $\mathrm{String}(n)$ and $\mathrm{Fivebrane}(n)$
cannot be finite-dimensional Lie groups and a smooth 
infinite-dimensional structure for $\mathrm{Fivebrane}$ is not known.
Moreover, even when such infinite-dimensional Lie group structures
on these higher covers exist 
(as was recently found for the $\mathrm{String}$-group 
\cite{NikolausSachseWockel}) by themselves they lead to the wrong 
smooth cohomological refinement (discussed in section 4.1 of
\cite{nactwist}).
However, $\mathrm{String}(n)$ does have a natural incarnation as a smooth \emph{2-group}
\cite{BCSS} \cite{H} \cite{BaezSchr} \cite{BBK} \cite{BSt}
\cite{SchommerPries}\cite{NikolausSachseWockel}, a higher categorical version of a Lie group
(see \cite{Lurie} for all general matters of higher category theory needed here
and \cite{nactwist} for smooth higher geometry).
Similarly, $\mathrm{Fivebrane}(n)$ does naturally exist as a 
\emph{smooth 6-group} (see section 4.1 of \cite{nactwist}). 
Generally, there are \emph{smooth $\infty$-group}-refinements
of all higher connected covers of Lie groups.
Being smooth, these
spaces have \emph{infinitesimal approximations}
by $L_\infty$-algebras in generalization of 
how any ordinary Lie group has a Lie algebra associated with it. 
Therefore, after passing to the smooth $\infty$-groupoid incarnation of the
objects in the Whitehead tower of $B{\rm O}$ there is a chance of obtaining differential refinements of
$\mathrm{String}(n)$- and $\mathrm{Fivebrane}(n)$-structures (and beyond) 
that are expressed
in terms of higher smooth bundles with smooth $L_\infty$-algebra-valued connection
forms on them, and indeed these structures naturally
exist \cite{SSS1}\cite{FSSI}.

\vspace{3mm}
The general refinement of cohomology classes to \emph{differential} cohomology
classes for the case of \emph{abelian} (Eilenberg-Steenrod-type) 
generalized cohomology theories has been discussed by
Hopkins and Singer \cite{HS} and shown by Freed \cite{Freed} 
to encode various differential (and twisted) structures in String theory.
However,
 the  cohomological structures that appear in the 
Freed-Witten \cite{FW} and
in the Green-Schwarz anomaly cancellation mechanism \cite{GS}, as well
as in the magnetic dual Green-Schwarz mechanism \cite{SSS2} 
themselves originate from --and are controlled by-- nonabelian
structures. These are the $\mathrm{O}(n)$-principal bundle
underlying the tangent bundle of spacetime and the $\mathrm{U}(n)$-principal
bundle underlying the gauge bundle on spacetime, as well as their
lifts to the higher connected structure groups: a map 
$X \to B\mathrm{Spin}(n)$ (smooth or not) gives a cocycle in 
\emph{nonabelian} cohomology, and so do its lifts
such as $X \to B {\rm String}(n)$ (and beyond). 

\vspace{3mm}
Therefore, here we provide applications for the theory of (twisted) differential 
\emph{nonabelian} cohomology, that builds on
\cite{BaezSchr} \cite{SWI}  \cite{SWII} \cite{SWII} \cite{SSS1}  
and is discussed in more detail in
\cite{nactwist}. We show that the Freed-Witten and the Green-Schwarz mechanisms,
as well as the magnetic dual Green-Schwarz mechanism, define 
differential twisted \emph{nonabelian}
 cocycles that may be interpreted as
 differential  twisted $\mathrm{Spin}^c$-,
$\mathrm{String}$- and $\mathrm{Fivebrane}$-structures, respectively.
Equivalent anomalies differ by a coboundary, so that they are given by 
cohomologous nontrivial cocycles. Equivalence classes of anomalies are
captured by the relevant cohomology.
We thus have a refinement of the treatment in \cite{SSS2} to the 
twisted, smooth and differential cases.  

\vspace{3mm}
In particular, the various abelian background fields
appearing in the theory, such as the B-field and the supergravity
3-form field, are unified into a natural coherent structure with the
nonabelian background fields -- the spin- and gauge-connections -- with which they interact. 
For instance, the relations between the abelian and the nonabelian differential
forms that govern the Green-Schwarz mechanism \cite{GS} are  realized here as 
a (twisted) Bianchi identity of a single nonabelian $L_\infty$-algebra valued
connection on a smooth twisted $\mathrm{String}(n)$-principal 2-bundle (cf. \cite{SSS1}).
 The explicit derivation of the twisted Bianchi
identites of $L_\infty$-algebra connections
corresponding to the Green-Schwarz mechanism and its magnetic dual
is in section \ref{Twisted differential String- and Fivebrane structures},
with more details in section 6.2 of \cite{FSSI}.

\begin{figure}[h]
$$
  \underbrace{
 \begin{array}{c}
  \xymatrix{
    {B}\mathrm{String}
     \ar[dd]
     &&
     X
     \ar@{-->}[ll]_{\hat g_{TX}}
     \ar@{-->}[rr]^{\hat g_E}
     \ar[ddll]|{g_{T X}}
     \ar[ddrr]|{g_E}
     \ar@/^1.5pc/[dd]|>>>>>>{\frac{1}{2}p_1(TX)}_{\ }="s"
     \ar@/_1.5pc/[dd]|>>>>>>{\frac{1}{6}c_2(E)}^{\ }="t"
     &&
    {B}U\langle 8\rangle
     \ar[dd]
     \\
     \\
{B}\mathrm{Spin}
     \ar[rr]^{\frac{1}{2}p_1}
     &&
{B}^{3}U(1)
     &&
{B}U\langle 6\rangle
     \ar[ll]_{\frac{1}{6}c_2}
     \ar@{<=>}_{H_3} "s"; "t"
  }
  \\
  \underbrace{\hspace{60pt}}_{\mbox{abelian cohomology}}
  \end{array}
  }_{\mbox{nonabelian cohomology}}
$$
\caption{
{\bf Abelian versus nonabelian cohomology.}
Since the groups $\mathrm{String}(n)$ as well as $\mathrm{Fivebrane}(n)$ are
shifted central extension of nonabelian groups, cohomology with coefficients
in these groups has abelian components but also components 
in nonabelian cohomology \cite{Toen2}. This appears as 
abelian cohomology twisted by nonabelian cocycles in a certain way.
  The Green-Schwarz mechanism implies that two classes in ordinary
  abelian cohomology, namely 
  in degree four differential integral cohomology, coincide. 
  But these classes are particularly obstruction classes to 
  $\mathrm{String}$-lifts in 
  nonabelian cohomology. The middle part of Figure 1 , 
  labeled ``abelian cohomology",
  identifies the cocycle representative in $H^4(X,\mathbb{Z})$ 
  and the coboundary between them,
  but does not specify where these cocycles come from. 
  The outer part of the diagram, labeled ``nonabelian
  cohomology" does specify the object whose class is the one identified by
  the middle part. 
  We can interpret this in ordinary homotopy theory, 
  where it describes topological obstruction theory,
  but we can also interpret this after differential refinement in
  the $\infty$-topos \cite{Lurie} of smooth $\infty$-groupoids \cite{nactwist}. 
  In any case the morphisms in the above diagram may be interpreted
  as \emph{cocycles}.
  The smooth and differential refinement we discuss in 
  section \ref{Twisted differential String- and Fivebrane structures}.
}
\end{figure}


\vspace{3mm}
\noindent {\bf Summary.}
In this paper we achieve the following goals:
\begin{enumerate}
\item 
generalize a Fiverbrane structure to the twisted case, and
similarly for a twisted String structure;
\item provide differential cohomology versions of these twisted structures;
\item provide a description of the Green-Schwarz anomaly cancellation and its dual
using these structures; 
\item describe the M-theory $C$-field and its dual in this context.
\end{enumerate}
The first two are purely mathematical results that are
of independent interest in developing higher (algebraic, geometric, 
topological, categorical) phenomena \cite{FSSI}. The third and fourth are
applications to (heterotic) string theory and to 
 M-theory, respectively.  We hope this will add to a better understanding 
of structures appearing in these theories, which in turn is hoped to result in 
identification of yet more rich mathematical structures within them
(see \cite{S1} \cite{S2} \cite{S3} for concrete examples).
From 
a mathematical point of view, they serve as interesting concrete
examples of the formalism we have developed in the first two points.

\vspace{3mm}
Section \ref{TwistedTopologicalStructuresInStringTheory} 
discusses
    how the anomaly cancellation mechanisms in
    String-theory can be understood topologically in terms of  twisted higher structures
    given by twisted nonabelian topological cocycles.
 The physical examples of most relevance here arise in various
anomaly cancellations in string theory. The Freed-Witten condition
\cite{FW} in type IIA string theory says that the third integral
Stiefel-Whitney class $W_3$ of a D-brane $Q$ has to be \emph{trivial
relative to the Neveu-Schwarz field $H_3|_Q$} restricted to the
D-brane, in that the two classes agree:
$W_3 = [H_3|_Q]$.

\vspace{3mm} 
Higher versions of this example --within our point of view-- are the Green-Schwarz mechanism
and its magnetic dual version.
 Recall the notion of String structures, e.g. from 
  \cite{SSS2}, as  maps from a space $X$ to $B{\rm String}(n)$, the 3-connected cover of $B{\rm Spin}(n)$.
 In \cite{Wang} the notion of twist for a String
structure was considered: a space $X$ can have a twisted
String structure without having a String structure, i.e.
the fractional Pontrjagin class $\frac{1}{2}p_1(TX)$ 
of the tangent bundle
can be nonzero while the modified class is $\frac{1}{2}p_1(TX) + [\beta]=0$,
where $\beta : X \to K(\Z, 4)$ is a fixed \emph{twist} for the
String structure. The Green-Schwarz mechanism in string theory may 
be understood as defining a twisted $\mathrm{String}$ structure
on target space, with twist given in terms of  a classifying map for the gauge bundle.

\vspace{3mm}
Since a $\mathrm{String}$ structure is refined by a $\mathrm{Fivebrane}$ structure
in analogy to how  a $\mathrm{String}$ structure itself refines a
$\mathrm{Spin}$ structure, it is
natural to consider twists of $\mathrm{Fivebrane}$ structures in the
above sense. In this paper we give a definition of {\it twisted
Fivebrane structures} and show that the dual Green-Schwarz mechanism 
in heterotic String theory, reviewed and formalized in detail in \cite{SSS2}, 
provides an example. 
 Hence, variations on the twisted Fivebrane condition do
  in fact appear in string theory and in M-theory and correspond, as we will see, to anomaly cancellation conditions
 for the heterotic fivebrane \cite{DDP} \cite{LT} and for the M-fivebrane
\cite{Effective} \cite{Among} \cite{FHMM} \cite{DFM}, respectively.
We
discuss these two cases in section \ref{HeteroticStringAndTwistedStringStructures} and
section \ref{NS5 and twisted Fivebrane Structures} in terms of topological cocycles
(maps to the appropriate classifying spaces)
and describe  their
differential refinements in section \ref{G4 twist} and \ref{G8 twist}.

Some time has passed between the original inception of the
discussion presented in this article and the present form. The
differential structures discussed here have motivated us with
collaborators to further expand the general theory of higher smooth
stacks and the formulation of differential cohomology in terms of
these. In section 3 we survey this context and refer to various
technical results obtained meanwhile.

\section{Twisted topological structures in String theory}
\label{TwistedTopologicalStructuresInStringTheory}

We discuss here cohomological conditions arising from
anomaly cancellation in String theory, for various $\sigma$-models. 
In each case, we introduce a corresponding notion of  
topological \emph{twisted structures} in terms of which we interpret the 
corresponding anomaly 
cancellation condition.
This prepares the
ground for the material in section \ref{Twisted differential String- and Fivebrane structures},
where the differential refinement of these twisted structures is considered leading to the derivation 
of the \emph{differential} anomaly-free 
field configurations.

\medskip

The physics of all the cases we consider involves a 
manifold $X$ --the \emph{target space}-- or a submanifold $Q \hookrightarrow X$ 
 thereof --a \emph{brane}--
equipped with 
\begin{itemize}
\item two principal bundles with their canonically associated vector bundles:
\begin{itemize}
\item a $\mathrm{Spin}$-principal bundle underlying the tangent bundle
$T X$ (and we will write $T X$ also to denote that $\mathrm{Spin}$-principal bundle), 
\item and a complex vector bundle $E \to X$
-- the ``gauge bundle'' --  associated to a $SU(n)$-principal bundle or
to an $E_8$-principal bundle with respect to a unitary representation of $E_8$;
\end{itemize}
  \item  an $n$-gerbe / circle $(n+1)$-bundle with class $ H^{n+2}(X,\mathbb{Z})$
  representing  the higher background gauge field and denoted $[H_i]$ or $[G_i]$ 
	or similar in the following.
\end{itemize}
All these structures are equipped with a suitable notion of 
\emph{connetions}, locally given by some differential-form data.
The connection on the $\mathrm{Spin}$-bundle encodes the field of 
gravity, that on the gauge bundle a Yang-Mills field and that 
on the $n$-gerbe a higher analog of the electromagnetic field.

The $\sigma$-model quantum field theory of a super-brane 
propagating in such 
a background (for instance the superstring, or the super 5-brane)
has an effective action functional on its bosonic 
worldvolume fields that takes values, in general, in the fibers of 
the Pfaffian line bundle of a worldvolume Dirac operator,
tensored with a line bundle that encodes the electric and magnetic charges
of the higher gauge field. Only if this 
tensor product \emph{anomaly line bundle} is trivializable is the
effective bosonic action a 
well-defined starting point for
quantization of the $\sigma$-model. Therefore, the Chern class of this
line bundle over the bosonic configuration space is called the 
 \emph{global anomaly} of the system.
Conditions on the background gauge fields that ensure that this class
vanishes are called \emph{global anomaly cancellation conditions}.
These turn out to be conditions on cohomology classes that are characteristic 
of the above background fields. 
This is what we discuss in this section.
Moreover, the anomaly line bundle is canonicaly equipped with a 
\emph{connection},
induced from the connections of the background gauge
fields, hence induced from their \emph{differential cohomology} data. 
The curvature 2-form of this connection over
the bosonic configuration space is called the \emph{local anomaly}
of the $\sigma$-model. Conditions on the differential data of the 
background gauge field that canonically induce a trivialization of this
2-fom are called \emph{local anomaly cancellation conditions}. We 
consider these below in section \ref{Twisted differential String- and Fivebrane structures}.

\medskip

The phenomenon of anomaly line bundles of $\sigma$-models 
induced from background field differential cohomology is 
classical in the physics literature, albeit in broad terms. 
A clear exposition is in \cite{Freed}. Only recently 
has the special case of the heterotic string $\sigma$-model for trivial 
background gauge bundle been made fully precise in 
\cite{Bunke}, using a certain model \cite{Waldorf} 
for the differential string structures
that we discuss in section 
\ref{Twisted differential String- and Fivebrane structures}.

\subsection{Twisted $\mathrm{Spin}^c$ structures and physical applications} 
\label{Type II superstring and twisted Spinc structures}

As a preparation for the twisted String- and Fivebrane structures discussed in the 
following sections,  we consider twisted $\mathrm{Spin}^c$-structures and 
their role in anomaly cancellation for the open type II string.

\subsubsection{Type II superstring on D-branes}

The open type II string propagating on a Spin-manifold $X$ in the presence of 
a) a background $B$-field with class $[H_3] \in H^3(X,\mathbb{Z})$
and b) with endpoints fixed on a D-brane
given by an oriented submanifold $Q \hookrightarrow X$,
has a global worldsheet anomaly that vanishes
if \cite{FW} the condition
\(
  W_3(Q) + [H_3]|_Q = 0 \;\; \in H^3(Q;\Z)
  \;,
  \label{Freed Witten condition}
\)
holds.
Here $W_3(Q)$ is the third integral Stiefel-Whitney class of the tangent bundle $TQ$
of the brane and $[H_3]_Q$ denotes the restriction of $[H_3]$ to $Q$.
Sufficiency of the condition is discussed in \cite{ES2}.

Notice that $W_3(Q)$ is the obstruction to lifting the orientation structure on 
$Q$ to a $\mathrm{Spin}^c$ structure. We can 
formulate this in terms of homotopy theory as follows. There is a \emph{homotopy pullback diagram}
\(
  \raisebox{20pt}{
  \xymatrix{
    B \mathrm{Spin}^c \ar[r]_<{\hpull} \ar[d] & {*} \ar[d]
    \\
    B \mathrm{SO} \ar[r]^{\hspace{-2mm}W_3} & B^2 U(1)
  }
  }
  \label{SpincAsHomotopyFiber}
\)
of topological spaces, where $B \mathrm{SO}$ is the classifying space of the (stable)
special orthogonal group, $B^2 U(1) \simeq K(\mathbb{Z},3)$ is 
the Eilenberg-MacLane space that classifies degree-3
integral cohomology, and the continuous map denoted $W_3$ is a representative
of the universal class $W_3$ under this classification. 

\vspace{3mm}
Homotopy pullbacks and their properties play a central role in all of our 
discussion here, first in the category of topological spaces, and then
later, in section \ref{Twisted differential String- and Fivebrane structures},
in categories of higher stacks.  
The reader unfamiliar with the basics of abstract homotopy theory
might consult the review in section A.2 of \cite{Lurie}. 
A basic fact is that a \emph{homotopy fiber} 
of connected spaces as in (\ref{SpincAsHomotopyFiber}),
where the right vertical morphism is the point inclusion, may
be computed, up to weak homotopy equivalence, as an ordinary pullback 
(an ordinary fiber product)
after replacing the point inclusion by the path fibration.
Specifically, for $X$ any pointed topological space, write $P X$ for the 
space of continuous paths in $X$ that end at the base point $x \in X$, and
write $P X \to X$ for the projection to the other endpoint of the path. Then the 
ordinary pullback
\(
  \xymatrix{
    \Omega_x X \ar[r]_<{\pull}\ar[d] & P X \ar[d]
	\\
	{*} \ar[r]^x & X
  }
\)
is a model for the loop space of $X$. Moreover, for $\phi : Y \to X$ any 
map, the ordinary pullback
\(
  \xymatrix{
    \mathrm{hofib}(\phi) \ar[r]_<{\pull}\ar[d] & P X \ar[d]
	\\
	Y \ar[r]^\phi & X
  }
\)
is a model for the homotopy fiber 
$\mathrm{hofib}(\phi)$
of $\phi$. 

\vspace{3mm}
The homotopy pullback (\ref{SpincAsHomotopyFiber}) exhibits the 
classifying space of the group $\mathrm{Spin}^c$ as the homotopy fiber of $W_3$.
The universal property of the homotopy pullback says that the space of continuous maps
$Q \to B \mathrm{Spin}^c$ is the same (is weak homotopy equivalent to) the space of
maps $o_Q : Q \to B \mathrm{SO}$ that are equipped with a homotopy from the composite
$\xymatrix{
  Q \ar[r]^{\hspace{-4mm}o_Q} &   B \mathrm{SO} \ar[r]^{\hspace{-2mm}W_3} & B^3 U(1)
}$ to the trivial cocycle 
$Q \to {*} \to B^3 U(1)$. In other words,  for every choice of homotopy
filling of the outer diagram
\(
  \xymatrix{
    Q 
    \ar@{-->}[dr]
    \ar@/^1pc/[drr]
    \ar@/_1pc/[ddr]_{o_Q}
    \\
    & B \mathrm{Spin}^c \ar[r] \ar[d] & {*} \ar[d]
    \\
    & B \mathrm{SO} \ar[r]^{\hspace{-2mm}W_3} & B^2 U(1)
  }
\)
there is a contractible space of 
choices for the dashed arrow such that everything commutes up to
homotopy. Since a choice of map $o_Q : Q \to B \mathrm{SO}$ is
an \emph{orienation structure} on $Q$, and a choice of map $Q \to B \mathrm{Spin}^c$ is
a \emph{$\mathrm{Spin}^c$ structure}, this implies that $W_3(o_Q)$ is the obstruction
to the existence of a $\mathrm{Spin}^c$structure on $Q$ (equipped with  $o_Q$).
Moreover, since $Q$ is a manifold (hence a CW-complex), 
the functor $\mathrm{Maps}(Q,-)$ that forms mapping spaces out of $Q$ preserves
homotopy pullbacks. Since $\mathrm{Maps}(Q, B \mathrm{SO})$
is the \emph{space} of orientation structures, we can refine the discussion so far by 
noticing that the \emph{space of $\mathrm{Spin}^c$ structures on $Q$}, 
denoted $\mathrm{Maps}(Q, B \mathrm{Spin}^c)$, is itself the homotopy pullback in the diagram
\(
  \raisebox{20pt}{
  \xymatrix{
    \mathrm{Maps}(Q,B \mathrm{Spin}^c) \ar[rr]_<{\hpull}
	\ar[d] && {*} \ar[d]
    \\
    \mathrm{Maps}(Q,B \mathrm{SO}) \ar[rr]^{\hspace{-2mm}\mathrm{Maps}(Q,W_3)} && 
     \mathrm{Maps}(Q,B^2 U(1))
  }
  }
  \label{SpaceOfSpincStructures}
  \,.
\)
A variant of this characterization will be crucial for the definition of 
(spaces of) \emph{twisted} such structures below.

These kinds of arguments,  although elementary in homotopy theory, are of importance
for the interpretation of anomaly cancellation conditions that we consider here. Variants
of these arguments (first for other topological structures, then with twists, 
then refined to smooth
and differential structures) will appear over and over again in our discussion.

In the case that the class of the $B$-field vanishes on the D-brane, 
$[H_3]|_Q = 0$, hence that its representative $H_3 : Q \to K(\mathbb{Z},3)$
factors through the point up to homotopy,  
condition (\ref{Freed Witten condition}) states that the oriented D-brane $Q$ must 
admit a  $\mathrm{Spin}{}^c$ structure, namely a choice of null-homotopy $\eta$ in
\footnote{Beware that there are homotopies filling \emph{all} our diagrams, but only in some
cases, such as here, do we want to make them explicit and give them a name.}
\(
  \raisebox{20pt}{
  \xymatrix{
     Q
     \ar[rr]^{o_Q}_>>>{\ }="s"
     \ar[drr]_{{H_3}|_Q \simeq {*}}^{\ }="t"
	 \ar[d]
     &&
     B \mathrm{SO}
     \ar[d]^{W_3}
     \\
	 X
	 \ar[rr]_{H_3}
     &&
     K(\mathbb{Z},3)
     \ar@{=>}^\eta "s"; "t"
  }
  }
  \,.
\)
If, more generally, $[H_3]|_Q$ does not necessarily vanish, 
then condition (\ref{Freed Witten condition})
still is equivalent to the existence of a homotopy $\eta$ in a diagram of the above form:
\(
  \raisebox{20pt}{
  \xymatrix{
     Q
     \ar[rr]^{o_Q}_>>>{\ }="s"
     \ar[drr]_{{H_3}|_Q }^{\ }="t"
	 \ar[d]
     &&
     B \mathrm{SO}
     \ar[d]^{W_3}
     \\
	 X
	 \ar[rr]_{H_3}
     &&
     K(\mathbb{Z},3)
     \ar@{=>}^\eta "s"; "t"
  }
  }
  \label{HomotopyOfTwistedSpinCStructure}
  \,.
\)
We may think of this as saying that $\eta$ still ``trivializes'' $W_3(o_Q)$, but not
with respect to the canonical trivial cocycle, but with respect to the given reference background 
cocycle ${H_3}|_Q$ of the $B$-field. Accordingly, following \cite{Wang}, 
we may say that such an $\eta$ exhibits not a $\mathrm{Spin}^c$-structure on $Q$, but
an \emph{$[H_3]_Q$-twisted $\mathrm{Spin}^c$ structure}.

For this notion to be useful, we need to say what an equivalence or homotopy
between two twisted $\mathrm{Spin}^c$ structures is, what a homotopy between such
homotopies is, etc., and hence what the \emph{space} of twisted 
$\mathrm{Spin}^c$ structures is. 
However, by generalization of (\ref{SpaceOfSpincStructures}) 
we naturally do have such a space.
\begin{definition}
\label{SpaceOfSPinCStructures}
For $X$ a manifold and $[c] \in H^3(X, \mathbb{Z})$ a degree-3 cohomology
class, we say that the space
$W_3\mathrm{Struc}(Q)_{[c]}$ defined as the homotopy pullback
\(
 \raisebox{30pt}{
  \xymatrix{
    W_3 \mathrm{Struc}(Q)_{[H_3]|_Q}
	\ar[rr]_<{\hpull}
	\ar[d]
	&&
	{*}
	\ar[d]^{c}
	\\
	\mathrm{Maps}(Q, B \mathrm{SO})
	\ar[rr]^{\hspace{-3mm}\mathrm{Maps}(Q,W_3)}
	&&
	\mathrm{Maps}(Q, B^2 U(1))
  }
  }
  \label{space of twisted Spinc structures}
  \,,
\)
is the \emph{space of $[c]$-twisted $\mathrm{Spin}^c$ structures} on $X$,
where the right vertical morphism picks any representative
$c : X \to B^2 U(1) \simeq K(\mathbb{Z},3)$ of $[c]$.
\end{definition}

In terms of this notion, the anomaly cancellation condition
(\ref{Freed Witten condition}) is now read as being
a requirement of  
\emph{existence of structure}:
\begin{observation}
On an oriented manifold $Q$, condition (\ref{Freed Witten condition}) 
implies the existence of 
\emph{$[H_3]|_Q$-twisted $W_3$ structure}, 
provided by a lift of the orientation structure $o_Q$
on $T Q$ through the left vertical morphism in 
def. \ref{SpaceOfSPinCStructures}.
\end{observation}
This makes good sense, because that extra structure is the extra 
structure of the background field  of the 
$\sigma$-model background, subjected to the
condition of anomaly freedom. We will see this in more detail in the
following examples, and then again
in section \ref{Twisted differential String- and Fivebrane structures}.


\subsection{Twisted String structures and physical applications}

We discuss \emph{twisted String structures} and their role in anomaly
cancellation of
\begin{enumerate}
  \item the heterotic string;
  \item M-theory in the bulk;
  \item the boundary in M-theory.
\end{enumerate}

\subsubsection{The heterotic string}
\label{HeteroticStringAndTwistedStringStructures}

The heterotic/type I string, propagating on a $\mathrm{Spin}$-manifold
$X$ and coupled to a gauge field given by a Hermitean complex vector bundle $E \to X$,
has a global anomaly that vanishes if the 
\emph{Green-Schwarz anomaly cancellation condition} \cite{GS}
\(
  \frac{1}{2}p_1(TX) - {\rm ch}_2(E)=0 \;\;\in H^4(X;\Z)
  \label{Green Schwarz anomaly}
\)
holds. Here $\mathrm{ch}_2(E)$ is the second Chern character of $E$
which reduces to the second Chern class $c_2(E)$ in the cases we consider, 
and $\frac{1}{2}p_1$ is given by the following
classical fact (see \cite{Bott}).
\begin{fact}
  The first Pontrjagin class $p_1 \in H^4(B \mathrm{SO}, \mathbb{Z})$
  becomes divisible precisely by 2 when pulled back to $B \mathrm{Spin}$.
  The corresponding preimage under multiplication by two, denoted
  $\frac{1}{2}p_1$,
  is a generator of the group $H^4(B \mathrm{Spin}, \mathbb{Z}) \cong \mathbb{Z}$. 
  \label{FirstFractionalPontrjagin}
\end{fact}

As in the Spin${}^c$ case discussed above, 
this means that at the level of cocycles a certain homotopy
exists. Here it is this homotopy which is the representative of the 
$B$-field to which the string couples.
In detail, write $\frac{1}{2}p_1 : B \mathrm{Spin} \to B^3 U(1)$ for 
a representative of the
universal first fractional Pontrjagin class, and similarly
${c}_2 : B \mathrm{SU} \to B^3 U(1)$ for a representative of the
universal second Chern class, where now
$B^3 U(1) \simeq K(\mathbb{Z},4)$ is equivalent to the 
Eilenberg-MacLane space that classifies degree-4 integral cohomology. 
Then if $T X : X \to B \mathrm{Spin}$ is a classifying
map of the $\mathrm{Spin}$-bundle and $E : X \to B \mathrm{SU}$ is one of the
gauge bundles, the above anomaly cancellation condition says that there
is a homotopy, denoted $H_3$, in the diagram
\(
  \raisebox{20pt}{
  \xymatrix{
    X \ar[r]^{E}_>{\ }="s" \ar[d]_{T X} & B \mathrm{SU} \ar[d]^{c_2}
	\\
	B \mathrm{Spin}
	\ar[r]_{\frac{1}{2}p_1}^<{\ }="t"
	&
	B^3 U(1)
  \ar@{=>}_{H_3} "s"; "t"
  } }
  \label{homotopy in twisted string structure}
  \,.
\)
Notice that if both $\frac{1}{2}p_1(T X)$ as well as $c_2(E)$
happen to be trivial, such a homotopy is equivalently a map 
$H_3 : X \to \Omega B^3 U(1) \simeq B^2 U(1)$. So in this special case the 
B-field in the background of the heterotic string is a
$U(1)$-gerbe, or a circle 2-bundle, as in the previous case of the type II string
in section \ref{Type II superstring and twisted Spinc structures}.
Generally, the homotopy $H_3$ in the above diagram exhibits
the B-field as a \emph{twisted} gerbe, whose twist is the
difference class $\frac{1}{2}p_1(TX)- c_2(E)$.
This is essentially the perspective adopted in \cite{Freed} (in a somewhat different language).

For the general discussion of interest here it is useful to slightly
shift the perspective on the twist. To that end, 
first consider, by analogy with 
(\ref{SpincAsHomotopyFiber}), the following definition.
\begin{definition}
The \emph{String group} is the loop space of the
homotopy fiber $B \mathrm{String}$ of the universal class
$\frac{1}{2}p_1$ from Proposition \ref{FirstFractionalPontrjagin}:
$$
  \raisebox{20pt}{
  \xymatrix{
    B \mathrm{String} 
	\ar[r]_<{\hpull}
	\ar[d]
	&
	{*}
	\ar[d]
	\\
	B \mathrm{Spin}
	\ar[r]^{\frac{1}{2}p_1}
	&
	B^3 U(1)
  }}
  \,.
$$
 \label{StringGroup}
\end{definition}
For more details on this and a collection of
references see section 4.1 of \cite{nactwist}.
Accordingly, one says that a \emph{String structure} on 
the Spin bundle $T X : X \to B \mathrm{Spin}$ is a homotopy 
filling the outer square of the diagram
\(
  \raisebox{30pt}{
  \xymatrix{
    X
	\ar@{-->}[dr]
	\ar@/^1pc/[drr]
	\ar@/_1pc/[ddr]_{T X}
    \\
    &
	B \mathrm{String} 
	\ar[r]
	\ar[d]
	&
	{*}
	\ar[d]
	\\
	& B \mathrm{Spin}
	\ar[r]^{\frac{1}{2}p_1}
	&
	B^3 U(1)
  }}
  \,,
\)
or, equivalently-- by  the universal property of homotopy pullbacks--
a choice of dashed morphism filling the interior of this square, as indicated.
\footnote{
Originally, this condition was considered in its weaker incarnation after transgression
to loop space \cite{Ki}.}
Therefore, now by analogy with (\ref{HomotopyOfTwistedSpinCStructure}),
we say that a $c_2(E)$-\emph{twisted String structure}
is a choice of homotopy $H_3$ filling the diagram
(\ref{homotopy in twisted string structure}).
This notion of twisted String structures was originally
suggested in \cite{Wang}. For it to be useful, we need to 
say what homotopies of twisted $\mathrm{String}$ structures
are, homotopies between these, etc. Hence we need to say
what the \emph{space} of twisted $\mathrm{String}$ structures
is. This is what the following definition provides, analogously
to  Definition \ref{SpaceOfSPinCStructures}.
\begin{definition}
  \label{TwistedStringStructures}
  For $X$ a manifold,  
  and for $[c] \in H^4(X,\mathbb{Z})$ a degree-4 cohomology class,
  we say that the space of \emph{$c$-twisted $\mathrm{String}$ structures}
  on $X$ is the homotopy pullback $\frac{1}{2}p_1 \mathrm{Struc}_{[c]}(X)$ in
  the diagram
  \(
    \raisebox{20pt}{
    \xymatrix{
	  \frac{1}{2}p_1 \mathrm{Struc}_{[c]}(X)
	  \ar[rr]_<{\hpull}
	  \ar[d]
	  &&
	  {*}
	  \ar[d]^{c}
	  \\
	  \mathrm{Maps}(X, B \mathrm{Spin})
	  \ar[rr]^{\hspace{-2mm}\mathrm{Maps}(X,\frac{1}{2}p_1)}
      &&
      \mathrm{Maps}(X, B^3 U(1))	  
	}
	}
	\,,
  \)
  where the right vertical morphism picks a representative $c$ of $[c]$.
  \label{SpaceOfTwistedStringStructures}
\end{definition}

In terms of this then, we find
\begin{observation} 
 The anomaly cancellation condition 
(\ref{Green Schwarz anomaly}) is, for a fixed
gauge bundle $E$, precisely the condition that 
ensures a lift of the given $\mathrm{Spin}$ structure to a
$c_2(E)$-twisted $\mathrm{String}$ structure on $X$,
through the left vertical morphism of def. \ref{SpaceOfTwistedStringStructures}.
\end{observation}
Of course the full background field content involves more than just 
this topological data; it also consists of local differential form data,
such as a 1-form connection on the bundles $E$ and on $T X$ 
and a connection 2-form on the 2-bundle whose curvature is $H_3$. Below, in 
section \ref{Twisted differential String- and Fivebrane structures}, we
identify this \emph{differential} anomaly-free field content
with a \emph{differential} twisted $\mathrm{String}$ structure.

\vspace{3mm}
The 10-dimensional string theory backgrounds discussed so far
are supposed to be part of a bigger picture, which schematically looks 
as follows
 \(
  \raisebox{30pt}{
  \xymatrix{
  \fbox{UV-complete theory}
  \ar[rrrr]^{\mbox{low energy approximation}}
  &&&&
  \fbox{effective field theory}
  \\
  \mbox{M-theory}
  \ar[rrrr]^{\ell_1}
   \ar[d]^{{\partial}_2}
  &&&&
 \mbox{$D=11$ Supergravity}
 \ar[d]^{{\partial}_1}
 \\
 \mbox{Heterotic~String}
 \ar[rrrr]^{\ell_2}
 &&&&
  \mbox{$D=10$ Supergravity SYM}\;.
  }
  }
  \label{UVAndBoundaryCompletion}
\)
Here the bottom and right structures and maps are fairly well-defined.
The right vertical map denotes the restriction of 11-dimensional supergravity
on a manifold with boundary to 10-dimensional heterotic supergravity
on the boundary. The bottom left morphism indicates that heterotic
supergravity is the effective low energy field theory whose UV-completion
is heterotic string theory.
The top left entry denotes M-theory, which is what meaningfully completes this 
schematic diagram. Accordingly, the  2-brane and 5-brane
solutions of 11-dimensional supergravity have incarnations as fundamental
objects in M-theory, the \emph{M2-brane} and the \emph{M5-brane}.

We discuss the topological anomalies for the M2-brane, first in the bulk theory,
then after restriction to the boundary, where it becomes the heterotic string
of section \ref{HeteroticStringAndTwistedStringStructures}.
At the level of fields, this is the $C$-field in the bulk becoming the 
$B$-field on the boundary.

\subsubsection{M-theory in the bulk}

The bosonic  field content of 11-dimensional supergravity 
on a Spin 11-manifold $Y$
consists of a
$\mathrm{Spin}$-bundle $T Y$ (with connection), 
as well as of the \emph{$C$-field}, which has underlying
it a 2-gerbe -- or \emph{circle 3-bundle} -- with class $[G_4] \in H^4(Y, \mathbb{Z})$.
The M2-brane that couples to these background fields has an anomaly that
vanishes if \cite{Flux} 
\(
  2[G_4] = \frac{1}{2}p_1(T Y) - 2 a(E)
  \;\;
  \in H^4(Y, \mathbb{Z})
  \label{WittenQuantizationCondition}
  \,,
\)
where $E \to Y$ is an auxiliary $E_8$-principal bundle, whose class $a(E)$ is
defined by this condition. This is also called 
the \emph{quantization condition for the $C$-field}.

In the absence of smooth or differential structure,
since $E_8$ is 15-coskeletal, 
one could therefore replace the $E_8$-bundle here by a U(1)-2-gerbe, hence by
a $B^2 U(1)$-principal bundle, and
replace condition (\ref{WittenQuantizationCondition}) by
\(
  2[G_4] = \frac{1}{2}p_1(T Y) - 2 \mathrm{DD}_{2}
  \,,
\label{eq DD2}
\)
where $\mathrm{DD}_2$ is the canonical 4-class of this 2-gerbe 
(the ``second Dixmier-Douady class'', hence the notation).
While topologically this condition is equivalent, over an 11-dimensional $X$,
to (\ref{WittenQuantizationCondition}), the spaces of solutions of 
smooth refinements of these
two conditions will differ, because the space of smooth gauge transformations
between $E_8$ bundles is quite different from that of smooth gauge transformations
between circle 2-bundles.
In the Ho{\v r}ava-Witten reduction \cite{HW} of the 11-dimensional theory 
down to the  heterotic string in 10 dimensions, this difference is supposed to be 
relevant, since the heterotic string in 10 dimensions sees the
smooth $E_8$-bundle with connection.

In either case, we can understand the situation as a refinement of that described by
(twisted) $\mathrm{String}$-structures, as above,
via a higher analogue of the passage from
$\mathrm{Spin}$-structures to $\mathrm{Spin}^{c}$-structures, as in section
\ref{Type II superstring and twisted Spinc structures}. 
To that end,
notice the following
fact, which provides an alternative perspective on 
(\ref{SpincAsHomotopyFiber}), and which uses the point of view on twisted
structures advocated recently in \cite{S2} \cite{S3}.

\begin{proposition}
  \label{SpinCByHomotopyPullback}
  The classifying space $B \mathrm{Spin}^c$ is the homotopy fiber product of 
  a representative of the 
  universal second Stiefel-Whitney class $w_2 \in H^2(B \mathrm{SO}, \mathbb{Z}_2)$ 
  with a representative of the mod 2-reduction of the
  universal first Chern class $c_1 \in H^2(B U(1), \mathbb{Z})$
  $$
    \raisebox{20pt}{
    \xymatrix{
      B \mathrm{Spin}^c 
	  \ar[r]_<{\hpull}
	  \ar[d] 
	  & B U(1) \ar[d]^{c_1~ \mathrm{mod} ~2}
      \\
      B \mathrm{SO} \ar[r]^{w_2} & B^2 \mathbb{Z}_2
    }
	}
    \,.
  $$
\end{proposition}
This remains true after refinement to smooth and differential structures, as
we discuss below in \ref{Twisted differential String- and Fivebrane structures}.
Due to the universal property of the homotopy
pullback, this says, in particular, 
that a lift from an orientation structure to a $\mathrm{Spin}^c$-structure
is a cancelling by a Chern class of the class obstructing a $\mathrm{Spin}$-structure. 
In this way, lifts from orientation structures to  
$\mathrm{Spin}^c$-structures are analogous to the divisibility condition
(\ref{WittenQuantizationCondition}), since in both cases the obstruction to 
a further lift through the Whitehead tower of the orthogonal group is
absorbed by a universal ``unitary'' class.

In higher analogy to this situation, and in the spirit of \cite{S2} \cite{S3},
 we therefore have the following definition.
\begin{definition}
  \label{StringC}
 For $G$ some topological group, and $\alpha : B G \to K(\mathbb{Z},4)$
 a universal 4-class, we say that $\mathrm{String}^\alpha$ is the 
loop group of the homotopy pullback
\footnote{We are using the notation String${}^\alpha$ to distinguish from the 
notion of String${}^c$ structure (related to Spin${}^c$ in the same way that String is related to Spin)
studied in \cite{S3}. The latter is a special case of the 
former, as explained in op. cit.}
 $$
   \xymatrix{
     B \mathrm{String}^\alpha \ar[r]_<{\hpull} 
	 \ar[d] 
	 & 
	 B G \ar[d]^{\alpha}
     \\
     B \mathrm{Spin} \ar[r]^{\frac{1}{2}p_1} & B^3 U(1)
   }
 $$
of $\alpha$ along a representative of the first fractional Pontrjagin class
$\frac{1}{2}p_1 \in H^4(B \mathrm{Spin}, \mathbb{Z})$.
 \label{Stringc}	
\end{definition}

Of relevance for the present purpose are the following cases.
Notice that on the space $B^2 U(1)\simeq K(\mathbb{Z},3)$, which classifies circle 2-bundles / $U(1)$-bundle gerbes, 
the canonical 3-class is traditionally called the \emph{Dixmier-Douady class}, denoted $\mathrm{DD}$. Accordingly,
 it makes sense to speak of the canonical 4-class on $B^3 U(1) \simeq K(\mathbb{Z},4)$, which classifies 
 circle 3-bundle / $U(1)$-bundle 2-gerbes
as the \emph{second Dixmier-Douady class} $\mathrm{DD}_2$. 
We are interested in the case when the first Spin characteristic class $\lambda=\frac{1}{2}p_1$ is divisible by 2,
so that we can `divide by 2' in equation \eqref{eq DD2}.
For $\alpha = \mathrm{DD}_2$ we have
that a $\mathrm{Spin}$-structure lifts to a $\mathrm{String}^{2\mathrm{DD}_2}$-structure
precisely if $\frac{1}{2}p_1$ is further divisible by 2. Equivalently,
a $\mathrm{Spin}$-principal bundle $P$ lifts to a 
$\mathrm{String}^{2 \mathrm{DD}_2}$-principal bundle precisely if the 
corresponding \emph{Chern-Simons 2-gerbe} $\frac{1}{2}p_1(P)$ is 
twice some other 2-gerbe.
Equivalently this says that, the
fourth Stiefel-Whitney class vanishes, $w_4=0$,
because, as explained in \cite{S1} \cite{S2}, $w_4$ 
is the mod 2 reduction of the 
integral first Spin characteristic class $\lambda=\frac{1}{2}p_1$.
In summary, we have for a given $\mathrm{Spin}$-structure
$Y \to B \mathrm{Spin}$ the following diagram
\(
  \xymatrix{
  Y
  \ar@/^1.3pc/[rrrr]^{\frac{1}{4}p_1(X)}
  \ar@{-->}[rr]
\ar[rrd]
  &&
   {B} \mathrm{String}^{2\mathrm{DD}_2}
   \ar[rr]_<{\hpull}|{\frac{1}{4}p_1}
   \ar[d]
   &&
   K(\Z, 4)
   \ar[d]^{\times 2}
   \ar[drr]^0
  &&
   \\
   &&
   (BO)\langle 4 \rangle=B{\rm Spin}
    \ar[rr]^{\hspace{5mm}\frac{1}{2}p_1}
  \ar[d]
  &&
K(\Z, 4)
\ar[rr]
&&
K(\Z_2, 4)
\ar[d]^{=}
\\
&&
(BO)\langle 2\rangle =BSO
\ar[rrrr]^{w_4}
&&
&&
K(\Z_2, 4),
  }
\)
where the dashed arrow is a lift of the given $\mathrm{Spin}$-structure.
\begin{observation}
The class $\frac{1}{4}p_1$ is the obstruction to lifting
 a $\mathrm{String}^{2\mathrm{DD}_2}$ bundle,  
 to a $\mathrm{String}$-bundle. Here,
  $\mathrm{String}^{2\mathrm{DD}_2}$ is the loop space of 
  the space $B\mathrm{String}^{2\mathrm{DD}_2}$, defined above.
\end{observation}
\proof
  The pasting law for homotopy pullback implies that the 
  homotopy fiber of $\frac{1}{4}p_1$ is $B \mathrm{String}$.
  $$
    \raisebox{40pt}{
    \xymatrix{
	  B \mathrm{String}
	  \ar[r]_<{\hpull}
	  \ar[d]
	  &
	  {*}
	  \ar[d]^0
	  \\
	  B \mathrm{String}^{2 \mathrm{DD}_2}
	  \ar[r]_<{\hpull}^{\hspace{4mm}\frac{1}{4}p_1}
	  \ar[d]
	  &
	  K(\mathbb{Z}, 4)
	  \ar[d]^{\times 2}
	  \\
	  B \mathrm{Spin}
	  \ar[r]
	  &
	  K(\mathbb{Z},4)
	}
	}
	\,.
  $$
\endofproof
This observation is analogous to proposition 2 in \cite{SSS2} for the Fivebrane
case, where there we were considering the comparison of $\frac{1}{48}p_2$
to the obstruction to Fivebrane structure, given by $\frac{1}{6}p_2$. We will discuss this 
further later in the paper.

Similarly, with $a : B E_8 \to B^3 U(1)$ the canonical universal 4-class for 
$E_8$-bundles and $X$ a manifold of dimension $\mathrm{dim}X \leq 14$ we have
that a
$\mathrm{Spin}$-structure on $X$ lifts to a $\mathrm{String}^{2a}$-structure
precisely if $\frac{1}{2}p_1$ is further divisible by 2. 
\(
  \xymatrix{
   & B \mathrm{String}^{2a} \ar[d] \ar[dr]^{\frac{1}{4}p_1} \ar[r] 
   & B E_8 \ar[d]^{2 a}
   \\
   X \ar[r] \ar@{-->}[ur] & B \mathrm{Spin} \ar[r]^{\hspace{-2mm}\frac{1}{2}p_1} & B^3 U(1)
  }
  \,.
\)
Using this we can now reformulate the anomaly cancellation condition 
(\ref{WittenQuantizationCondition})
as follows.
\begin{definition}
  \label{TwistedString2sStructures}
  For $X$ a manifold and for 
  $[\alpha] \in H^4(X,\mathbb{Z})$ a cohomology class, the space
  $(\frac{1}{2}p_1-2a)\mathrm{Struc}_{[\alpha]}(X)$ 
  of \emph{$[\alpha]$-twisted $\mathrm{String}^{2a}$-structures} on $X$
  is the homotopy pullback
  $$
    \raisebox{20pt}{
	\xymatrix{
	  (\frac{1}{2}p_1\!-\!2a)\mathrm{Struc}_{[\alpha]}(X)
	  \ar[rr]_<{\hpull}
	  \ar[d]
	  &&
	  {*}
	  \ar[d]^\alpha
	  \\
	  \mathrm{Maps}(X, B (\mathrm{Spin} \times E_8) )
	  \ar[rr]^{\hspace{3mm}\frac{1}{2}p_1 - 2a}
	  &&
	  \mathrm{Maps}(X, B^3 U(1) )
	}
	}
	\,,
  $$
  where the right vertical map picks a cocycle $\alpha$ representing the class $[\alpha]$.
\end{definition}
This can be viewed as a ``twist for the twisted String structure".
  In terms of this definition, we have
\begin{observation}
 Condition (\ref{WittenQuantizationCondition}) 
  is precisely the conditon guaranteeing a 
  lift of the given $\mathrm{Spin}$- and the given $E_8$-principal
  bundle to a \emph{$[G_4]$-twisted $\mathrm{String}^{2a}$-structure}
  through the left vertical map from def. \ref{TwistedString2sStructures}.
\end{observation}

\vspace{3mm}
\subsubsection{M-theory with boundary: Heterotic M-theory}

In \cite{HW}, Horava and Witten carefully analyzed the map $\partial_1$
in diagram (\ref{UVAndBoundaryCompletion})
and gave arguments on how it must extend to $\partial_2$. 
If we denote by
$Q := \partial Y \hookrightarrow Y$ the boundary inclusion,
then the condition they find is a boundary
condition on the $C$-field, saying that the restriction of its
4-class to $Q$ has to vanish,
\(
  [G_4]|_Q = 0
  \,.
\)
This implies that over $Q$ the anomaly-cancellation condition 
(\ref{WittenQuantizationCondition}) becomes
\(
  \frac{1}{2}p_1(TY)|_Q = 2 a(E)|_Q \;\; \in H^4(Q,\mathbb{Z})
  \label{CFieldConditionOverBoundary}
  \,.
\)
This 
is equivalent, in direct analogy with (\ref{HomotopyOfTwistedSpinCStructure}),
to the existence of a homotopy of the form
\(
  \label{DiagramForStringTrivializationOnBrane}
  \raisebox{20pt}{
  \xymatrix{
    Q \ar[rr]_>{\ }="s"
	\ar[d]
	\ar[drr]_{2 a(E)|_{Q}}^{\ }="t"
	&&
	B \mathrm{Spin}
	\ar[d]^{\frac{1}{2}p_1}
	\\
	X
	\ar[rr]_{2a(E)}
	&&
	K(\mathbb{Z}, 4)
	\ar@{=>} "s"; "t"
  }}
  \,.
\)
Essentially in this form twisted string structures in the context of
string theory 
were proposed in (the first eprint version of) \cite{Wang}.
We see here a general pattern of twisted structures occuring as
\emph{relative trivializations on branes}.

Notice that on $Q$ this is the Green-Schwarz anomaly cancellation condition
(\ref{Green Schwarz anomaly}) of the heterotic string, but 
refined by a further cohomological divisibility condition.
The following statement says that this may equivalently be reformulated 
in terms of $\mathrm{String}^{2a}$ structures.
\begin{proposition}
  For $E \to Y$ a fixed $E_8$-bundle, 
  we have an equivalence
  $$
    \mathrm{Maps}(Y, B \mathrm{String}^{2a})|_{E}
    \simeq
    (\frac{1}{2}p_1)\mathrm{Struc}(Y)_{[2a(E)]}
  $$
  between, on the right, the space of $2a(E)$-twisted $\mathrm{String}$-structures
from def. \ref{SpaceOfTwistedStringStructures}, and, on the left, the space of 
 $\mathrm{String}^{2a}$-structures with fixed class $2a$, hence the homotopy pullback
 $\mathrm{Maps}(Y, B \mathrm{String}^{2a}) \times_{\mathrm{Maps}(Y, B E_8)} \{E\}$.
 \label{EquivalenceOfString2aAndTwistedStringStructure}
\end{proposition}
\proof
  Consider the diagram
  $$
    \xymatrix{   
	  \mathrm{Maps}(Y, \mathrm{String}^{2a})|_E
      \ar[rr]_<{\hpull}
	  \ar[d]
	  &&
	  {*}
	  \ar[d]^{E}
	  \\
	  \mathrm{Maps}(Y, \mathrm{String}^{2a})
	  \ar[rr]_<{\hpull}
	  \ar[d]
	  &&
	  \mathrm{Maps}(Y, B E_8)
	  \ar[d]^{\mathrm{Maps}(X, 2a)}
	  \\
	  \mathrm{Maps}(Y, B \mathrm{Spin})
	  \ar[rr]^{\mathrm{Maps}(Y, \frac{1}{2}p_1)}
	  &&
	  \mathrm{Maps}(Y, B^3 U(1))
	}
  $$
  The top square is a homotopy pullback by definition.
  Since $\mathrm{Maps}(Y,-)$ preserves homotopy pullbacks
  (for $Y$ a manifold, hence a CW-complex), 
  the bottom square is a homotopy pullback
  by definition \ref{Stringc}. Therefore, by the pasting law,
also the total rectangle is a homotopy pullback. With
 def. \ref{SpaceOfTwistedStringStructures} this implies the claim.    
\endofproof

Therefore the boundary anomaly cancellation condition for the M2-brane
has the following equivalent formulation.
\begin{observation}
  For $X$ a $\mathrm{Spin}$-manifold equipped with a
  complex vector bundle $E \to Y$, condition 
  (\ref{CFieldConditionOverBoundary}) precisely guarantees the
  existence of a lift to a $\mathrm{String}^{2a}$-structure
  through the left vertical map in the proof of prop.
  \ref{EquivalenceOfString2aAndTwistedStringStructure}.
\end{observation}


\subsection{Twisted Fivebrane structures and physical applications}
\label{twisted fivebrane}

We discuss \emph{twisted fivebrane structures} and their role in anomaly
cancellation in
\begin{enumerate}
  \item the NS-5-brane and dual heterotic string theory;
  \item the M5-brane.
\end{enumerate}

\subsubsection{The NS-5-brane}
\label{NS5 and twisted Fivebrane Structures}

The magnetic dual of the (heterotic) string is the NS-5-brane.
Where the string is electrically charged under the $B_2$-field
with class $[H_3] \in H^3(X,\mathbb{Z})$, the NS-5-brane is electrically
charged under the $B_6$-field with class $[H_7] \in H^7(X, \mathbb{Z})$
\cite{Cham}.
As we discuss in detail shortly,
in the presence of a $\mathrm{String}$-structure,
hence when $\frac{1}{2}p_1(T X) = 0$, the anomaly
of the 5-brane $\sigma$-model vanishes  
if the background fields satisfy
\( 
   \frac{1}{6}p_2(T X)
   =   
   8\hspace{0.5mm}{\rm ch}_4(E)
   \;\;
   \in H^8(X, \mathbb{Q})
  \label{NS5BraneAnomalyCancellationCondition} 
  \,,
\) 
where $E$ is an 
$E_8 \times E_8$- or $\mathrm{Spin}(32)/\mathbb{Z}_2$-principal bundle,
and $\mathrm{ch}(E)$ denotes its Chern character, and where
$\frac{1}{6}p_2$ denotes the \emph{second fractional Pontrjagin class}  
given by the following classical fact (see \cite{Bott}).
\begin{fact}
  The second Pontrjagin class $p_2 \in H^8(B \mathrm{SO}, \mathbb{Z})$
  becomes divisible precisely by 6 when pulled back to $B \mathrm{String}$.
  The corresponding preimage under multiplication by six, denoted
  $\frac{1}{6}p_2$,
  is a generator of the group $H^8(B \mathrm{Spin}, \mathbb{Z}) \cong \mathbb{Z}$. 
  \label{SecondFractionalPontrjagin}
\end{fact}

It is clear now that a discussion entirely analogous to that of
section \ref{HeteroticStringAndTwistedStringStructures} applies. 
For the untwisted case the following terminology was introduced in
\cite{SSS2}.
\begin{definition}
  Write $\mathrm{Fivebrane}$ for the loop group of the homotopy
  fiber $B \mathrm{Fivebrane}$ of a representative $\frac{1}{6}p_2$
  of the universal second fractional Pontrjagin class
  $$
    \raisebox{20pt}{
    \xymatrix{ 
	  B \mathrm{Fivebrane} \ar[r]_<{\hpull} \ar[d] & {*} \ar[d]
	  \\
	  B \mathrm{String}
	  \ar[r]^{\frac{1}{6}p_2}
	  &
	  B^7 U(1)
	}
	}
	\,.
  $$  
\end{definition}
In direct analogy with def. \ref{SpaceOfTwistedStringStructures} we therefore have the following
notion.
\begin{definition}
  \label{TwistedFivebraneStructures}
  For $X$ a manifold and $[c] \in H^8(X,\mathbb{Z})$ a class, 
  we say that the 
  \emph{space of $[c]$-twisted $\mathrm{Fivebrane}$-structures}
  on $X$, denoted $(\frac{1}{6}p_2)\mathrm{Struc}_{[\alpha]}(X)$, is the 
  homotopy pullback
  $$
    \raisebox{20pt}{
    \xymatrix{
      (\frac{1}{6}p_2)\mathrm{Struc}_{[c]}(X)
      \ar[rr]_<{\hpull}
      \ar[d]
      &&
      {*}
      \ar[d]^{c}
      \\
      \mathrm{Maps}(X, B\mathrm{String})
      \ar[rr]^{\mathrm{Maps}(X,\frac{1}{6}p_2)}
      &&
      \mathrm{Maps}(X, B^7 U(1))
    }
	}
    \,,
  $$
\end{definition}
Explicitly, a $[c]$-twisted Fivebrane structure
on a brane $\iota : Q \to X$ equipped 
with String structure $f : Q \to B \mathrm{String}$ is
a homotopy $\eta$ in a diagram analogous to (\ref{DiagramForStringTrivializationOnBrane})
  \(
  \raisebox{20pt}{
  \xymatrix{
    Q \ar[rr]_>{\ }="s"
	\ar[d]
	\ar[drr]_{c|_{Q}}^{\ }="t"
	&&
	B \mathrm{String}
	\ar[d]^{\frac{1}{6}p_2}
	\\
	X
	\ar[rr]_{c}
	&&
	K(\mathbb{Z}, 8)
	\ar@{=>}^\eta "s"; "t"
  }}
    \,.
  \label{5tX}
  \)
Two  $[c]$-twisted Fivebrane
structures  $\eta$ and
$\eta'$ on $Q$ are regarded as equivalent 
if there is a homotopy between $\eta$ and $\eta'$.
In the case that $[c] = 0$ this reduces to the untwisted Fivebrane
structures considered in \cite{SSS2}.

In terms of these notions we now have
\begin{observation}
 For $X$ a manifold with $\mathrm{String}$-structure
 and with a background gauge bundle $E \to X$ fixed such that 
 $8\mathrm{ch}(E)$
 is integral, 
 condition (\ref{NS5BraneAnomalyCancellationCondition})
 is precisely the condition for the existence of 
\emph{$8\,\mathrm{ch}(E)$-twisted $\mathrm{Fivebrane}$-structure} on $X$.
\end{observation}

We now consider the above
anomaly cancellation condition 
in more detail. 

In
\cite{SSS2} the main example of a Fivebrane structure came from the
dual formulation \cite{SS} \cite{Gates} of the Green-Schwarz anomaly
cancellation mechanism \cite{GS}, using the dual $H$-field $H_7$ of
\cite{Cham}. The expression is given by 
\( 
dH_7= 2\pi \left[ {\rm ch}_4(F_A)
-\frac{1}{48} p_1(F_{\omega}) {\rm ch}_2(F_A) + \frac{1}{64} p_1(F_{\omega})^2
-\frac{1}{48}p_2(F_{\omega}) 
\right]\;, \label{dH72} 
\) 
where $F_A$ and $F_\omega$ are the curvatures of the connections $A$ and 
$\omega$ on the gauge bundle $E$ and the tangent (or Spin) bundle of the ten-manifold
$M$, 
respectively. 
In order to define a
Fivebrane structure, we assume we already have a String structure,
 so we require $\frac{1}{2}p_1(TM)=0$. Then the expression
(\ref{dH72}) becomes
\(
dH_7= 2 \pi \left[ {\rm ch}_4(F_A) -
\frac{1}{48}p_2(F_{\omega}) \right]\; .
\label{dH7v2}
\)
In \cite{SSS2} we had
to find ways to get rid of the extra terms to isolate the
non-decomposable terms. In the twisted formalism in this paper we
see that the presence of such terms amounts to a piece of the twist
and that it does not matter how many terms we have, as long as they
are integral and 
have the same total degree and hence provide a map to $K(\Z, 8)$.
Indeed, if we can define
 \( 
 [\beta] := - {\rm ch}_4 (E) : \xymatrix{ M
\ar[r]^{\hspace{-3mm}!} & K(\Z, 8)}\;, \label{!} 
\) 
i.e. require factorization 
\(
\xymatrix{
M 
\ar[rr] 
\ar@{..>}[dr]_{[\beta]}
&&
K(\Q, 8)
\\
&
K(\Z, 8)
\ar@{^{(}->}[ur]
&
}\;,
\)
then we can reinterpret
expression (\ref{dH7v2}) as $\frac{1}{48}p_2(TM) + [\beta] =0$,
since $[dH_7]=0$, the cohomology class of an exact form.

\vspace{3mm}
 We discuss the validity of the map in (\ref{!}). The Chern character is
 in general not an integral expression, but rather
 \(
{\rm  ch} : K^0(X)  \to H^{\rm even}(X; \Q).
 \)
One way out of this is to first define a rational version of the twist, for
which the map in (\ref{!}) is replaced by a map from $M$ to
the rational Eilenberg-MacLane space
\(
[\beta] : = -{\rm ch}_4(E): M \to K(\Q, 8),
\)
which gives that indeed ${\rm ch}_4(E)$ is in general in $[M, K(\Q, 8)]=H^8(M, \Q)$.
Hence
\begin{definition}
A rational Fivebrane twist on $M$ is a map from $M$ to $K(\Q, 8)$, i.e.
an element of $H^8(M; \Q)$.
\end{definition}
However, we can also give conditions under which the map in (\ref{!})
is valid.
The degree four Chern character is given by
\(
{\rm ch}_4=\frac{1}{24} \left( c_1^4 -4c_1^2 c_2 + 4c_1 c_3 + 2c_2^2 -4c_4 \right)\;.
\)
The Chern classes are integral classes and so the Chern character is
a priori integral up to a factor of 24.

\vspace{3mm}
We describe this as follows.
 The
Chern character is not integral in $BU$ but it will be 
integral in some
lift, say $\mathcal{BU}$, of $BU$. Then we ask: when can we
lift to this new space? This is given in terms of the following diagram
\(
\raisebox{40pt}{
\xymatrix{
&&
&&
K(\Z_{24}, 7)
\ar[d]
&&
\\
&&
\mathcal{BU}
\ar[d]
\ar[rr]^{{\rm ch}_4}_<{\hpull}
&&
K(\Z, 8)
\ar[d]^{ \times 24}
&&
\\
M
\ar@{..>}[urr]^{f}
\ar[rr]
&&
BU
\ar[rr]^{24 {\rm ch}_4}
&&
K(\Z, 8)
\ar[rr]
&&
K(\Z_{24}, 8) \; .
}
}
\)
The right-most factor $K(\Z_{24}, 8)$ represents
the obstruction: there is a class $k$ in $H^8(M ;\Z_{24})$
which measures this obstruction. The top-most factor $K(\Z_{24}, 7)$
represents the different labeling of lifts $f$ to the new space $\mathcal{BU}$.
If we take connected covers of $BU$ rather than $BU$ itself in the diagram, then
we have that the space $\mathcal{BU}$ is isomorphic to another space
in which $\frac{1}{6}c_4$, instead of ${\rm ch}_4$, is integral.
The relevance of the unitary groups here is because they provide the adjoint representation
for our structure groups and this is the representation relevant for Yang-Mills theory.
For $E_8$,  the adjoint representation is ad : $E_8 \to SU(248)$, so that the
adjoint representation of $G=E_8 \times E_8$ is
$
({\rm ad}, {\rm ad}) : E_8 \times E_8 \to SU(248) \times SU(248)
\hookrightarrow SU(496)
$.

\vspace{3mm}
Note that the above general discussion can be simplified.
For both structure groups $E_8 \times E_8$ and ${\rm Spin}(32)/\Z_2$
we have $c_1(E)=0$, so that for these groups
$
{\rm ch}_4(E)=\frac{1}{12} (c_2(E)^2 - 2c_4(E))
$.
We now consider  the case when,  in addition, we have
$c_2(E)=0$. In this case,
the formula for the Chern character ${\rm ch}_4(E)$ further simplifies to
\(
{\rm ch}_4(E) = -\frac{1}{6}c_4(E).
\label{c4}
\)

\vspace{3mm}
Here what we have really done is lifted the unitary group to its connected cover
$BU \langle 8 \rangle$. Indeed let us consider the result from \cite{Sing} where
the mod $p$ ($p$ an odd prime) cohomology of the connective cover $BU\langle 2n \rangle$
was calculated. From that result and the result of Stong \cite{Stong} for $p=2$,
the following divisibility result was deduced for all primes $p$ in \cite{Sing}. Let
$c_k \in H^{2k}(BU ; \Z)$ be the universal Chern class in $BU$, then the
Chern class $r_n^* (c_k)$ in $BU\langle 2n \rangle$ where
$r_n : BU \langle 2n \rangle \to BU$ be the canonical projection
is divisible by \cite{Sing}
\(
\prod_p p^q
\)
where $q$ is the least integer part of $\frac{(n-1) -\sigma_p(k-1)}{p-1}$,
with $\sigma_p(n)=\sum a_i$ the sum of the coefficients in the unique
decomposition of the integer $n$ as $n=a_0 + a_p + \cdots + a_k p^k$,
with $a_i < p$. Applying this result for $n=4$, $p=2, 3$, and using
$\sigma_2(3)=2$, $\sigma_3(3)=1$, we get that $r_4^*(c_4)$ is divisible
by
\(
2^{\frac{2 - \sigma_2(3)}{1}} \cdot 3^{\frac{3 - \sigma_3(3)}{2}}=6\;.
\label{6}
\)

\vspace{3mm}
We will give an example where this occurs and where
the expression (\ref{c4}) is integral.

\vspace{3mm}
\noindent {\bf Example.}
Consider a complex vector bundle $E$ on the eight-sphere $S^8$. For
ten-manifold we can simply take $S^8 \times \R^2$ for example.
The index of the Dirac operator on $S^8$
coupled to the vector bundle $E$ is given by the evaluation of the
twisted $\widehat{A}$-genus
$\widehat{A}(S^8, E):= \left( {\rm ch}(E)\cdot \widehat{A}(S^8)\right)[S^8]$
on the fundamental class $[S^8]$ of $S^8$
\(
{\rm Index} D_E= \left( \widehat{A} (TS^{8}) \cdot {\rm ch}(E) \right)
= {\rm ch}(E) [S^8],
\)
as $\widehat{A}(TS^8)=1$, since spheres have stably trivial tangent
bundles. Since $S^8$ is a Spin manifold, the index should be
an integer. This then gives the requirement
\(
{\rm ch}_4(E)[S^8]= -\frac{1}{6} c_4(E) [S^8] \in \Z.
\)


\subsubsection{The M5-brane and the dual $C$-field}
\label{Twisted Fivebrane2a2aStructures}
\label{The dual C-field}

The magnetic dual of the M2-brane is the M5-brane. Where the M2-brane is electrically charged
under the $C_3$-field with class $[G_4] \in H^4(X, \mathbb{Z})$,
the M5-brane is electrically charged under the dual $C_6$-field with 
class $[G_8] \in H^8(X,\mathbb{Z})$.
If $X$ admits a $\mathrm{String}$-structure, 
then, as we discuss in more detail in a moment, one finds an 
anomaly cancellation condition on these background fields analogous to (\ref{WittenQuantizationCondition})
which reads
\(
  8 [G_8] = 4 a(E) \cup a(E) -  \frac{1}{6}p_2(T X)
  \label{M5braneAnomalyCancellationCondition}
  \,,
\)
where $a : B E_8 \to K(\mathbb{Z},4)$ is a representative of the 
canonical degree-4 class on the classifying space of the exceptional
Lie group $E_8$.

The homotopy-theoretic interpretation of this
condition involves the following $\mathrm{Fivebrane}$-analog 
of $\mathrm{Spin}^c$ as it appeared in prop. \ref{SpinCByHomotopyPullback}, 
and of $\mathrm{String}^\alpha$ as it was considered in def. \ref{StringC}.
\begin{definition}
  For $G$ a topological group and $[\alpha] \in H^8(B G, \mathbb{Z})$
  a universal 8-class, we say that $\mathrm{Fivebrane}^\alpha$ is the loop group
  of the homotopy pullback
  $$
    \raisebox{20pt}{
    \xymatrix{
	  B \mathrm{Fivebrane}^\alpha \ar[r]_<{\hpull} \ar[d] 
	   & B G \ar[d]^\alpha
	  \\
	  B \mathrm{String}
	  \ar[r]^{\frac{1}{6}p_2}
	  &
	  B^7 U(1)
	}
	}
	\,.
  $$
\end{definition}
As before with  $\mathrm{String}^\alpha$, 
this notion in particular subsumes spaces on which 
the class $\frac{1}{6}p_2$ is further divisble. For let 
$G = B^6 U(1)$ in the above and take $\alpha := \mathrm{DD}_6$ 
to be the canonical
degree-8 class on $B^7 U(1) \simeq K(\mathbb{Z},8)$. This is 
the canonical class of circle 7-bundles / U(1)-bundle 6-gerbes.
Then $\mathrm{Fivebrane}^{8\mathrm{DD}_6}$ is the loop group of
the homotopy pullback
\(
  \xymatrix{
    B \mathrm{Fivebrane}^{8 \mathrm{DD}_6}
	\ar[r]^>>>>>{\frac{1}{48}p_2}_<{\hpull}
	\ar[d]
	&
	K(\mathbb{Z}, 8)
	\ar[d]^{\times 8}
	\\
	B \mathrm{String}
	\ar[r]^{\frac{1}{6}p_2}
	&
	K(\mathbb{Z}, 8)
  }
  \,.
\)
Accordingly, we have that a $\mathrm{String}$-bundle
lifts to a  $\mathrm{Fivebrane}^{8\mathrm{DD}_6}$-bundle precisely
if the class of $\frac{1}{6}p_2$ is further divisible by 8, 
hence if $\frac{1}{48}p_2$ exists. In \cite{SSS2} this space
was denoted 
\(
  \mathcal{F}^{\langle 8\rangle}
  :=
  B \mathrm{Fivebrane}^{8\mathrm{DD}_6}
  \,.
\)
Moreover, the class $\frac{1}{48}p_2 \in H^8(\mathrm{Fivebrane}^{8 \mathrm{DD}_6}; \mathbb{Z})$ is the universal obstruction to lifting to a genuine 
$\mathrm{Fivebrane}$-bundle.
Accordingly, we have
a notion of twisted $\mathrm{Fivebrane}$-structures induced
by $\frac{1}{48}p_2$. They form the space given by the
homotopy pullback.
  \(
    \xymatrix{
	  \frac{1}{48}p_2 \mathrm{Struc}_{[\beta]}(X)
	  \ar[r]_<{\hpull} 
	  \ar[d]
	  & 
	  {*}
	  \ar[d]^{\beta}
	  \\
	  \mathrm{Maps}(X, B \mathrm{Fivebrane})
	  \ar[r]^{\hspace{3mm}\frac{1}{48}p_1}
	  &
	  \mathrm{Maps}(X, B^7 U(1))
	}\;.
 \)
Such a $[\beta]$-twisted $\mathrm{Fivebrane}$-structure
on a brane $\iota : Q \hookrightarrow X$ 
 is a homotopy $\eta$  
in the diagram
\(
\xymatrix{
Q
\ar[rr]^{\nu}_>{\ }="s"
\ar[d]^{\iota}
&&
   B \mathrm{Fivebrane}^{8 \mathrm{DD}_6}
\ar[d]^{\frac{1}{48}p_2}
\\
X
\ar[rr]_{\beta}^<{\ }="t"
&&
K(\Z, 8)
\ar@{=>}^\eta "s"; "t"
}
\,,
\)
which exists precisely if 
 \(
 \frac{1}{48}p_2(X) + \iota^* ([\beta])=0
 \,.
 \) 
 We may also regard the situation 
 in analogy with def. \ref{TwistedString2sStructures} and consider the following.
\begin{definition}
  \label{TwistedFivebrane2a2aStructures}
  For $X$ a manifold and for 
  $[c] \in H^8(X;\mathbb{Z})$ a degree 8 cohomology class, the space
  \newline
  $(\frac{1}{6}p_2- 2a\cup 2a)\mathrm{Struc}_{[c]}(X)$ 
  of \emph{$[c]$-twisted $\mathrm{Fivebrane}^{2a \cup 2a}$-structures} on $X$
  is the homotopy pullback
  $$
    \raisebox{20pt}{
	\xymatrix{
	  (\frac{1}{6}p_2-2a\cup 2a)\mathrm{Struc}_{[c]}(X)
	  \ar[rr]
	  \ar[d]
	  &&
	  {*}
	  \ar[d]^c
	  \\
	  \mathrm{Maps}(X, B \mathrm{String} \times E_8 )
	  \ar[rr]^{~~~~\frac{1}{6}p_2 - 2a\cup 2a}
	  &&
	  \mathrm{Maps}(X, B^7 U(1) )
	}
	}
	\,,
  $$
  where the right vertical map picks a cocycle $c$ representing the class $[c]$.
\end{definition}
In terms of these notions we thus see that 
\begin{observation}
Over a manifold $X$ with $\mathrm{String}$-structure
and with a fixed gauge bundle $E$,
condition (\ref{M5braneAnomalyCancellationCondition}) 
is precisely the condition that guarantees 
existence of a lift to 
\emph{$[8G_8]$-twisted $\mathrm{Fivebrane}^{2a \cup 2a}$-structure}
through the left vertical morphism in def. 
\ref{TwistedFivebrane2a2aStructures}.
\end{observation}

We discuss now how (\ref{M5braneAnomalyCancellationCondition})
arises in more detail.
Locally the $C_3$-field is given by a 3-form, which 
is traditionally denoted by the same symbol.
The equation of motion for $C_3$ is
obtained from varying the action 
\( S(C_3)= \int_{Y} \left[ G_4
\wedge * G_4 + \frac{1}{6}G_4 \wedge G_4 \wedge C_3 - I_8 \wedge C_3
\right] \) 
on an eleven-dimesional Spin manifold $Y$ to obtain 
\( d*G_4 =
- \frac{1}{2} G_4 \wedge G_4 + I_8\;. \label{EOM} 
\) 
Here $I_8$ is the
one-loop polynomial \cite{VW} \cite{DLM} given in terms of the Pontrjagin
classes of the tangent bundle $TY$ to $Y$
\(
I_8=\frac{p_2(TY) -
\frac{1}{2}(\frac{1}{2}p_1(TY))^2}{48}\;,
\label{1loop}
\)
and $*$ is the Hodge duality operator for a given metric in eleven dimensions.

\vspace{3mm}
The integral lift of (\ref{EOM}) leads to a class defined in
\cite{DFM}
\bea
[G_8] &=&\left[  \frac{1}{2} G_4^2 -I_8 \right]
\nonumber\\
&=& \frac{1}{2} a (a -\lambda) + \frac{7 \lambda^2 - p_2}{48},
\label{TH}
\eea
where $\lambda=\frac{1}{2}p_1$, and $a$ is the degree four class
of an $E_8$ bundle coming from Witten's shifted
quantization condition for $G_4$ \cite{Flux}
\(
\label{g4 quant}
[G_4]= a -\frac{1}{2}\lambda= a -\frac{1}{4} p_1\;.
\)

\vspace{3mm}
In \cite{Among} Witten interpreted the vanishing of a certain torsion
class $\theta$ on the M-fivebrane worldvolume as a necessary condition for the
decoupling of the fivebrane from the ambient space ( ``the bulk"). Hence the
vanishing of $\theta$ meant that the fivebrane can have a well-defined partition
 function.
  Consider the embedding $\iota : W \hookrightarrow Y$
 of the fivebrane with six-dimensional worldvolume $W$ into eleven-dimensional
 spacetime $Y$. Consider the ten-dimensional
 unit sphere bundle $\pi : X \to W$ of $W$ with fiber $S^4$ associated to the
 normal bundle $N \to W$ of the embedding $\iota$. Then it was shown in
 \cite{DFM} that the integration of $G_8$ over the fiber of $X$ gives exactly
 the torsion class $\theta$ on the fivebrane worldvolume
 \(
 \theta= \pi_*(G_8) \in H^4(W; \Z).
 \)
 Therefore, the vanishing of $G_8$ is a necessary condition for
 the existence of a non-zero partition function \cite{DFM}.

 \vspace{3mm}
 We now proceed with the interpretation. Since we have Fivebrane
 structures in mind, we assume that $Y$ already admits a String
 structure, i.e. that $\frac{1}{2}p_1(Y)=0$. Then, from (\ref{TH})
we see that the class $G_8(Y)$ simplifies to
 \(
   \label{dual of c}
 G_8 (Y) =\frac{1}{2}a^2 - \frac{1}{48} p_2(Y).
 \)
 The class $a$ is an integral class of an $E_8$ bundle and hence defines
 a map to $K(\Z, 4)$. Then the square of $a$ defines a map to $K(\Z, 8)$,
 and hence defines a twist for us. As we also have the class
 $\frac{1}{48}p_2$, then we have a twist for the modified Fivebrane
 structure. 

\vskip3ex
\noindent{\bf Necessity of the Fivebrane condition?}
The Fivebrane condition is stronger than simply the requirement that the one-loop
term $I_8$ to vanish. For the former we require the obstructions $\frac{1}{2}p_1$ 
and $\frac{1}{6}p_2$ vanish {\it separately}, whereas for the latter we only 
require the combination to vanish. This has been studied in \cite{IPW}
\cite{IP} \cite{Witt}.
For instance, following \cite{Witt},
a Riemannian 8-dimensional spin manifold $M^8$ is said to be {\em doubly
supersymmetric} if and only if the tangent bundle $TM^8$ and
the spinor bundles $\Delta_+M^8$ and $\Delta_-M^8$ are associated with a
principal $G$-fiber bundle such that there exist $G$-invariant
isomorphisms between any two of the three bundles ,
 i.e. $TM^8=\Delta_+M^8=\Delta_-M^8$.
If $M^8$ is doubly supersymmetric,
\(
w_1=w_2=0,\qquad e=0, \qquad
4p_2=p_1^2\;,
\)
where $e$ is the Euler class. 
Then this implies for the signature
$
sgn(M^8)=16\hat{A}[M^8]
$.
In particular, $sgn(M^8)\equiv 0\:\mod 16$.
One example is $PSU(3)$-structure for which
\bea
w_i&=&0 ~(i \neq 4), \qquad \quad w_4^2=0 \\
e&=&0, \qquad  \qquad \qquad  p_1^2=4p_2.
\eea
In particular, all Stiefel-Whitney numbers vanish.

\vspace{3mm}
A second example is a differentiable 8-fold $M^8$ with 
an odd topological generalized ${\rm Spin}(7)$-structure
(in the sense of \cite{Witt}) for which 
\begin{equation}\label{obgenspin72}
\chi(M^8)=0,\qquad p_1(M^8)^2-4p_2(M^8)=0.
\end{equation}
The 7-sphere admits a Spin structure and therefore admits a generalized
$G_2$-structure. The tangent bundle of the 8-sphere is stably
trivial and therefore all the Pontrjagin classes vanish. Since the
Euler class is non-trivial, there exists no generalized
$Spin(7)$-structure on an 8-sphere. However,
equation~(\ref{obgenspin72}) is
automatically satisfied for manifolds of the form $M^8=S^1\times
N^7$ with $N^7$ Spin.


\section{Twisted differential structures in String theory}
\label{Twisted differential String- and Fivebrane structures}

We now indicate a theory of nonabelian \emph{differential cohomology}
in which the topological structures considered in section 
\ref{TwistedTopologicalStructuresInStringTheory} have smooth and 
differential refinements. A full account of this theory is 
given in \cite{nactwist}.
The complete differential refinements are based on a smooth refinement
of topological spaces by structures called \emph{higher smooth stacks}
or \emph{higher smooth groupoids}.
The full constructions of these, for the cases considered here, 
are discussed in \cite{FSSI}. 
As shown there, underlying any such differential cohomological structure is, 
locally, explicit differential geometric data, given by differential forms
with values in $L_\infty$-algebras. 
This \emph{$L_\infty$-algebra valued connection} data by itself was
discussed in \cite{SSS1}, a brief collection of relevant $L_\infty$-algebraic
notions is in the appendix section \ref{LinfinityNotions}. 

Here, after a brief survey of the general theory, we concentrate on a discussion
of this local differential form data for the differential refinements
of the twisted $\mathrm{String}$- and $\mathrm{Fivebrane}$-structures
from section \ref{TwistedTopologicalStructuresInStringTheory}. 
We show that these reproduce the equations on differential forms that 
are traditional in the physics literature on the Green-Schwarz anomaly
cancellation as well as its magnetic dual. 

The key fact is the existence of definition \ref{TwistedDifferentialStructures}
below, which gives smooth and differential refinements of the
homotopy pullbacks that defined twisted $\mathrm{String}$- and
$\mathrm{Fivebrane}$-structures in def. \ref{TwistedStringStructures}
and def. \ref{TwistedFivebraneStructures}, respectively.

\subsection{Differential twisted cohomology}
\label{differential twisted cohomology}

Differential twisted cohomology is the pairing of the notions of 
\emph{twisted cohomology}
with \emph{differential cohomology} (see \cite{HS}).
We have that 
\begin{itemize}
\item a cocycle in differential cohomology is to the underlying bare cocycle as a 
connection on an $\infty$-bundle is to the underlying principal $\infty$-bundle;

\item a cocycle in twisted cohomology is to an ordinary cocycle as a 
twisted bundle is to a principal bundle.
\end{itemize}
We indicate now briefly a formal definition of such objects
as described in full detail in \cite{nactwist}, showing how they
connect via the constructions from \cite{FSSI}
the cohomological discussion that we had so far 
 to the 
$L_\infty$-algebraic differential form data of 
\cite{SSS1} in terms of which we shall obtain the relevant twisted
Bianchi identities.

The basic mechanism is to refine the homotopy theory 
by passing from the $\infty$-topos \cite{Lurie} 
$\mathrm{Top} \simeq \infty \mathrm{Grpd}$
of topological spaces -- equivalently: geometrically discrete
$\infty$-groupoids -- to that of \emph{smooth $\infty$-groupoids},
defined as the $\infty$-category of $\infty$-stacks over 
the site of smooth manifolds.

\medskip

It is a familiar fact that \emph{Lie groupoids}, 
such as for instance orbifolds, are naturally to be thought of 
as \emph{stacks} on the site of smooth manifolds, often called
\emph{differentiable stacks} \cite{BehrendXu}. Such differentiable
stacks form a common generalization of smooth manifolds and a 
small fragment of homotopy theory. For instance, for each 
Lie group $G$ there is a Lie groupoid, denoted ${*}//G$ or $\mathbf{B}G$, 
with a single object
and $G$ as its space of morphisms, such that for $X$ any smooth manifold,
regarded as a Lie groupoid, the collection $\mathbf{H}(X, \mathbf{B}G)$ 
of stack morphisms 
$X \to \mathbf{B}G$ forms the groupoid of smooth $G$-principal bundles
on $X$ and smooth gauge transformations between them. So the Lie groupoid
$\mathbf{B}G$ serves as a smooth refinement of the topological classifying
space $BG$.

But in the context of ordinary stacks, not all the classifying spaces
that we considered in the previous section, such as 
$B^n U(1) \simeq K(\mathbb{Z},n+1)$ for higher $n$,
have a smooth refinement to stacks on manifolds. 
But there is a natural generalization of the
theory of stacks to a general theory of \emph{higher stacks}, often called
\emph{$\infty$-stacks}. The collection of all of these
over a given site is an \emph{$\infty$-topos} \cite{Lurie}
and this notion serves as a complete joint unification of 
the geometry over the given site with homotopy theory. 

For instance, over the trivial site an $\infty$-stack is the same
as an $\infty$-groupoid (a Kan complex), which in turn is equivalently
-- in the sense of homotopy theory -- a topological space. We say that
the collection of all of these
$$
  \mathrm{Top} \simeq \infty \mathrm{Grpd}
$$
is the canonical base $\infty$-topos. 
Next, in generalization to
the relation between Lie groupoids and differentiable stacks, we speak
of the $\infty$-stacks on the site of smooth manifolds
as being \emph{smooth $\infty$-groupoids}.
We write $\mathbf{H} := \mathrm{Smooth}\infty \mathrm{Grpd}$
for the $\infty$-topos of smooth $\infty$-groupoids. 
For emphasis, we will say that an object in $\infty \mathrm{Grpd}$
is a \emph{bare} or \emph{discrete} $\infty$-groupoid
(not equipped with nontrivial smooth structure). Notice that therefore
also topological spaces are, in the sense of homotopy theory, 
\emph{discrete} $\infty$-groupoids. 

It turns out that the 
$\infty$-topos $\mathrm{Smooth}\infty \mathrm{Grpd}$ 
sits over that of bare $\infty$-groupoids
by an adjoint quadruple of $\infty$-functors
$$
  \xymatrix{
    \mathrm{Smooth \infty \mathrm{Grp}}
     \ar@<+12pt>[rr]^{\Pi}
     \ar@{<-}@<+4pt>[rr]|{\mathrm{Disc}}
     \ar@{->}@<-4pt>[rr]|{\Gamma}
     \ar@{<-}@<-12pt>[rr]_{\mathrm{coDisc}}     
    &&
    \infty\mathrm{Grpd}
  }
$$
and this controls the notion of smooth refinement of 
bare cohomology. The geometric interpretation is this:
\begin{itemize}
  \item ``$\mathrm{Disc}$" produces smooth $\infty$-groupoids with
    discrete smooth structure;
  \item ``$\Gamma$" forgets the smooth structure on a smooth $\infty$-groupoid;
  \item ``$\Pi$" sends a smooth $\infty$-groupoid $X$ to its 
   \emph{fundamental path $\infty$-groupoid}; combined with the
   equivalence $\infty \mathrm{Grpd} \simeq \mathrm{Top}$ this is
   the operation of \emph{geometric realization}.
\end{itemize}
For any object or diagram in $\infty \mathrm{Grpd}$
by a \emph{smooth lift} of it we mean a lift through $\Pi$ to an object
or diagram, respectively, in $\mathrm{Smooth}\infty\mathrm{Grpd}$.

\vspace{3mm}
A smooth $\infty$-groupoid with a single object 
we write $\mathbf{B}G$, where $G$ is the 
\emph{smooth $\infty$-group} of automorphisms of that single object.
The boldface $\mathbf{B}$ denotes delooping in 
$\mathrm{Smooth}\infty \mathrm{Grpd}$ as opposed to 
in $\infty \mathrm{Grpd}$.
This completely defines both structures: pointed connected 
smooth $\infty$-groupoids are equivalent to smooth $\infty$-groups.
For instance there is for each $n \in \mathbb{N}$ a smooth $\infty$-groupoid
$\mathbf{B}^n U(1)$, defined this way inductively from the ordinary
smooth circle group $U(1) = \mathbf{B}^0 U(1)$. This is such that
$\Pi \mathbf{B}^n U(1) \simeq B^{n+1} \mathbb{Z} \simeq K(\mathbb{Z}, n+1)$
is the Eilenberg-MacLane space that classifies integral cohomology in 
degree $n$. Hence $\mathbf{B}^n U(1)$ is a smooth refinement of 
the classifying space $B^n U(1) \simeq K(\mathbb{Z},n+1)$.

A morphism $X \to \mathbf{B}^n U(1)$ of smooth $\infty$-groupoids
classifies a \emph{circle $n$-bundle} on $X$. For $n = 1$ and 
$X$ an ordinary manifold this are ordinary circle bundle, for 
$n = 2$ this are bundle gerbes, for $n = 3$ this are bundle 2-gerbes, etc. 
Generally, a morphism
$X \to \mathbf{B}G$ classifies a 
\emph{$G$-principal $\infty$-bundle}.

Note that $X$ here can be much more general than a smooth manifold. 
Notably, we can have $X = \mathbf{B}G$
for $G$ a Lie group or more general smooth $\infty$-group. A morphisms
$\mathbf{c} : \mathbf{B}G \to \mathbf{B}^n U(1)$ then defines a 
cocycle in generalized (Segal-Brylinski-)smooth group cohomology
on $G$ with coefficients in $U(1)$ in degree $n$. If
$G$ is a Lie group or an $\infty$-group presented by a 
simplicial Lie group, then $\Pi \mathbf{B}G \simeq B G$
is the ordinary classifying space of $G$. Therefore such a cocycle
maps under $\Pi$ to an ordinary integral cocycle 
\(
  [\Pi(\mathbf{B}G \stackrel{\mathbf{c}}{\to} \mathbf{B}^n U(1))]
  \simeq
  [B G \stackrel{c}{\to} K(\mathbb{Z}, n+1)]
  \in 
   H^{n+1}(B G, \mathbb{Z})
\)
on the classifying space. 
So the cocycle $\mathbf{c}$ is a smooth refinement of the topological
characteristic map $c$. We say that $\mathbf{c}$ lives in the \emph{smooth cohomology}
$H^n_{\mathrm{smooth}}(\mathbf{B}G, U(1))$.

Given any such cocycle $\mathbf{c}$, we say that the 
\emph{$\mathbf{c}$-twisted cohomology} on $X$ is the connected components of the
homotopy pullback
\(
  \raisebox{20pt}{
  \xymatrix{
    \mathbf{c}\mathrm{Struc}_{\mathrm{tw}}(X)
    \ar[d]
    \ar[r]^{\hspace{-4mm}\mathrm{tw}}_<{\hpull}
    & 
    H^n_{\mathrm{smooth}}(X,U(1))
    \ar[d]
    \\
    \mathbf{H}(X, \mathbf{B}G)
    \ar[r]^{\hspace{-3mm}\mathbf{c}}
    &
    \mathbf{H}(X, \mathbf{B}^n U(1))
  }
  }
  \,,
\)
where the right vertical morphism is the canonical
effective epimorphism that picks one cocycle in each
cohomology class. The map $\mathrm{tw}$
here sends twisted cocycles to their twist.
For instance for $\mathbf{c} : \mathbf{B}PU(n) \to \mathbf{B}^2 U(1)$
the cocycle that classifies the extension 
$\mathbf{B}U(1) \to \mathbf{B} U(n) \to \mathbf{B}PU(n)$ we have that
$\mathbf{c}\mathrm{Struc}_{\mathrm{tw}}(X)$ is the groupoid of 
twisted complex vector bundles of rank $n$ on $X$, those that appear
in the geometric model for twisted K-theory in degree 0.

The reader may be more familiar with twisted cohomology formulated in terms
of sections of certain bundles. We briefly indicate how this is equivalently
another perspective on the above setup. 

Consider first the example of a Lie group $G$ acting on a vector space $V$.
The weak quotient of this action is a Lie groupoid $V//G$ whose objects
form the space $V$ and where there is precisely one morphism for every
ordered pair of points related by the group action. This Lie groupoid
is equipped with a cononical projection $\mathbf{\rho} : V//G \to \mathbf{B}G$.
We may think of this as the smooth incarnation of the vector bundle
that is associated via $\rho$ to the universal $G$-bundle over $\mathbf{B}G$.
In fact we have a fiber sequence $V \to V//G \stackrel{\mathbf{\rho}}{\to}
\mathbf{B}G$ of Lie groupoids, 
and this may equivalently be taken to define the action 
$\rho$ of $G$ on $V$.

Now consider a morphism $g : X \to \mathbf{B}G$ classifying a $G$-principal
bundle $P \to X$, as above. By inspection one finds that a lift $\sigma$ of 
this morphism along this projection
\(
  \raisebox{20pt}{
  \xymatrix{
    & V // G
	\ar[d]^{\mathbf{\rho}}
    \\
    X \ar[r]^g \ar[ur]^\sigma & \mathbf{B}G
  }
  }
\)
is precisely a section of the vector bundle $P \times_\rho V$ that is 
associated to $P$ by the given representation. On the other hand,
in terms of the above notion of twisted cohomology, such a lift 
is also precisley an element in the \emph{$\mathbf{\rho}$-twisted cohomology}
of $X$ with coefficients in $V//G$, where the twist is the class of $P$:
we have an equivalence
\(
  \mathbf{\rho}\mathrm{Struc}_{[P]}(X)
  \simeq
  \Gamma_X(P \times_\rho V)
\)
of cocycles in $\mathbf{\rho}$-twisted cohomology with sections of the
$\rho$-associated vector bundle.

This perspective generalizes verbatim to all twisting
cocycles $\mathbf{c}$ on all $\infty$-groups. 
For those of the form $\mathbf{c} : \mathbf{B}G \to \mathbf{B}^{n}U(1)$
considered before, we may think of their homotopy fiber $\mathbf{B}\hat G$
in the fiber sequence
\(
  \label{FiberSequenceAndRepresentation}
  \mathbf{B}\hat G \to \mathbf{B}G \stackrel{{\mathbf{c}}}{\to}
  \mathbf{B}^n U(1)
\)
as the delooping of the shifted central extension $\hat G$ of $G$
classified by the cocycle $\mathbf{c}$, or equivalently think of
$\mathbf{B}G$ as a universal 
$\mathbf{c}$-associated $\mathbf{B}\hat G$-bundle over
$\mathbf{B}^n U(1)$. Then for
$P \to X$ classified by $g : X \to \mathbf{B}^n U(1)$ a given
circle $n$-bundle on $X$, the $\infty$-groupoid 
$\mathbf{c}\mathrm{Struc}_{[g]}$ may be thought of equivalently
as the space of sections of the associated $\mathbf{B}\hat G$-bundle.
If $P$ is trivial, such sections are just maps from $X$ to $\mathbf{B}\hat G$.

\vspace{3mm}
In order to equip such structures of smooth and twisted 
cohomology with connections,
we reflect $\Pi$ back to smooth $\infty$-groupoids by
defining $\mathbf{\Pi} := \mathrm{Disc}\; \Pi  : 
\mathrm{Smooth}\infty \mathrm{Grpd} \to \mathrm{Smooth}\infty \mathrm{Grpd}$.
For $X$ a smooth $\infty$-groupoid we say that $\mathbf{\Pi}(X)$
is its \emph{smooth path $\infty$-groupoid}. 
The adjunction unit gives a canonical morphism
$X \to \mathbf{\Pi}(X)$ which includes $X$ as the 
\emph{constant} paths in $X$.
A morphism
$\mathbf{\Pi}(X) \to \mathbf{B}G$ therefore has an 
underlying $G$-principal $\infty$-bundle $X \to \mathbf{\Pi}(X) \to \mathbf{B}G$
but also assigns equivalences between fibers over endpoints of paths,
that are equivalent around disks, etc. This is 
the parallel transport of a 
\emph{flat $\infty$-connection} on the smooth $G$-principal bundle.

From this in turn derives a notion of $G$-valued flat
differential forms: these are the flat $G$-connections whose underlying
$\infty$-bundle is trivial. We write
\(
  \mathbf{\Pi}_{\mathrm{dR}}X := \mathbf{\Pi}X \coprod_{X} * 
\)
for the canonical homotopy pushout and speak 
of the \emph{de Rham homotopy type} of $X$. A morphism
$\mathbf{\Pi}_{\mathrm{dR}}(X) \to \mathbf{B}G$ is a closed
$G$-valued form on $X$. In particular 
a morphism $\mathbf{\Pi}_{\mathrm{dR}}(X) \to \mathbf{B}^n U(1)$
is a closed $n$-form.

A key observation now is that there is a canonically induced
morphism of cocycle $\infty$-groupoids
\(
  \mathrm{curv} : \mathbf{H}(X,\mathbf{B}^n U(1))
  \to 
  \mathbf{H}(\mathbf{\Pi}_{\mathrm{dR}}(X), \mathbf{B}^{n+1}U(1))
\)
that sends each circle $n$-bundle to a \emph{curvature}
characteristic form. We define ordinary
\emph{differential cohomology} 
$\mathbf{H}_{\mathrm{diff}}(X, \mathbf{B}^n U(1))$
of $X$ to be the 
$\infty$-pullback of the canonical effective epimorphism
$H_{\mathrm{dR}}^n(X) \to
\mathrm{Smooth}\infty\mathrm{Grpd}(\mathbf{\Pi}_{\mathrm{dR}}(X), \mathbf{B}^{n+1}U(1))$ (that which picks one cocycle in each cohomology class) 
along this curvature characteristic map.
Cocyles in this homotopy pullback 
\(
  \raisebox{20pt}{
  \xymatrix{
    \mathbf{H}_{\mathrm{diff}, \mathbf{B}^n U(1)}
	\ar[r]_<{\hpull}
	\ar[d]
	&
	H^{n+1}_{\mathrm{dR}}(X)
	\ar[d]
	\\
	\mathbf{H}(X, \mathbf{B}^n U(1))
	\ar[r]
	&
	\mathbf{H}(\mathbf{\Pi}_{\mathrm{dR}} X, \mathbf{B}^{n+1} U(1))
  }
  }
\)
are smooth circle $n$-bundles
\emph{with connection}.

For $X$ an ordinary manifold, this reproduces the ordinary
notions in differential cohomology. More precisely, 
in this case the $n$-groupoid $\mathbf{H}_{\mathrm{diff}}(X, \mathbf{B}^n U(1))$
turns out to be that whose objects are cocycles 
in Deligne-Beilinson hypercohomology, whose morphisms
are smooth gauge transformations of these, whose 2-morphisms
are higher gauge transformations of those, and so on. 
But as before, we can 
apply this over any smooth $\infty$-groupoid.
In particular we may consider differential cohomology
over moduli stacks $\mathbf{B}G$ of $G$-principal 
$\infty$-bundles that differentially refine smooth
lifts $\mathbf{c} : \mathbf{B}G \to \mathbf{B}^{n+1} U(1)$
of characteristic map. These are \emph{Chern-Simons $n$-gerbes
with connection}.

This finally gives rise to the notion of $\infty$-connections
on general nonabelian $G$-principal $\infty$-bundles: these
are structures that lift  
universal curvature classes from de Rham cohomology to 
differential cohomology.

With all these concepts thus abstractly given, we
can look for explicit constructions of these. 
By standard theory \cite{Lurie} 
every smooth $\infty$-groupoid $A$ is presented 
a simplicial presheaf on the category of Cartesian spaces
and smooth maps: a functor
\(
  A : (U =\mathbb{R}^n , [k]) \mapsto A_k(U) \in \mathrm{Set}
\)
that we read as assigning to each test space $U$ 
and each $k \in \mathbb{N}$ the 
set of $k$-morphisms of the 
$\infty$-groupoid of possible ways of \emph{probing} $A$ with $U$, 
or equivalently the set of $U$-parameterized smooth families of 
$k$-morphisms in $A$.

We consider the construction of such simplicial
presheaves from infinitesimal data. For $\mathfrak{g}$ an 
an $L_\infty$-algebra
(see \cite{SSS1} and the appendix for a review or relevant notions) 
we can define such a presheaf by
setting, in evident generalization of the construction in \cite{H},
\(
  \exp(\mathfrak{g}) : (U,[k]) \mapsto
  \left\{
    \xymatrix{
      \Omega^\bullet(U \times \Delta^k)_{\mathrm{vert}}
      \ar@{<-}[r]^<<<{ \hspace{4mm}A_{\mathrm{vert}}}
      &
      \mathrm{CE}(\mathfrak{g})
    }
  \right\}
  \,,
\)
where on the right we have the set of dg-algebra homomorphisms
from the \emph{Chevalley-Eilenberg algebra} of $\mathfrak{g}$
to the de Rham complex of vertical forms on the trivial simplex 
bundle $U \times \Delta^k \to U$. The Chevalley-Eilenberg dg-algebra
$\mathrm{CE}(\mathfrak{g})$ is the free graded-commutative algebra
on the the degreewise dual of the graded vector space underlying 
$\mathfrak{g}$ and equipped with the differential obtained by 
dualizing all the brackets on $\mathfrak{g}$. For $\mathfrak{g}$
an ordinary Lie algebra it reduces to the ordinary Chevalley-Eilenberg
algebra, hence its name.

This 
simplicial presheaf $\exp(\mathfrak{g})$ 
presents the smooth $\infty$-groupoid $\mathbf{B}G$ for 
$G$ the ``universal $\infty$-connected'' Lie integration of 
$\mathfrak{g}$. Smooth Postnikov truncations $\tau_n$ of this object 
yield smooth $n$-groups
integrating $\mathfrak{g}$. For instance for $\mathfrak{g}$
an ordinary Lie algebra we have that $\tau_1 \exp(\mathfrak{g})
\simeq \mathbf{B}G$ for $G$ the ordinary simply connected Lie group 
integrating $\mathfrak{g}$.

The crucial point for our discussion here is that
one obtains from this also a model for the 
smooth $\infty$-groupoid that classifies $G$-principal $\infty$-connections
for such Lie integrated $G$ \cite{FSSI}. 
If we write $\mathrm{W}(\mathfrak{g})$ for the \emph{Weil algebra}
of $\mathfrak{g}$, the unique free dg-algebra on the dual graded 
vector space of $\mathfrak{g}$ such that the canonical projection 
is a dg-algebra homomorphism
to $\mathrm{CE}(\mathfrak{g})$, and $\mathrm{inv}(\mathfrak{g})$ for
the subalgebra of closed elements formed out of shifted generators
-- the \emph{invariant polynomials} on $\mathfrak{g}$ --,
then this is 
given by the simplicial presheaf defined by
\(
  \label{expgconn}
 \exp(\mathfrak{g})_{\mathrm{conn}}
  :
  (U,[k])
  \mapsto
  \left\{
    \raisebox{40pt}{
    \xymatrix{
      \Omega^\bullet_{\mathrm{vert}}(U \times \Delta^k)
      \ar@{<-}[r]^<<<{\hspace{4mm}A_{\mathrm{vert}}}
      &
      \mathrm{CE}(\mathfrak{g})
      \\
      \Omega^\bullet(U \times \Delta^k)
      \ar@{<-}[r]^<<<{\hspace{8mm}(A,F_A)}
      \ar[u]
      &
      \mathrm{W}(\mathfrak{g})
      \ar[u]
      \\
      \Omega^\bullet(U)
      \ar@{<-}[r]^<<<{\hspace{13mm}\langle F_A\rangle}
      \ar[u]
      &
      \mathrm{inv}(\mathfrak{g})
      \ar[u]
    }
  }
  \right\}
  \,.
\) 
On the right we have the set of horizontal dg-algebra homomorphisms
that makes the diagram commute, see around def. \ref{expgconngroupoid}
in the appendix for more details.
This is the $L_\infty$-algebraic differential form
data discussed in detail in \cite{SSS1}, here parameterized 
over all test spaces $U$ and simplices $\Delta^k$ as discussed in 
\cite{FSSI}. The 
horizontal morphism in the middle constitutes $\mathfrak{g}$-valued
differential form data on $U \times \Delta^k$ and the fact that 
it sits in this commuting diagram encodes the $\infty$-analogs of the
two conditions of an ordinary Cartan-Ehresmann connection: the top square 
says that the vertical part of $A$ is flat, and the bottom square says
that the curvature forms ``transform covariantly'' and make all
invariant polynomials descent down to $U$.

For $X$ a smooth manifold a morphism
$X \to \exp(\mathfrak{g})_{\mathrm{conn}}$ is equivalently
\begin{enumerate}
  \item a choice of good open cover 
  $\{U_i \to X\}$;
  \item
    on each patch $U_i$, differential form data with 
    values in $\mathfrak{g}$;
  \item
    on each double intersection, a choice of 1-parameter
    gauge transformation between the corresponding differential form data;
  \item
    on each triple intersection, a choice of 2-parameter gauge-of-gauge
    transformation;
   \item and so on for higher intersections.
\end{enumerate}
Given then a differential refinement of a cocycle 
$\mathbf{c} : \mathbf{B}G \to \mathbf{B}^n U(1)$ to 
the corresponding moduli stacks of bundles with connection
\(
  \hat {\mathbf{c}}
  :
  \mathbf{B}G_{\mathrm{conn}}
  \to 
  \mathbf{B}^n U(1)_{\mathrm{conn}}\;,
\)
we can consider the differential refinement of the twisted
smooth cohomology discussed above:
the cocycle $\infty$-groupoid of
$\hat {\mathbf{c}}$-\emph{twisted differential cohomology}
is the homotopy pullback
\(
  \raisebox{20pt}{
  \xymatrix{
    \hat {\mathbf{c}}\mathrm{Struc}_{\mathrm{tw}}(X)
    \ar[d]
    \ar[r]^{\mathrm{tw}}_<{\hpull}
    & 
    H^{n+1}_{\mathrm{diff}}(X)
    \ar[d]
    \\
    \mathbf{H}(X, \mathbf{B}G_{\mathrm{conn}})
    \ar[r]^{ \hspace{-3mm}{\hat{\mathbf{c}}}}
    &
    \mathbf{H}(X, \mathbf{B}^n U(1)_{\mathrm{conn}})
  }
  }
  \,,
\)
where again the right vertical morphism is the canonical effective
epimorphism that picks one cocycle in each cohomology class.
Up to some slight technicalities which are discussed in
\cite{FSSI}, this homotopy pullback is modeled by the corresponding pullback
of the double square diagrams of  $L_\infty$-algebra data from above, which 
are  discussed in detail in the last part of \cite{SSS1}. This we will later 
unwind 
in sections \ref{G4 twist}
and \ref{G8 twist} for the case of twisted differential String- and
Fivebrane structures.

\begin{facts}
\label{fac}
  \begin{enumerate}
\item  There is, up to equivalence, a unique smooth refinement
  of the first fractional Pontrjagin class (from fact \ref{FirstFractionalPontrjagin})
  \(
    \frac{1}{2}\mathbf{p}_1 : \mathbf{B}\mathrm{Spin} \to 
    \mathbf{B}^3 U(1)
    \,.
   \label{thm fss}
  \)
  to the moduli stack of smooth $\mathrm{Spin}$-principal bundles
  with values in the moduli 3-stack of circle 3-bundle.
 \item 
  The smooth string 2-group, whose smooth delooping we denote
  $\mathbf{B}\mathrm{String}$, is the homotopy fiber of 
  $\frac{1}{2}\mathbf{p}_1$, sitting in a fiber sequence
  \(
    \mathbf{B}\mathrm{String} \to 
    \mathbf{B}\mathrm{Spin}
    \stackrel{\frac{1}{2}\mathbf{p}_1}{\longrightarrow}
    \mathbf{B}^3 U(1)
    \,.
    \label{SmoothString2Group}
  \)
\item  The models for the smooth string 2-group given in 
  \cite{H} and in \cite{BCSS} are indeed presentations of
  the abstract definition \eqref{SmoothString2Group}.
 \item There is a smooth refinement 
  $
    \frac{1}{6} \mathbf{p}_2 : \mathbf{B}\mathrm{String} \to 
    \mathbf{B}^7 U(1)
  $
  of the second fractional Pontrjagin class from 
  fact \ref{SecondFractionalPontrjagin} to the moduli 2-stack 
  of $\mathrm{String}$-principal 2-bundles with values in the
  moduli 7-stack of circle 7-bundles.
 \item  The homotopy fiber of $\frac{1}{2}\mathbf{p}_2$ is the
  smooth delooping of the \emph{smooth fivebrane 6-group}
  $$
    \mathbf{B}\mathrm{Fivebrane}
    \to 
    \mathbf{B}\mathrm{String}
    \stackrel{\frac{1}{6}\mathbf{p}_2}{\longrightarrow}
    \mathbf{B}^7 U(1)
    \,.
  $$

\item   Under geometric realization these smooth lifts indeed 
  reproduce the first steps in the Whitehead of $O$:
  $$
    \vert \mathbf{B}\mathrm{String}\vert \simeq 
    B \mathrm{String} \qquad
{\rm  and} \qquad
    \vert \mathbf{B}\mathrm{Fivebrane}\vert \simeq 
    B \mathrm{Fivebrane}
    \,.
  $$
\end{enumerate}
\end{facts}
This is shown in \cite{FSSI} and also discussed in section 4.1 of \cite{nactwist}.

In summary all this means that we obtain the following
canonical refinement of def. \ref{TwistedStringStructures}
and def. \ref{TwistedFivebraneStructures} to a notion of 
\emph{twisted differential string- and fivebrane
structures}.
\begin{definition}
  \label{TwistedDifferentialStructures}
  For $X$ a smooth manifold, the 2-groupoid of
  \emph{twisted differential String-structures} 
  $\frac{1}{2}\hat {\mathbf{p}}_1\mathrm{Struc}_{\mathrm{tw}}(X)$ on $X$ is the homotopy 
  pullback
  \(
    \raisebox{20pt}{
    \xymatrix{
      \frac{1}{2}\hat {\mathbf{p}}_1\mathrm{Struc}_{\mathrm{tw}}(X)      
      \ar[d]
      \ar[r]^{\mathrm{tw}}_<{\hpull}
      &
      H^4_{\mathrm{diff}}(X)
      \ar[d]
      \\
      \mathbf{H}(X, \mathbf{B}\mathrm{Spin}_{\mathrm{conn}})
      \ar[r]^{\frac{1}{2}\hat {\mathbf{p}}_1}
      &
      \mathbf{H}(X, \mathbf{B}^3 U(1)_{\mathrm{conn}})
    }
	}
    \,.
  \) 
  Analogously, the 6-groupoid of
  \emph{twisted differential fivebrane-structures} 
  $\frac{1}{6}\hat {\mathbf{p}}_2\mathrm{Struc}_{\mathrm{tw}}(X)$ 
  on $X$ is the homotopy 
  pullback
  \(
    \raisebox{20pt}{
    \xymatrix{
      \frac{1}{6}\hat {\mathbf{p}}_2\mathrm{Struc}_{\mathrm{tw}}(X)      
      \ar[d]
      \ar[r]^{\mathrm{tw}}_<{\hpull}
      &
      H^8_{\mathrm{diff}}(X)
      \ar[d]
      \\
      \mathbf{H}(X, \mathbf{B}\mathrm{String}_{\mathrm{conn}})
      \ar[r]^{\frac{1}{6}\hat {\mathbf{p}}_2}
      &
      \mathbf{H}(X, \mathbf{B}^7 U(1)_{\mathrm{conn}})
    }
	}
    \,.
  \)
\end{definition}
Notice that since these constructions had been announced in \cite{SSS1}
the article \cite{Waldorf} has appeared which defines 
a presentation of twisted differental String-structures
in terms of bundle 2-gerbes for the case that the underlying
topological twist vanishes.

We now describe the local $L_\infty$-algebraic differential form data 
of such twisted differential structures, following \cite{SSS1}.
Details of the following fact are in \cite{FSSI}.

\begin{theorem}
(i) The local differential form data of a 
twisted $\mathrm{String}(n)$-bundle with connection is that 
known from the Green-Schwarz mechanism (section \ref{HeteroticStringAndTwistedStringStructures}).

\noindent (ii) The local differential form data of a 
twisted $\mathrm{Fivebrane}(n)$-bundle with connection is that of the dual 
Green-Schwarz mechanism (section \ref{NS5 and twisted Fivebrane Structures}).

\end{theorem}

\subsection{Twisted $\mathfrak{string}(n)$ 2-connections} 
\label{G4 twist}

One finds that the smooth first fractional Pontryagin map
\(
  \frac{1}{2}p_1 : \mathbf{B}\mathrm{Spin}(n) \to \mathbf{B}^3(U(1))
\)
from the moduli stack of $\mathrm{Spin}$-principal bndles
to the moduli 3-stack of circle 3-bundles
may be modeled in terms of simplicial presheaves 
by the span of smooth 2-groupoids 
of the form
\(
  \mathbf{B}\mathrm{Spin}
  \stackrel{\simeq}{\leftarrow}
  \mathbf{B}(\mathbf{B}U(1) \to \mathrm{String})
  \to 
  \mathbf{B}^3 U(1)
\)
given by a span of smooth crossed complexes of the form
\(
  (
    \xymatrix{
   1 \ar[r] &  1 \ar[r] & \mathrm{Spin}(n) 
   })
  \stackrel{\simeq}{\longleftarrow}
  (
   \xymatrix{
   U(1) \ar[r] &  
   \hat \Omega \mathrm{Spin}(n) 
   \ar[r] &  
   P \mathrm{Spin}(n) 
   }
   )
  \stackrel{\frac{1}{2}p_1}{\longrightarrow}
  (
    \xymatrix{
    U(1) \ar[r] &  1 \ar[r] &  1 
	}
   )
  \,,
\)
where $\hat \Omega \mathrm{Spin}$ denotes the Kac-Moody central extension
of the loop group of $\mathrm{Spin}$.

The corresponding span of $L_\infty$-algebras is of the form
\(
  \mathfrak{so}(n)
  \stackrel{\simeq}{\longleftarrow}
  (b \mathfrak{u}(1) \to \mathfrak{string}(n))
  \to
  b^2 u(1)
  \,,
\)
where in the middle we we have the mapping cone of the inclusion 
of the line Lie 2-algebra, example \ref{LineLienAlgebra},
into the string Lie 2-algebra, example \ref{StringLie2Algebra}.

Therefore for $C(\{U_i\}) \stackrel{\simeq}{\to} X$ the {\v C}ech
nerve projection out of a sufficiently good open cover of $X$
and for $X \stackrel{\simeq}{\leftarrow} C(\{U_i\}) \stackrel{\mathrm{tw}}{\to} \mathbf{B}^3 U(1)$ a cocycle for the twisting circle 3-bundle, 
a corresponding twisted $\mathrm{String}$-2-bundle is a lift
$\hat g$ in the diagram of simplicial presheaves
$$
  \xymatrix{
     & \mathbf{B}(\mathbf{B}U(1) \to \mathrm{String})
	 \ar[dr]
     \\
     C(\{U_i\})
	 \ar[d]^{\simeq}
	 \ar[rr]^{\mathrm{tw}}
	 \ar[ur]^{\hat g}
	 &&
	 \mathbf{B}^3 U(1)
	 \\
	 X
  }
  \,.
$$
After passing to the differential refinement of this
situation, the corresponding Cartan-Ehresmann $L_\infty$-connection is
given on $(k+1)$-fold intersections of the cover by compatible diagrams
\(
  \xymatrix{
    \Omega^\bullet_{vert}(U \times \Delta^k) 
     \ar@{<-}[r]^<<<<<<{A_{\mathrm{vert}}}
     \ar@{<-}[d]
    & \mathrm{CE}(b \mathfrak{u}(1) \to \mathfrak{string}(n))
     \ar@{<-}[d]
    \\
    \Omega^\bullet( U \times \Delta^k) 
    \ar@{<-}[r]^<<<<<<<<{ \hspace{-1mm}(A,F_A)}
    \ar@{<-}[d]
    & 
    \mathrm{W}(b \mathfrak{u}(1) \to \mathfrak{string}(n))    
    \ar@{<-}[d]
    \\
    \Omega^\bullet(U) \ar@{<-}[r]^<<<<<<<<{P(F_A)}
    & \mathrm{inv}(b \mathfrak{u}(1) \to \mathfrak{string}(n))    
  }
\)
such that the corresponding twist is the prescribed one. 
If we write $(C_3)$ for the local differential 3-form data 
of the prescribed twisting circle 3-bundle, then this 
means that over $(k+1)$-fold intersections the $L_\infty$-algebraic
data of a connection on the twisted $\mathrm{String}$-bundle
is given by a diagram of dg-algebras of the form
$$
  \xymatrix@C=30pt@R=40pt{
    &
    \mathrm{CE}(b \uu(1) \hookrightarrow\gg_\mu)
    \ar@{-->}[dl]|{\{A_{\mathrm{vert}}, B_{\mathrm{vert}}, (C_3)_{\mathrm{vert}}\}}
    \\
    \Omega^\bullet_{\mathrm{vert}}(U \times \Delta^k)
    &&
    \mathrm{CE}(b^2 \uu(1))
    \ar[ll]|>>>>>>>>>>>>>>{(C_3)_{\mathrm{vert}}}
    \ar@{_{(}->}[ul]
    \\
    &
    \mathrm{W}(b \uu(1) \hookrightarrow \gg_{\mu})
    \ar@{-->}[dl]|{\{A,F_A, B, \nabla B, C_3, G_4\}}
    \ar@{->>}[uu]|{\makebox(20,20){}}
    \\
    \Omega^\bullet(U \times \Delta^k)
    \ar@{->>}[uu]
    &&
    \mathrm{W}(b^2 \uu(1))
    \ar[ll]|>>>>>>>>>>>>>>>>{(C_3, G_4)}
    \ar@{_{(}->}[ul]
    \ar@{->>}[uu]
    \\
    &
    \mathrm{inv}(\mathrm{inn}(b \uu(1)) \hookrightarrow\mathrm{cs}_P(\gg))
    \ar@{-->}[dl]|{\{H_3,G_4,P(F_A)\}}
    \ar@{->>}[uu]|{\makebox(20,20){}}
    \\
    \Omega^\bullet(U)
    \ar@{^{(}->}[uu]^{p^*}
    &&
    \mathrm{inv}(b^2 \uu(1))
    \ar[ll]_>>>>>>>>>>>>>>>{G_4}
    \ar@{_{(}->}[ul]
    \ar@{^{(}->}[uu]
  }
  \,,
$$
where, for short, we write $\mathfrak{g}_\mu$ for $\mathfrak{string}$.
It may be helpful to think of this as forming sections of a higher
associated bundle, as in def. \ref{SectionsAndCovariantDerivatives} in 
the appendix. In terms of this we may read this diagram as
indicated in the following
$$
  
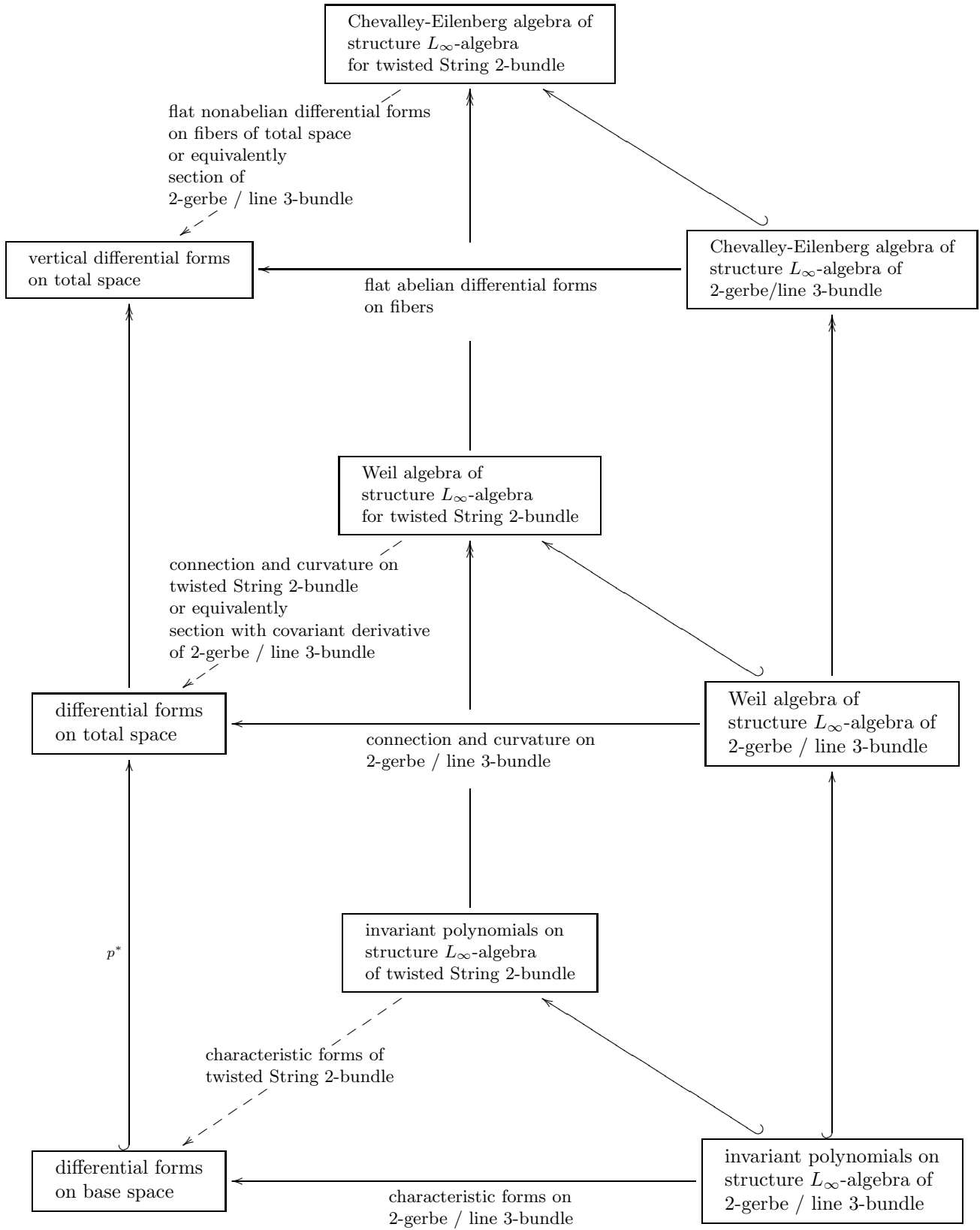
\hspace{-1cm}
  \xymatrix@C=30pt@R=68pt{
    &
    \fbox{
      \small
      \begin{tabular}{l}
        Chevalley-Eilenberg algebra of
        \\
        structure $L_\infty$-algebra
        \\
        for twisted String 2-bundle
      \end{tabular}
    }
    \ar@{-->}[dl]|{
      \mbox{
        \small
        \begin{tabular}{l}
          flat nonabelian differential forms 
          \\
          on fibers of total space
          \\
          or equivalently
          \\
          section of
          \\
          2-gerbe / line 3-bundle
        \end{tabular}
      }
    }
    \\
    \fbox{
      \small
      \begin{tabular}{l}
        vertical differential forms 
        \\
        on total space
      \end{tabular}
    }
    &&
    \fbox{
      \small
      \begin{tabular}{l}
        Chevalley-Eilenberg algebra of
        \\
        structure $L_\infty$-algebra of
        \\
        2-gerbe/line 3-bundle
      \end{tabular}
    }
    \ar[ll]^{
      \mbox{
        \small
        \begin{tabular}{l}
          flat abelian differential forms
          \\
          on fibers
        \end{tabular}
      }
    }
    \ar@{_{(}->}[ul]
    \\
    &
    \fbox{
      \small
      \begin{tabular}{l}
        Weil algebra of
        \\
        structure $L_\infty$-algebra
        \\
        for twisted String 2-bundle
      \end{tabular}
    }
    \ar@{-->}[dl]|{
      \mbox{
        \small
        \begin{tabular}{l}
          connection and curvature on
          \\
          twisted String 2-bundle
          \\
          or equivalently
          \\
          section with covariant derivative
          \\
          of 2-gerbe / line 3-bundle
        \end{tabular}
      }
    }
    \ar@{->>}[uu]|>>>>>>>>>>>>>>>>>>>>>>>>>>>>>>>>>>>{\makebox(20,50){}}
    \\
    \fbox{
      \begin{tabular}{l}
        differential forms \\
        on total space
      \end{tabular}
    }
    \ar@{->>}[uu]
    &&
    \fbox{
      \begin{tabular}{l}
        Weil algebra of
        \\
        structure $L_\infty$-algebra of
        \\
        2-gerbe / line 3-bundle
      \end{tabular}
    }
    \ar[ll]^{
      \mbox{
        \small
        \begin{tabular}{l}
          connection and curvature on
          \\
          2-gerbe / line 3-bundle
        \end{tabular}
      }
    }
    \ar@{_{(}->}[ul]
    \ar@{->>}[uu]
    \\
    &
    \fbox{
      \small
      \begin{tabular}{l}
        invariant polynomials on
        \\
        structure $L_\infty$-algebra
        \\
        of twisted String 2-bundle
      \end{tabular}
    }
    \ar@{-->}[dl]|{
      \mbox{
        \small
        \begin{tabular}{l}
          characteristic forms of
          \\
          twisted String 2-bundle
        \end{tabular}
      }
    }
    \ar@{->>}[uu]|>>>>>>>>>>>>>>>>>>>>>>>>>>>>>>>>>>>>{\makebox(20,40){}}
    \\
    \fbox{
      \begin{tabular}{l}
        differential forms \\ on base space
      \end{tabular}
    }
    \ar@{^{(}->}[uu]^{p^*}
    &&
    \fbox{
      \begin{tabular}{l}
        invariant polynomials on
        \\
        structure $L_\infty$-algebra of
        \\
        2-gerbe / line 3-bundle
      \end{tabular}
    }
    \ar[ll]^{
      \mbox{
        \small
        \begin{tabular}{l}
          characteristic forms on
          \\
          2-gerbe / line 3-bundle
        \end{tabular}
      }
    }
    \ar@{_{(}->}[ul]
    \ar@{^{(}->}[uu]
  }
$$
Chasing the generators of the graded-commutative algebras
through this diagram and recording the condition imposed by the
respect of the morphisms of dg-algebras for differentials, one finds that
in components the commutativity of this diagram encodes the following
differential form data and the following relations on that.
$$
  
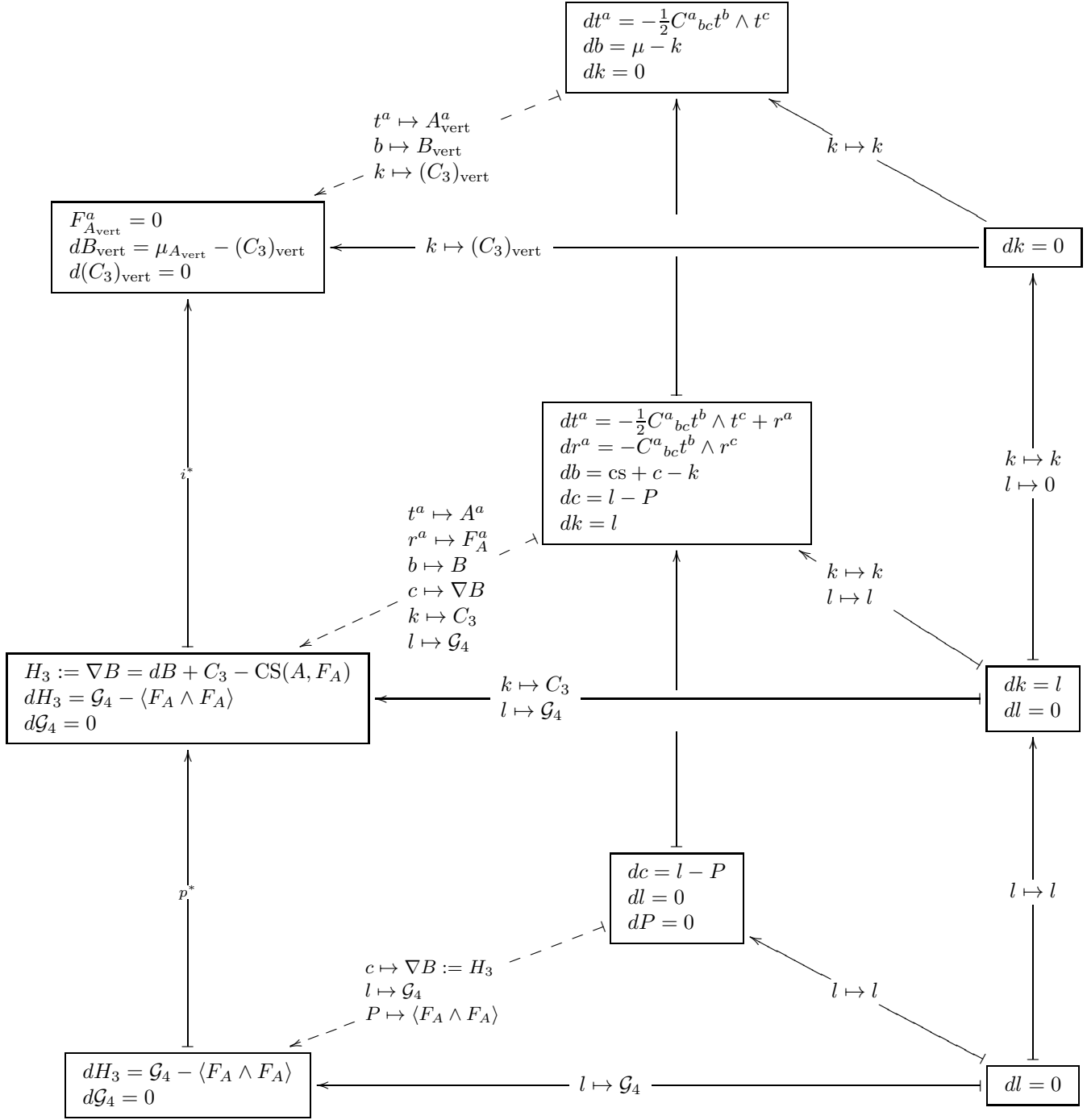
\hspace{-1cm}
  \xymatrix@C=75pt@R=45pt{
    &
    \fbox{$
      \begin{array}{l}
        dt^a = -\frac{1}{2}C^a{}_{bc}t^b \wedge t^c
        \\
        db = \mu  - k
        \\
        dk = 0
      \end{array}
    $}
    \ar@{|-->}[dl]|{
      \begin{array}{l}
        t^a \mapsto A^a_{\mathrm{vert}}
        \\
        b \mapsto B_{\mathrm{vert}}
        \\
        k \mapsto (C_3)_{\mathrm{vert}}
      \end{array}
    }
    \\
    \fbox{$
      \begin{array}{l}
        F_{A_{\mathrm{vert}}}^a = 0
        \\
        d B_{\mathrm{vert}} = \mu_{A_{\mathrm{vert}}} - (C_3)_{\mathrm{vert}}
        \\
        d(C_3)_{\mathrm{vert}} = 0
      \end{array}
    $}
    &
    &
    \fbox{$
      \begin{array}{l}
        dk = 0
      \end{array}
    $}
    \ar@{->}[ll]|>>>>>>>>>>>>>>>>>>>>>>>>>{
      \begin{array}{l}
        k \mapsto (C_3)_{\mathrm{vert}}
      \end{array}
    }
    \ar@{->}[ul]|{
      \begin{array}{l}
        k \mapsto k
      \end{array}
    }
    \\
    &
    \fbox{$
      \begin{array}{l}
        d t^a = - \frac{1}{2}C^{a}{}_{bc}t^b \wedge t^c + r^a
        \\
        d r^a = - C^a{}_{bc}t^b \wedge r^c
        \\
        db = \mathrm{cs} + c - k
        \\
        dc = l - P
        \\
        dk = l
      \end{array}
    $}
    \ar@{|-->}[dl]|<<<<<<<<<<<<<<<<{
      \begin{array}{l}
        t^a \mapsto A^a
        \\
        r^a \mapsto F_A^a
        \\
        b \mapsto B
        \\
        c \mapsto \nabla B
        \\
        k \mapsto C_3
        \\
        l \mapsto {\cal G}_4
      \end{array}
    }
    \ar@{|->}[uu]|<<<<<<<<<<<<<<<<<<<<<<<<<<{
      \makebox(20,20){}
    }
    \\
     \fbox{$
     \begin{array}{l}
       H_3 := \nabla B = dB + C_3 - \mathrm{CS}(A,F_A)
       \\
       d H_3 = {\cal G}_4 - \langle F_A \wedge F_A \rangle
       \\
       d {\cal G}_4 = 0
     \end{array}
    $}
    \ar@{|->}[uu]|{
      i^*
    }
    &
    &
    \fbox{$
     \begin{array}{l}
       dk = l
       \\
       dl = 0
     \end{array}
    $}
    \ar@{|->}[uu]|{
      \begin{array}{l}
        k \mapsto k
        \\
        l \mapsto 0
      \end{array}
    }
    \ar@{|->}[ll]|>>>>>>>>>>>>>>>>>>>>>>>>>>{
      \begin{array}{l}
        k \mapsto C_3
        \\
        l \mapsto {\cal G}_4
      \end{array}
    }
    \ar@{|->}[ul]|{
      \begin{array}{l}
        k \mapsto k
        \\
        l \mapsto l
      \end{array}
    }
  \\
  &
    \fbox{$
     \begin{array}{l}
       dc = l-P
       \\
       dl = 0
       \\
       dP = 0
     \end{array}
    $}
    \ar@{|-->}[dl]|{\small\begin{array}{l}
      c \mapsto \nabla B := H_3
      \\
      l \mapsto {\cal G}_4
      \\
      P \mapsto \langle F_A \wedge F_A\rangle
    \end{array}}
    \ar@{|->}[uu]|>>>>>>>>>>>>>>>>>>>>>>>>>>{\makebox(20,20){}}
  \\
    \fbox{$
     \begin{array}{l}
       d H_3 = {\cal G}_4 - \langle F_A \wedge F_A \rangle
       \\
       d {\cal G}_4 = 0
     \end{array}
     $}
    \ar@{|->}[uu]|{
      p^*
    }
    &&
    \fbox{$
     \begin{array}{l}
       dl = 0
     \end{array}
    $}
    \ar@{|->}[ll]|{
      \begin{array}{c}
        l \mapsto {\cal G}_4
      \end{array}
    }
    \ar@{|->}[ul]|{
      \begin{array}{l}
       l \mapsto l
      \end{array}
    }
    \ar@{|->}[uu]|{
      \begin{array}{l}
       l \mapsto l
      \end{array}
    }
  }
$$

\vspace{3mm}
Here, $P \in W(\gg)$ denotes the invariant polynomial on $\gg$
in transgression with with the cocycle $\mu \in \mathrm{CE}(\gg)$.
With $\{t^a\}$ a fixed chosen basis of $\gg^*$ in degree 1 and
$\{r^a\}$ the corresponding basis in degree 2, we have $P = P_{ab}
r^a \wedge r^b$ and $\mu = \mu_{abc} t^a \wedge t^b \wedge t^c$ and
$\mathrm{cs} =  P_{ab}t^b \wedge r^a + \frac{1}{6} \mu_{abc}t^a
\wedge t^b \wedge t^c$.
We have 

\vspace{3mm}
\begin{center}
\begin{tabular}{|l|l|}
  \hline
  curvature & $H_3 := d B + C_3 - \mathrm{CS}(A,F_A)$
  \\
  \hline
  Bianchi identity & $d H_3 = {\cal G}_4 - \langle F_A \wedge F_A\rangle$
  \\
  \hline
\end{tabular}
\end{center}

\vspace{3mm}
 In \cite{SSS1} this situation was considered from a
different perspective for the special case $B = 0$ and $\nabla B =
0$. There the dashed morphism was obtained as a twisted lift of a
$\gg$-connection to a $\gg_\mu$-connection and the $b^2
\uu(1)$-connection appeared as the corresponding obstruction. Here
now the perspective is switched: the $b^2 \uu(1)$-connection is
prescribed and the choice of dashed morphisms is a choice of twisted
$\gg_\mu$-connections with prescribed twist $G_4$.

\vspace{3mm}
The covariant derivative 3-form $\nabla B$ of the twisted
$\gg_\mu$-connection, which we denote by $H_3,$ measures the
difference between the prescribed $b^2 \uu(1)$-connection and the
twist of the chosen twisted $\gg_\mu$-connection. The
Bianchi identity
\(
  d H_3 = {\cal G}_4 - P(F_A)
\)
which appears in the middle on the left 
says that this difference
has to vanish in cohomology, as one expects. Indeed, this is
the structure of the differential forms in the Green-Schwarz mechanism,
constituting the differential refinement of the 
integral cohomology relation (\ref{Green Schwarz anomaly}).

\subsection{Twisted $\mathfrak{fivebrane}(n)$ 6-connections}
\label{G8 twist}
 
The discussion of twisted differential Fivebrane structures proceeds 
in direct analogy to the above discussion of twisted differential
string structures. One finds that the smooth second fractional 
Pontryagin class $\frac{1}{6}\mathbf{p}_2 : \mathbf{B}\mathrm{String}
\to \mathbf{B}^7 U(1)$ from the moduli 2-stack of $\mathrm{String}$-principal
2-bundles to the moduli 7-stack of circle 7-bundles is presented by
a span of simplicial presheaves of the form
$$
  \mathbf{B}\mathrm{String}
  \stackrel{\simeq}{\leftarrow}
  \mathbf{B}(\mathbf{B}^5 U(1) \to \mathrm{String})
  \to
  \mathbf{B}^7 U(1)
$$
whose infinitesimal version is a span of $L_\infty$-algebras
of the form
$$
  \mathfrak{string}
  \stackrel{\simeq}{\leftarrow}
  (b^5 \mathfrak{u}(1)
  \to 
  \mathfrak{fivebrane})
  \to 
  b^6 \mathfrak{u}(1)
  \,,
$$
where $\mathfrak{fivebrane} = (\mathfrak{so}_{\mu_3})_{\mu_7}$ 
denotes the Lie 6-algebra from example \ref{StringLie2Algebra},
and where the middle piece is the mapping cone of the defining extension
$$
  b^5 \mathfrak{u}(1) \to \mathfrak{fibrane} \to \mathfrak{string}
  \,.
$$
Therfore for $\mathrm{tw} : X \stackrel{\simeq}{\leftarrow} C(\{U_i\})
\to \mathbf{B}^7 U(1)$ a {\v C}ech cocycle for a twisting circle 7-bundle,
the corresponding \emph{twisted Fivebrane 6-bundles} are given by
lifts $\hat g$ in the disgram of simplicial presheaves
$$
  \xymatrix{
    & \mathbf{B}(\mathbf{B}^5 U(1) \to \mathrm{Fivebrane})
	\ar[dr]
    \\
    C(\{U_i\})
	\ar[d]^\simeq
	\ar[ur]^{\hat g}
	\ar[rr]^{\mathrm{tw}}
	&&
	\mathbf{B}^7 U(1)
	\\
	X
  }
  \,.
$$
As before, we consider now the analogous diagrams of local 
$L_\infty$-algebra valued forms in order to deduce the local
differental form data of twisted differential Fivebrane structures.
For transparency of the following diagrams we indicate both
the twist of the differential String-structure 
with local differential forms $(C_3)$ as before, as well as the
new twist of the differential Fivebrane structure by local
differential forms $(C_7)$. The former has to be taken to vanish,
but it is still instructive to display its differential incarnation,
for comparison with the previous case. 

Then the local differential cocycle data of a differential refinement
of the above twisted Fivebrane structure $\hat g$ is given over $k$-fold
intersections by a diagram of dg-algebras of the form.
  $$
  \xymatrix@C=30pt@R=40pt{
    &
    \mathrm{CE}((b \uu(1) \oplus b^5 \uu(1)) \hookrightarrow (\mathfrak{so}(n)_{\mu_3})_{\mu_7})
    \ar@{-->}[dl]|{\{A_{\mathrm{vert}}, (B_2)_{\mathrm{vert}}, (B_6)_{\mathrm{vert}}, (C_3)_{\mathrm{vert}}, (C_7)_{\mathrm{vert}}\}}
    \\
    \Omega^\bullet_{\mathrm{vert}}( U \times \Delta^k)
    &&
    \mathrm{CE}(b^2 \uu(1) \oplus b^6 \uu(1))
    \ar[ll]|>>>>>>>>>>>>>>{((C_3)_{\mathrm{vert}}, (C_7)_{\mathrm{vert}})}
    \ar@{_{(}->}[ul]
    \\
    &
     \mathrm{W}((b \uu(1) \oplus b^5\uu(1)) \hookrightarrow (\mathfrak{so}(n)_{\mu_3})_{\mu_7})
    \ar@{-->}[dl]|{\{A,F_A, B_2, B_6, \nabla B_2, \nabla B_6, C_3, C_7, G_4, G_8\}}
    \ar@{->>}[uu]|{\makebox(20,20){}}
    \\
    \Omega^\bullet(U \times \Delta^k)
    \ar@{->>}[uu]
    &&
    \mathrm{W}(b^2 \uu(1) \oplus b^6 \uu(1))
    \ar[ll]|>>>>>>>>>>>>>>>>{(C_3, C_7, G_4, G_8)}
    \ar@{_{(}->}[ul]
    \ar@{->>}[uu]
    \\
    &
    \mathrm{inv}(\mathrm{inn}(b \uu(1) \oplus b^5\uu(1)) \hookrightarrow \mathrm{cs}_{P_4 + P_8}(\mathfrak{so}(n)))
    \ar@{-->}[dl]|{\{H_3, H_7,G_4, G_8,P_4(F_A), P_8(F_A)\}}
    \ar@{->>}[uu]|{\makebox(20,20){}}
    \\
    \Omega^\bullet(U)
    \ar@{^{(}->}[uu]^{p^*}
    &&
    \mathrm{inv}(b^2 \uu(1) \oplus b^6 \uu(1))
    \ar[ll]_>>>>>>>>>>>>>>>{(G_4, G_8)}
    \ar@{_{(}->}[ul]
    \ar@{^{(}->}[uu]
  }
  \,.
$$
Here is again the meaning in words of the constituents of this diagram:
$$
  
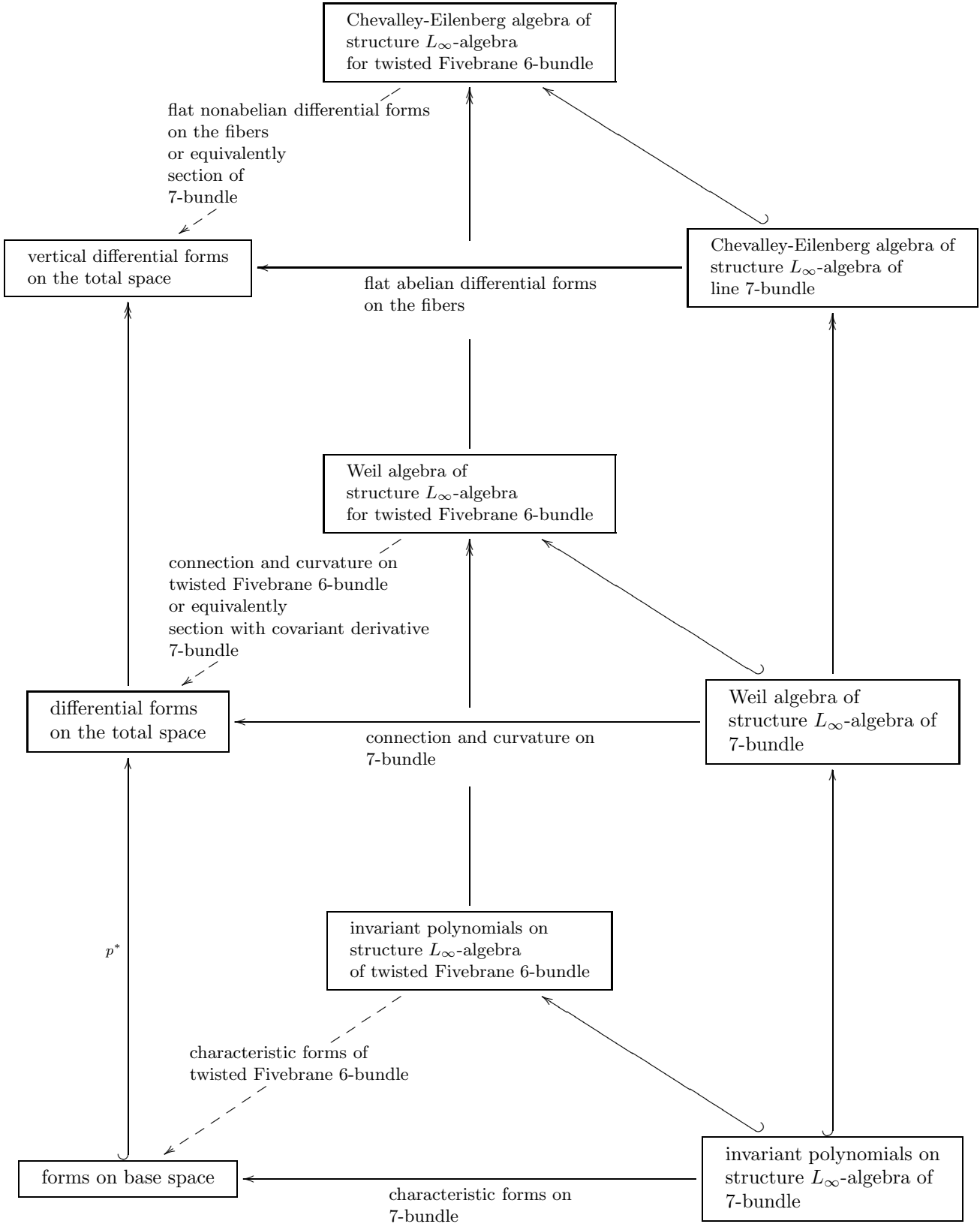
\hspace{-1cm}
  \xymatrix@C=30pt@R=68pt{
    &
    \fbox{
      \small
      \begin{tabular}{l}
        Chevalley-Eilenberg algebra of
        \\
        structure $L_\infty$-algebra
        \\
        for twisted Fivebrane 6-bundle
      \end{tabular}
    }
    \ar@{-->}[dl]|{
      \mbox{
        \small
        \begin{tabular}{l}
          flat nonabelian differential forms
          \\
          on the fibers
          \\
          or equivalently
          \\
          section of
          \\
          7-bundle
            \end{tabular}
      }
    }
    \\
    \fbox{
      \small
      \begin{tabular}{l}
        vertical differential forms 
        \\
        on the total space
      \end{tabular}
    }
    &&
    \fbox{
      \small
      \begin{tabular}{l}
        Chevalley-Eilenberg algebra of
        \\
        structure $L_\infty$-algebra of
        \\
        line 7-bundle
      \end{tabular}
    }
    \ar[ll]^{
      \mbox{
        \small
        \begin{tabular}{l}
          flat abelian differential forms
          \\
          on the fibers
        \end{tabular}
      }
    }
    \ar@{_{(}->}[ul]
    \\
    &
    \fbox{
      \small
      \begin{tabular}{l}
        Weil algebra of
        \\
        structure $L_\infty$-algebra
        \\
        for twisted Fivebrane 6-bundle
      \end{tabular}
    }
    \ar@{-->}[dl]|{
      \mbox{
        \small
        \begin{tabular}{l}
          connection and curvature on
          \\
          twisted Fivebrane 6-bundle
          \\
          or equivalently
          \\
          section with covariant derivative
          \\
          7-bundle
        \end{tabular}
      }
    }
    \ar@{->>}[uu]|>>>>>>>>>>>>>>>>>>>>>>>>>>>>>>>>>>>{\makebox(20,50){}}
    \\
    \fbox{
      \begin{tabular}{l}
        differential forms 
        \\
        on the total space
      \end{tabular}
    }
    \ar@{->>}[uu]
    &&
    \fbox{
      \begin{tabular}{l}
        Weil algebra of
        \\
        structure $L_\infty$-algebra of
        \\
        7-bundle
      \end{tabular}
    }
    \ar[ll]^{
      \mbox{
        \small
        \begin{tabular}{l}
          connection and curvature on
          \\
          7-bundle
        \end{tabular}
      }
    }
    \ar@{_{(}->}[ul]
    \ar@{->>}[uu]
    \\
    &
    \fbox{
      \small
      \begin{tabular}{l}
        invariant polynomials on
        \\
        structure $L_\infty$-algebra
        \\
        of twisted Fivebrane 6-bundle
      \end{tabular}
    }
    \ar@{-->}[dl]|{
      \mbox{
        \small
        \begin{tabular}{l}
          characteristic forms of
          \\
          twisted Fivebrane 6-bundle
        \end{tabular}
      }
    }
    \ar@{->>}[uu]|>>>>>>>>>>>>>>>>>>>>>>>>>>>>>>>>>>>>{\makebox(20,40){}}
    \\
    \fbox{
      \begin{tabular}{l}
        forms on base space
      \end{tabular}
    }
    \ar@{^{(}->}[uu]^{p^*}
    &&
    \fbox{
      \begin{tabular}{l}
        invariant polynomials on
        \\
        structure $L_\infty$-algebra of
        \\
        7-bundle
      \end{tabular}
    }
    \ar[ll]^{
      \mbox{
        \small
        \begin{tabular}{l}
          characteristic forms on
          \\
          7-bundle
        \end{tabular}
      }
    }
    \ar@{_{(}->}[ul]
    \ar@{^{(}->}[uu]
  }
$$
By again chasing elements through the diagram one finds the following data:
$$
  
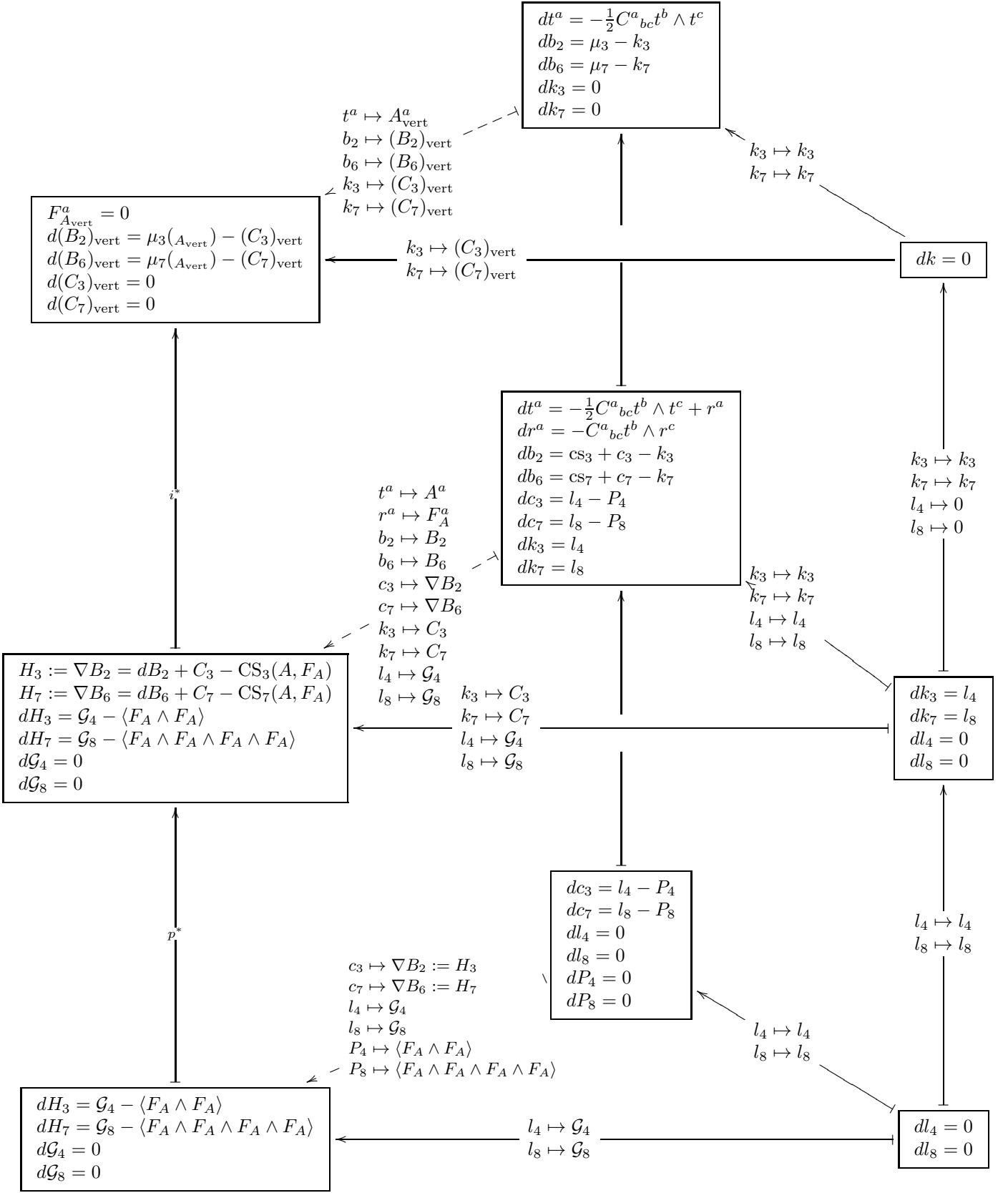
\hspace{-1cm}
  \xymatrix@C=75pt@R=30pt{
    &
    \fbox{$
      \begin{array}{l}
        dt^a = -\frac{1}{2}C^a{}_{bc}t^b \wedge t^c
        \\
        db_2 = \mu_3  - k_3
        \\
        db_6 = \mu_7  - k_7
        \\
        dk_3 = 0
        \\
        dk_7 = 0
      \end{array}
    $}
    \ar@{|-->}[dl]|{
      \begin{array}{l}
        t^a \mapsto A^a_{\mathrm{vert}}
        \\
        b_2 \mapsto (B_2)_{\mathrm{vert}}
        \\
        b_6 \mapsto (B_6)_{\mathrm{vert}}
        \\
        k_3 \mapsto (C_3)_{\mathrm{vert}}
        \\
        k_7 \mapsto (C_7)_{\mathrm{vert}}
      \end{array}
    }
    \\
    \fbox{$
      \begin{array}{l}
        F_{A_{\mathrm{vert}}}^a = 0
        \\
        d (B_2)_{\mathrm{vert}} = \mu_3(_{A_{\mathrm{vert}}}) - (C_3)_{\mathrm{vert}}
        \\
        d (B_6)_{\mathrm{vert}} = \mu_7(_{A_{\mathrm{vert}}}) - (C_7)_{\mathrm{vert}}
        \\
        d(C_3)_{\mathrm{vert}} = 0
        \\
        d(C_7)_{\mathrm{vert}} = 0
      \end{array}
    $}
    &
    &
    \fbox{$
      \begin{array}{l}
        dk = 0
      \end{array}
    $}
    \ar@{->}[ll]|>>>>>>>>>>>>>>>>>>>>>>>>>{
      \begin{array}{l}
        k_3 \mapsto (C_3)_{\mathrm{vert}}
        \\
        k_7 \mapsto (C_7)_{\mathrm{vert}}
      \end{array}
    }
    \ar@{->}[ul]|{
      \begin{array}{l}
        k_3 \mapsto k_3
        \\
        k_7 \mapsto k_7
      \end{array}
    }
    \\
    &
    \fbox{$
      \begin{array}{l}
        d t^a = - \frac{1}{2}C^{a}{}_{bc}t^b \wedge t^c + r^a
        \\
        d r^a = - C^a{}_{bc}t^b \wedge r^c
        \\
        db_2 = \mathrm{cs}_3 + c_3 - k_3
        \\
        db_6 = \mathrm{cs}_7 + c_7 - k_7
        \\
        dc_3 = l_4 - P_4
        \\
        dc_7 = l_8 - P_8
        \\
        dk_3 = l_4
        \\
        dk_7 = l_8
      \end{array}
    $}
    \ar@{|-->}[dl]|<<<<<<<<<<<<<<<<{
      \begin{array}{l}
        t^a \mapsto A^a
        \\
        r^a \mapsto F_A^a
        \\
        b_2 \mapsto B_2
        \\
        b_6 \mapsto B_6
        \\
        c_3 \mapsto \nabla B_2
        \\
        c_7 \mapsto \nabla B_6
        \\
        k_3 \mapsto C_3
        \\
        k_7 \mapsto C_7
        \\
        l_4 \mapsto {\cal G}_4
        \\
        l_8 \mapsto {\cal G}_8
      \end{array}
    }
    \ar@{|->}[uu]|<<<<<<<<<<<<<<<<<<<<<<<<<<{
      \makebox(20,20){}
    }
    \\
     \fbox{$
     \begin{array}{l}
       H_3 := \nabla B_2 = dB_2 + C_3 - \mathrm{CS}_3(A,F_A)
       \\
       H_7 := \nabla B_6 = dB_6 + C_7 - \mathrm{CS}_7(A,F_A)
       \\
       d H_3 = {\cal G}_4 - \langle F_A \wedge F_A \rangle
       \\
       d H_7 = {\cal G}_8 - \langle F_A \wedge F_A \wedge F_A \wedge F_A\rangle
       \\
       d {\cal G}_4 = 0
       \\
       d {\cal G}_8 = 0
     \end{array}
    $}
    \ar@{|->}[uu]|{
      i^*
    }
    &
    &
    \fbox{$
     \begin{array}{l}
       dk_3 = l_4
       \\
       dk_7 = l_8
       \\
       dl_4 = 0
       \\
       dl_8 = 0
     \end{array}
    $}
    \ar@{|->}[uu]|{
      \begin{array}{l}
        k_3 \mapsto k_3
        \\
        k_7 \mapsto k_7
        \\
        l_4 \mapsto 0
        \\
        l_8 \mapsto 0
      \end{array}
    }
    \ar@{|->}[ll]|>>>>>>>>>>>>>>>>>>>>>>>>>>{
      \begin{array}{l}
        k_3 \mapsto C_3
        \\
        k_7 \mapsto C_7
        \\
        l_4 \mapsto {\cal G}_4
        \\
        l_8 \mapsto {\cal G}_8
      \end{array}
    }
    \ar@{|->}[ul]|{
      \begin{array}{l}
        k_3 \mapsto k_3
        \\
        k_7 \mapsto k_7
        \\
        l_4 \mapsto l_4
        \\
        l_8 \mapsto l_8
      \end{array}
    }
  \\
  &
    \fbox{$
     \begin{array}{l}
       dc_3 = l_4-P_4
       \\
       dc_7 = l_8-P_8
       \\
       dl_4 = 0
       \\
       dl_8 = 0
       \\
       dP_4 = 0
       \\
       dP_8 = 0
     \end{array}
    $}
    \ar@{|-->}[dl]|>>>>>>>>>>>>>>>>>>>>>>>>>>>>>{\small\begin{array}{l}
      c_3 \mapsto \nabla B_2 := H_3
      \\
      c_7 \mapsto \nabla B_6 := H_7
      \\
      l_4 \mapsto {\cal G}_4
      \\
      l_8 \mapsto {\cal G}_8
      \\
      P_4 \mapsto \langle F_A \wedge F_A\rangle
      \\
      P_8 \mapsto \langle F_A \wedge F_A\wedge F_A\wedge F_A\rangle
    \end{array}}
    \ar@{|->}[uu]|>>>>>>>>>>>>>>>>>>>>>>>>>>{\makebox(20,20){}}
  \\
    \fbox{$
     \begin{array}{l}
       d H_3 = {\cal G}_4 - \langle F_A \wedge F_A \rangle
       \\
       d H_7 = {\cal G}_8 - \langle F_A \wedge F_A \wedge F_A\wedge F_A\rangle
       \\
       d {\cal G}_4 = 0
       \\
       d {\cal G}_8 = 0
     \end{array}
     $}
    \ar@{|->}[uu]|{
      p^*
    }
    &&
    \fbox{$
     \begin{array}{l}
       dl_4 = 0
       \\
       dl_8 = 0
     \end{array}
    $}
    \ar@{|->}[ll]|{
      \begin{array}{c}
        l_4 \mapsto {\cal G}_4
        \\
        l_8 \mapsto {\cal G}_8
      \end{array}
    }
    \ar@{|->}[ul]|{
      \begin{array}{l}
       l_4 \mapsto l_4
       \\
       l_8 \mapsto l_8
      \end{array}
    }
    \ar@{|->}[uu]|{
      \begin{array}{l}
       l_4 \mapsto l_4
       \\
       l_8 \mapsto l_8
      \end{array}
    }
  }
$$

Here, $P_4, P_8 \in W(\gg)$ are the invariant polynomials
on $\gg$ in transgression with with the cocycles $\mu_3, \mu_7 \in
\mathrm{CE}(\gg)$. 
The covariant derivative 7-form $\nabla B_6$ of the twisted
$(\mathfrak{so}(n)_{\mu_3})_{\mu_7}$-connection which we denote by
$H_7$ measures the difference between the prescribed $b^6
\uu(1)$-connection and the twist of the chosen twisted
$(\mathfrak{so}(n)_{\mu_3})_{\mu_7}$-connection. The
Bianchi identity
\(
  d H_7 = {\cal G}_8 - P_8(F_A)
\)
which appears in the middle on the left says that this difference
has to vanish in cohomology, being the local incarnation of the
differential refinement
of the anomaly cancellation condition discussed in section 
\ref{twisted fivebrane}.

\appendix
\section{$L_\infty$-algebraic notions}
\label{LinfinityNotions}

We collect here some $L_\infty$-algebraic definitions and constructions 
that are referred to 
in section \ref{Twisted differential String- and Fivebrane structures}.
Most of the following can be found in more detail in \cite{SSS1},
the main point here being the notion of representations of 
$L_\infty$-algebroids and the emphasis of the interpretation of 
twists as sections, as explained for smooth $\infty$-groupoids in 
section \ref{differential twisted cohomology}.
For a more conceptual account of $L_\infty$-algebra in the 
context of smooth $\infty$-groupoids see \cite{nactwist}.

\subsection{$L_\infty$-Algebras and $L_\infty$-Algebroids}
In direct generalization of how Lie algebras are infinitesimal
approximations to Lie groups, \emph{$L_\infty$-algebras} are
infinitesimal approximations to smooth $\infty$-groups. More
generally, \emph{$L_\infty$-algebroids} are infinitesimal
approximations to smooth $\infty$-groupoids.

\begin{definition}
  \label{Loo algebroid}
  An \emph{$L_\infty$-algebroid} $\mathfrak{a}$ of finite type
  over a smooth manifold $X$
  is  a non-positively graded 
  $A:= C^\infty(X)$-module degreewise of finite rank, 
  together with a degree +1 derivation
  \(
    d : \wedge^\bullet_{A} \mathfrak{a}^*
      \to
      \wedge^\bullet_{A} \mathfrak{a}^*\;,
  \)
  linear over the ground field (not necessarily over
  $A$)
  on the free (over $A$) graded-symmetric algebra
  generated from the $\mathbb{N}$-graded
  dual $\mathfrak{a}^*$ (over $A$),
  such that $d^2 = 0$.
The quasi-free (over $A$) differential graded-commutative algebra
  \(
    \mathrm{CE}(\mathfrak{a}) :=
    (\wedge^\bullet_{A} \mathfrak{a}^*,d)
  \)
  defined this way we call the
  \emph{Chevalley-Eilenberg algebra} of the $L_\infty$-algebroid
  $\mathfrak{a}$.
  
  We say the category $L_\infty \mathrm{Algd}$ of $L_\infty$-algebroids
  is the opposite of the full subcategory of dg-algebras
  on those of the above form.
\end{definition}
It is useful to distinguish the following special cases of this definition.
\begin{itemize}
  \item
    For $X = \mathrm{pt}$ and $\mathfrak{a}$ concentrated in degree 0
	on a vector space $\mathfrak{g}$
    we have $\mathrm{CE}(\mathfrak{a}) = (\wedge^\bullet \gg, d_\gg)$
    where $\wedge^\bullet \gg$ 
    is the Grassmann algebra on $\gg^*$ 
    $d_\mathfrak{a}$ is the Chevalley-Eilenberg differential uniquely 
    corresponding to the structure of a Lie algebra on $\gg$.
	$L_\infty$-algebroids arising this we we write 
	$\mathfrak{a} = b \mathfrak{g}$.
  \item
   For $X = \mathrm{pt}$ and $\mathfrak{a}$ concentrated 
   in arbitrary (non-positive) degree
   the above definition is that of an 
   \emph{$L_\infty$-algebra} structure (of finite type).
   For $\mathfrak{g}$ any $L_\infty$-algebra, we write
   $\mathfrak{a} = b \mathfrak{g}$ for the $L_\infty$-algebroid
   corresponding to it.
  \item For $X = \mathrm{pt}$ and $d$ at most co-binary (sending generators
    to wedge products of at most word length 2 in the generators) we have
	a \emph{dg-Lie algebra}.
  \item
   For arbitrary $X$ and $\mathfrak{a}$ concentrated in degree 0
   (being finitely generated and projective as a module over
   $C^\infty(X)$)
    this is equivalent to the usual
    definition of Lie algebroids as vector bundles $E \to X$ with anchor 
   map \cite{Mac} $\rho : E \to TX$: we have $\gg = \Gamma(E)$ and
   the anchor is encoded as
   $d_\gg|_{C^\infty(X)} : f \mapsto \rho(\cdot)(g) $.
  \item
    If $\mathfrak{a}$ is concentrated in degrees 0 through $-(n-1)$, then we
    speak of a {\it Lie $n$-algebroid} (Lie $n$-algebra if $X = {*}$).
\end{itemize}
\begin{example}[line Lie $n$-algebra]
  \label{LineLienAlgebra}
  For $n \in \mathbb{N}$ let $b^n \mathbb{R}$ or
equivalently $b^n \mathfrak{u}(1)$  
  be the Lie $n$-algebra defined by the fact that its 
  corresponding Lie $n$-algebroid $b^{n+1}\mathfrak{R}$ 
  has a Chevalley-Eilenberg algebra coming from a single generator
  in degree $(n+1)$ with vanishing differential. 
  We call this the \emph{line Lie $n$-algebra}.
\end{example}
\begin{example}[tangent Lie algebroid]
  For $X$ a smooth manifold, 
  the de Rham complex $\Omega^\bullet(X)$ is 
  the Chevalley-Eilenberg algebra of a Lie 1-algebroid over $X$, called
  the \emph{tangent Lie algebroid} $T X$,
  $
    \mathrm{CE}(T X) = (\Omega^\bullet(X), d_{\mathrm{dR}})
	\,.
  $
\end{example}
For $\mathfrak{g}$ an $L_\infty$-algebra, we say that a closed element
$\mu \in \mathrm{CE}(\mathfrak{g})$ of degree $n+1$ is an 
$(n+1)$-\emph{cocycle} on 
$\mathfrak{g}$. Cocycles are equivalently morphisms of $L_\infty$-algebras
$
  \mu : \mathfrak{g} \to b^n \mathbb{R}
$
or, equivalently, morphisms of $L_\infty$-algebroids
$
  \mu : b\mathfrak{g} \to b^{n+1} \mathbb{R}
 $.
\begin{example}[higher central extensions]
 \label{StringLInfinityAlgebras}
 \label{String-like extensions}
Every cocycle $\mu$ on an $L_\infty$-algebra 
$\mathfrak{g}$ induces a new $L_\infty$-algebra,
to be denoted $\mathfrak{g}_\mu$, which is defined by its $\mathrm{CE}$-algebra
being that of $\mathfrak{g}$ with a single generator $b$ in degree $n$ adjoined
and the differential extended to this generator by the formula
$
  d_{\mathfrak{g}_\mu} b = \mu
  \,.
$
This yields a sequence
\(
  b^n \mathbb{R} \to b \mathfrak{g}_\mu \to b\mathfrak{g}
\)
exhibiting the shifted central extension classified by $\mu$.
\end{example}
\begin{example}
  \label{StringLie2Algebra}
 For $\mathfrak{g}$ a semisimple Lie algebra and 
 $\mu_3 = \langle -,[-,-]\rangle$ the canonical 3-cocycle, the Lie
 2-algebra $\mathfrak{g}_{\mu_3}$ is called the 
 corresponding \emph{String Lie 2-algebra} \cite{H}\cite{BCSS}.
 For $\mathfrak{g} = \mathfrak{so}$ there is also the canonical 
 7-cocycle $\mu_7 \in \mathrm{CE}(\mathfrak{so})$. This is still a cocycle
 on $\mathfrak{so}_{\mu_3}$, too, and so there is a Lie 6-algebra
 $\mathfrak{so}_{\mu_3, \mu_7}$, called the
 \emph{fivebrane} Lie 6-algebra in \cite{SSS1}.
\end{example}

For brevity we state several constructions only for $L_\infty$-algebras.
The generalization to general $L_\infty$-algebroids is immediate. In
particular for $\mathfrak{g}$ an $L_\infty$-algebra we shall 
write $\mathrm{CE}(\mathfrak{g})$ as shorthand for $\mathrm{CE}(b\mathfrak{g})$.

Differential form data on a manifold $X$ 
with values in an $L_\infty$-algebra
$\mathfrak{g}$ is a graded algebra homomorphism (\emph{not necessarily} respecting the
differentials) from $\mathrm{CE}(\mathfrak{g})$ into the differential forms on $X$: 
\(
  \Omega^\bullet(X,\mathfrak{g})
  :=
  \mathrm{Hom}_{\mathrm{grAlg}}(\mathrm{CE}(\mathfrak{g}),\Omega^\bullet(X))
  \,.
\) 
The space of graded algebra homomorphisms is a subspace of the space
of linear maps of graded vector spaces from $\mathrm{CE}(\gg)$ to
$\Omega^\bullet(X)$ and, since $\mathrm{CE}(\gg)$ is freely
generated as a graded algebra and of finite type, this is isomorphic to the
space of grading preserving homomorphisms 
 $\mathrm{Hom}_{\mathrm{Vect}[{\mathbb{Z}]}}(\gg^*,\Omega^\bullet(X))$
from the graded vector space
$\gg^*$ of dual generators to $\Omega^\bullet(X)$. 
By the usual relation in
$\mathrm{Vect}[\mathbb{Z}]$ for $\gg$ of finite type, this is
isomorphic to the space of elements of total degree 1 in
forms tensored with $\gg$:
 \(
  \Omega^\bullet(X,\gg)
  \simeq
  (\Omega^\bullet(X)\otimes \gg)_0
  \label{omega g}
  \,.	
\) (Recall that $\gg$ is non-positively graded.)

If instead we  consider the corresponding homomorphisms 
of dg-algebras from $\mathrm{CE}(\gg)$
into forms, we find that respecting the differentials what deserves to be called
\emph{flatness} 
\(
  \Omega^\bullet_{\mathrm{flat}}(X,\gg) :=
  \mathrm{Hom}_{\mathrm{dgAlg}}(\mathrm{CE}(\gg),\Omega^\bullet(X))
  \,.
\)
The inclusion 
$  \xymatrix{
    \Omega^\bullet_{\mathrm{flat}}(X,\gg)
    =
    \mathrm{Hom}_{\mathrm{dgAlg}}(\mathrm{CE}(\gg),\Omega^\bullet(X))
    ~\ar@{^{(}->}[r]
    &
    \Omega^\bullet(X) \otimes \gg
  }
$
 realizes flat $L_\infty$-algebra valued forms as elements
$
  A \in \Omega^\bullet(X)\otimes \gg
$
of forms of total degree 0 with the special property that they
satisfy a flatness constraint of the form
    \(
      d A + \partial A  + [A \wedge A] + [A\wedge A\wedge A ]
      + \cdots = 0
      \,,
    \)
    where $d$ and $\wedge$ are the operations in $A \in \Omega^\bullet(X)\otimes
    \gg$ and where
 $[\cdot, \cdot, \cdots]$
 are the $n$-ary brackets in the $L_\infty$-algebra and $ \partial$ is the
    differential in the chain complex $\gg$.
    For $\gg$ a dg-Lie algebra, only the binary
    bracket is present and $A$ is an ordinary Maurer-Cartan element
    $
    D A +  [A \wedge A]  = 0
    $,
where $D= d  + \partial$.
This equation of course has a long and honorable history in various guises
called a \emph {Maurer-Cartan equation}. 

We would like to describe also non-flat $L_\infty$-algebra valued forms
by homomorphisms of
\emph{differential} graded algebras. This is accomplished by passing
from the Chevalley-Eilenberg algebra $\mathrm{CE}(\gg)$ to the Weil
algebra $\mathrm{W}(\gg)$.
\begin{definition}
  For $\mathfrak{g}$ an $L_\infty$-algebra, let $\mathrm{W}(\mathfrak{g})$
  be the unique dg-algebra free on the underlying graded vector space 
  $\mathfrak{g}^*$ such that the canonical projection morphism
  $
    W(\mathfrak{g}) \to \mathrm{CE}(\mathfrak{g})
  $
  is a dg-homomorphism.
\end{definition}
Due to the freeness of $\mathrm{W}(\mathfrak{g})$ we have 
an isomorphism
\(
  \Omega^\bullet(X, \mathfrak{g}) 
   = 
   \mathrm{Hom}_{\mathrm{grAlg}}(\mathrm{CE}(\mathfrak{g}), \Omega^\bullet(X))
   \simeq
   \mathrm{Hom}_{\mathrm{dgAlg}}(\mathrm{W}(\mathfrak{g}), \Omega^\bullet(X))
   \,.
\)
\begin{definition}
  \label{expgconngroupoid}
  The $\infty$-groupoid of $\mathfrak{g}$-valued forms on 
  a smooth manifold $X$ 
  is the Kan complex whose $k$-simplices are $\mathfrak{g}$-valued forms
  $
    \Omega^\bullet_{\mathrm{si}}(X \times \Delta^k)
	\leftarrow
	\mathrm{W}(\mathfrak{g})
	:
	A
  $
  on $X \times \Delta^k$ with \emph{sitting instants} 
  (becoming constant perpendicularly towards the faces of the $k$-simplex)
  and fitting into commutative diagrams of dg-algebras of the form
  \(
    \raisebox{30pt}{
    \xymatrix{
	  \Omega^\bullet_{\mathrm{vert}}(X \times \Delta^k)
	  \ar@{<-}[r]^{\hspace{0.8cm}A_{\mathrm{vert}}}
	  &
	  \mathrm{CE}(\mathfrak{g})
	  \\
	  \Omega^\bullet(X \times \Delta^k)
	  \ar@{<-}[r]^{\hspace{0.8cm}A}
	  \ar[u]
	  &
	  \mathrm{W}(\mathfrak{g})
	  \ar[u]
	  \\
	  \Omega^\bullet(X)
	  \ar@{<-}[r]
	  \ar@{^{(}->}[u]
	  &
	  \mathrm{inv}(\mathfrak{g})
         \ar@{^{(}->}[u]
	}
	}\;,
  \)
where the vertical morphisms are the canonical ones.
For more details on this definition see \cite{FSSI}.
The parameterized version of this construction leads to the 
smooth $\infty$-groupoid $\exp(\mathfrak{g})_{\mathrm{conn}}$
defined in (\ref{expgconn}).
\end{definition}

\subsection{$L_\infty$-algebra representations and section}
\label{L reps}

The $L_\infty$-analog of the notion of \emph{representations}
of smooth $\infty$-groupoids as in (\ref{FiberSequenceAndRepresentation})
is the following.
\begin{definition}[representations of $L_\infty$-algebroids]
\label{representation}
  A representation of an $L_\infty$-algebroid
  $\mathfrak{a}$ over $X$ on a cochain complex $V$ of finite rank 
  $(A := C^\infty(X))$-modules
  is an $L_\infty$-algebroid $V//_\rho\mathfrak{a}$ whose
  Chevalley-Eilenberg algebra $\mathrm{CE}(V//_\rho\mathfrak{a})$
  is an extension of $\mathrm{CE}(\mathfrak{a})$ by 
  $\wedge^\bullet_{\mathrm{Sym}_{A} V_0^*} V^*$
  \(
    \xymatrix{
       \wedge^\bullet_{\mathrm{Sym}_{A}V_0^*} V^*
       &&
       \mathrm{CE}(V//_\rho\mathfrak{a})
       \ar@{->>}[ll]
       &&
       \mathrm{CE}(\mathfrak{a})
       \ar@{_{(}->}[ll]
       \ar@/^2pc/[llll]^0
    }
	\,.
  \)
\end{definition}
\noindent This means that the differential is
 $d|_{\mathfrak{a}^*}  = d_\mathfrak{a}$, 
  $d|_{V^*} = d_{V^*} + d_\rho$,
where $ d_\rho$  
encodes the action of $\gg$ on $V$.

\begin{remark}
In roughly this latter form, the definition appears in \cite{BlockI}, where it is
called a \emph{superconnection}. Indeed, in cases where the
$L_\infty$-algebroid in question is similar to a tangent Lie
algebroid of some space,
its representations behave like (flat)
connections on that space. 
In the work of \cite{Crainic}, for the special case of 1-Lie algebroids,
such representations are called \emph{representations up to homotopy}.

For $\mathfrak{a} = b \mathfrak{g}$ coming from an 
$L_\infty$-algebra, the above notion of representation
reproduces the notion of  \emph{sh-representations} of \cite{Jim} \cite{LM}. 
\end{remark}

\noindent{\bf Example: ordinary adjoint representation.} Let $\gg$
be an ordinary Lie algebra with basis $\{t_a\}$ and structure
constants $\{C^a{}_{bc}\}$. Write $\{\underbrace{t^a}_{\mathrm{deg}
= +1}\}$ for the corresponding dual basis elements and
$\{\underbrace{\chi^a}_{\mathrm{deg} = 0}\}$ for the corresponding
basis elements of $V_\gg^*$. Then we have 
$d_\rho \chi^a
  =
  \sigma^{-1}(d_\gg t^a)
  =
  \sigma^{-1}(-\frac{1}{2}C^a{}_{bc}t^b \wedge t^c)
  =
  C^a{}_{bc}t^b \chi^c
  $.

\begin{definition}[sections and covariant derivatives]
\label{SectionsAndCovariantDerivatives}
Let $\gg$ be an $L_\infty$-algebra, and let $V//_\rho \mathfrak{g}$ be a
representation of $\gg$, def. \ref{representation}.

Then for $A \in \Omega\bullet(X \times \Delta^k, \mathfrak{g})$
a $k$-morphism of $\mathfrak{g}$-valued form data on  a smooth
manifold $X$, def. \ref{expgconngroupoid}, 
we say that a \emph{section} of the associated
$V$-connection is a choice of the dotted arrows in
\(
  \xymatrix{
    &
    \mathrm{CE}_\rho(\gg,V)
    \ar@{..>}[dl]|{(s,A_{\mathrm{vert}})}
    \\
    \Omega^\bullet_{\mathrm{vert}}(X \times \Delta^k)
    &&
    \mathrm{CE}(\gg)
    \ar[ll]|<<<<<<<<<<{A_{\mathrm{vert}}}
    \ar@{_{(}->}[ul]
    \\
    &
    \mathrm{W}_\rho(\gg,V)
    \ar@{..>}[dl]|{(s,\nabla_A s, A, F_A)}
    \ar@{->>}[uu]|{\makebox(12,12){}}
    \\
    \Omega^\bullet(X \times \Delta^k)
    \ar@{->>}[uu]
    &&
    \mathrm{W}(\gg)
    \ar[ll]|<<<<<<<<<<{(A,F_A)}
    \ar@{->>}[uu]
    \ar@{_{(}->}[ul]
    \\
    &
    \mathrm{inv}_\rho(\gg,V)
    \ar@{..>}[dl]
    \ar@{^{(}->}[uu]|{\makebox(12,12){}}
    \\
    \Omega^\bullet(X)
    \ar@{^{(}->}[uu]
    &&
    \mathrm{inv}(\gg)\;.
    \ar[ll]
    \ar@{^{(}->}[uu]
    \ar@{_{(}->}[ul]
  }
\)
\label{def cov deriv}
\end{definition}
\noindent Here we say that
$s$ is the section itself
whereas $\nabla_A s$ is its \emph{covariant derivative}.

\begin{example}
Let $\gg$ be an ordinary Lie algebra with Lie group $G$, let $V$ be
a vector space (a chain complex concentrated in degree 0) and $\rho$
an ordinary representation of $\gg$ on $V$, let $P$ be a principal
$G$-bundle and $(A,F_A)$ an ordinary Cartan-Ehresmann connection on
$P$. Then the dotted morphism in \(
  \xymatrix{
    &
    \mathrm{CE}_\rho(\gg,V)
    \ar@{..>}[dl]|{(s,A_{\mathrm{vert}})}
    \\
    \Omega^\bullet_{\mathrm{vert}}(P)
    &&
    \mathrm{CE}(\gg)
    \ar[ll]|<<<<<<<<<<{A_{\mathrm{vert}}}
    \ar@{_{(}->}[ul]
    \\
  }
\)
is dual to a $V$-valued function on the total space of the bundle
(not on base space!)
$
  s : P \to V
 $,
which is covariantly constant along the fibers in that the covariant derivative
\(
  \nabla_A s := d s + (\rho\circ A) s
\)
 vanishes when evaluated on vertical vectors,
 where $(\rho\circ A)s$ denotes the action of $A$ on the section $s$
using the representation $\rho$.
 This means that
$s$ descends to a section of the associated vector bundle $P
\times_G V$. The covariant derivative 1-form $\nabla_A s$ of
the section $s$ is one component of the extension in the middle
part of our diagram \(
  \xymatrix{
    &
    \mathrm{W}_\rho(\gg,V)
    \ar@{..>}[dl]|{(s,\nabla_A s,A,F_A)}
    \\
    \Omega^\bullet(P)
    &&
    \mathrm{W}(\gg)
    \ar[ll]|<<<<<<<<<<{A_{\mathrm{vert}}}
    \ar@{_{(}->}[ul]
    \\
  }
  \,.
\) The equation \(
  \nabla_A \nabla_A s = (\rho\circ F_A)\wedge s
\) is the Bianchi identity for $\nabla_A s$. If $s$ is
everywhere non-vanishing, this says that the curvature $F_A$ of our
bundle is covariantly exact on $P$. In the case that $\gg = \uu(1)$
it follows that $F_A$ is an exact 2-form on $P$ and the choice of
the non-vanishing section amounts to a trivialization of the bundle.
\end{example}

In sections \ref{G4 twist} and \ref{G8 twist} 
we see twisted differential
String-structures and twisted differential Fivebrane structures as 
parameterized examples of this notion of $L_\infty$-sections.

\bigskip\bigskip
\noindent
{\bf \large Acknowledgements}

\vspace{2mm}
\noindent H. S. and U. S. would like to thank the Hausdorff Institute for Mathematics in Bonn
for hospitality and the organizers of the ``Geometry and Physics" Trimester Program
at HIM for the inspiring atmosphere during the initial stages of this project. H.S. thanks Matthew
Ando for useful discussions. U. S. also thanks the Max-Planck institute for Mathematics in 
Bonn for hospitality later during this work; and the crew of the $n$Lab, where some of
the material presented here was first exposed. This research is supported in parts
by the FQXi mini-grant ``QFT and Nonabelian Differential Cohomology'' and NSF Grant PHY-1102218.
H. S. thanks the Department of Mathematics at Hamburg University for hospitality during the writing of this paper.
J. S. would like to thank the Department of Mathematics of the University of Pennsylvania for support of the 
Deformation Theory Seminar enabling the three authors to have at least one meeting in person.
The authors are indebted to the referee for many useful remarks and suggestions that led to major
improvements of the paper.


%




%



\begin{thebibliography}{99}

\bibitem{Crainic}
C. Abad and  M. Crainic, {\it Representations up to homotopy of Lie algebroids },
[{\tt arXiv:0901.0319}] [{\tt math.DG}].




\bibitem{AJ}
P. Aschieri and B. Jur{\v c}o,
{\it Gerbes, M5-brane anomalies and $E_8$ gauge theory},
J. High Energy Phys. {\bf 0410} (2004) 068,
[{\tt arXiv:hep-th/0409200}].

%

\bibitem{BBK}
N. A. Baas, M. B{\"o}kstedt, and T. A. Kro,
{\it Two-categorical bundles and their classifying spaces}, 
[{\tt arXiv:math/0612549}] [{\tt math.AT}]. 
%

\bibitem{BCSS} J. Baez, A. Crans, U. Schreiber,
   and D. Stevenson, {\it From loop groups to 2-groups},
    Homology, Homotopy Appl. {\bf  9} (2007), no. 2, 101- 135,
  [{\tt arXiv:math/0504123}] [{\tt math.QA}].

\bibitem{BaezSchr} J. Baez, and U. Schreiber,
  {\it Higher gauge theory},
   Categories in Algebra, Geometry and Mathematical Physics,
   7--30, Contemp. Math., {\bf 431}, Amer. Math. Soc., Providence, RI, 2007,
[{\tt arXiv:math/0511710v2}] [{\tt math.DG}].

%

\bibitem{BSt}
J. Baez and D. Stevenson, 
{\it The classifying space of a topological 2-group},
[{\tt arXiv:0801.3843}] [{\tt math.AT}].

\bibitem{BehrendXu}
K.~Behrend, P.~Xu, {\it Differentiable Stacks and Gerbes},
J. Symplectic Geom. {\bf 9} (2011), 285--341,
[{\tt arXiv:math/0605694}].

\bibitem{BlockI}
J. Block,
{\it Duality and equivalence of module categories in noncommutative geometry I},
[{\tt  arXiv:math/0509284}] [{\tt math.QA}].


\bibitem{Bott}
R.~Bott,
\newblock {\it The space of loops on a Lie group},
\newblock Michigan Math. J., 5:35–61, (1958)


 \bibitem{BCMMS}
 P. Bouwknegt , A. L. Carey, V. Mathai, M. K. Murray, and D. Stevenson,
 {\it  Twisted K-theory and K-theory of bundle gerbes},
 Commun. Math. Phys. {\bf 228} (2002) 17-49,
 [{\tt arXiv:hep-th/0106194v2}].

\bibitem{BM}
P. Bouwknegt and V. Mathai,
{\it D-branes, B-fields and twisted K-theory},
J. High Energy Phys. {\bf 03} (2000) 007,
[{\tt arXiv:hep-th/0002023}].

\bibitem{BMRS}
J. Brodzki, V. Mathai, J. Rosenberg, and R. J. Szabo,
{\it D-branes, RR-fields and duality on noncommutative manifolds},
Commun. Math. Phys. {\bf 277}, no.3 (2008) 643-706,
[{\tt 	arXiv:hep-th/0607020v3}].

\bibitem{Bunke}
U.~Bunke,
\newblock {\it String structures and trivialisations of a Pfaffian line bundle},
\newblock {\tt arXiv:0909.0846}


\bibitem{Cham}
A. H. Chamseddine,
{\it Interacting supergravity in ten-dimensions: the role of the six-index gauge field},
Phys. Rev. {\bf D24} (1981) 3065.

\bibitem{Ch}
A. H. Chamseddine,
{\it A family of dual $N=2$ supergravity actions in ten-dimensions},
Phys. Lett. {\bf B367} (1996) 134-139,
[{\tt arXiv:hep-th/9510100}].












\bibitem{DFM}
E. Diaconescu, D. S. Freed and G. Moore, {\it The M-theory 3-form
and $E_8$ gauge theory}, Elliptic cohomology, 44--88,
London Math. Soc. Lecture Note Ser., 342, Cambridge Univ. Press, Cambridge, 2007,
[{\tt arXiv:hep-th/0312069}].


 \bibitem{DMW}
E.~Diaconescu, G.~Moore and E.~Witten,
{\it $E_8$ gauge theory, and a derivation of K-theory from M-theory},
Adv. Theor. Math. Phys. {\bf 6} (2003) 1031,
[{\tt arXiv:hep-th/0005090}].

\bibitem{DDP}
J. A. Dixon, M. J. Duff, and J. C. Plefka,
{\it Putting string/fivebrane duality to the test},
Phys. Rev. Lett. {\bf 69} (1992) 3009-3012,
[{\tt arXiv:hep-th/9208055v1}].



\bibitem{DLM}
M. J. Duff, J. T. Liu, and R. Minasian,
{\it Eleven-dimensional origin of string-string duality: A one loop test},
Nucl. Phys. {\bf B452} (1995) 261-282,
[{\tt arXiv:hep-th/9506126}].


\bibitem{ES2}
J. Evslin and H. Sati, 
{\it Can D-branes wrap nonrepresentable cycles?},
J. High Energy Phys. {\bf 0610} (2006), 050,
[{\tt 	arXiv:hep-th/0607045}].

\bibitem{FSSI}
D.~Fiorenza, U.~Schreiber, J.~Stasheff,
\newblock {\it {\v C}ech cocycles for differential characteristic classes -- an $\infty$-Lie theoretic construction}, 
\newblock [{\tt arXiv:1011.4735}]


 \bibitem{Freed}
D.~S.~Freed,
\newblock {\it Dirac charge quantization and generalized differential cohomology},
\newblock Surv. Differ. Geom., VII, 129--194, Int. Press, Somerville, MA, 2000,
\newblock [{\tt arXiv:hep-th/0011220}].

\bibitem{FHMM}
D. Freed, J. A. Harvey, R. Minasian, and G. Moore,
{\it Gravitational anomaly cancellation for M-theory fivebranes},
Adv. Theor. Math. Phys. {\bf 2} (1998) 601-618,
[{\tt arXiv:hep-th/9803205}].


\bibitem{FW}
D.~S.~Freed and E.~Witten, {\it Anomalies in string theory with
D-Branes}, Asian J. Math. {\bf 3} (1999) 819,
[{\tt arXiv:hep-th/9907189}].

\bibitem{Gates}
S.J. Gates, Jr., and H. Nishino,
{\it New $D=10$, $N=1$ supergravity coupled to Yang-Mills supermultiplet
and anomaly cancellations}, Phys. Lett. {\bf B157} (1985) 157.

\bibitem{GS}
M. B. Green and J. H. Schwarz,
{\it Anomaly cancellation in supersymmetric $D=10$ gauge theory and
superstring theory}, Phys. Lett. {\bf B149} (1984) 117.



\bibitem{H} A. Henriques,
{\it Integrating $L_\infty$-algebras},
Compos. Math. {\bf 144} (2008), no. 4, 1017--1045,
[{\tt arXiv:math/0603563}] [{\tt math.AT}].





\bibitem{HS}
M. J. Hopkins and I. M. Singer,
{\it Quadratic functions in geometry, topology, and M-theory}, J. Diff. Geom. {\bf 70}
(2005) 329-452,
[{\tt arXiv:math.AT/0211216}].


\bibitem{HW}
P. Horava and E. Witten,
{\it Eleven-dimensional supergravity on a manifold with boundary},
Nucl. Phys. {\bf B475} (1996) 94-114,
[{\tt arXiv:hep-th/9603142}].


\bibitem{IP}
C.J. Isham and C.N. Pope,
{\it Nowhere vanishing spinors and topological obstructions 
to the equivalence of the NSR and GS superstrings},
Class. Quant. Grav. {\bf 5} (1988) 257.

\bibitem{IPW}
C.J. Isham, C.N. Pope, and N.P. Warner, 
{\it Nowhere vanishing spinors and triality rotations in eight manifolds},
Class. Quant. Grav. {\bf 5} (1988) 1297.



 \bibitem{Kap}
A.~Kapustin, {\it D-branes in a topologically nontrivial B-field},
Adv.~Theor.~Math.~Phys.~{\bf 4} (2000) 127,
[{\tt arXiv:hep-th/9909089}].

\bibitem{Ki}
T. P. Killingback, {\it World-sheet anomalies and loop geometry},
Nucl. Phys. {\bf B288} (1987) 578.



\bibitem{KS}
L. Kramer and S. Stolz,
{\it A diffeomorphism classification of manifolds which are like projective planes},
J. Differential Geom. {\bf 77} (2007), no. 2, 177--188,
[{\tt arXiv:math/0505621v3}] [math.GT].


\bibitem{LM}
T. Lada and M. Markl,
{\it Strongly homotopy Lie algebras},
 Comm. Algebra  {\bf 23}  (1995), no. 6, 2147--2161.

 \bibitem{LT}
K. Lechner and M. Tonin,
{\it World volume and target space anomalies in the $D=10$
superfivebrane sigma model},
Nucl. Phys. {\bf B475} (1996) 545-561,
[{\tt arXiv:hep-th/9603094}].



\bibitem{Lurie} 
J.~Lurie, 
\newblock {\it Higher topos theory},
\newblock Annals of mathematics studies, no. 170, Princeton University Press, Princeton, NJ,	2009,
\newblock [{\tt arXiv:math/0608040}] [{\it math.CT}].





\bibitem{Mac}
K. C. H. Mackenzie,
{\it General theory of Lie groupoids and Lie algebroids},
Cambridge University Press, Cambridge, 2005.







\bibitem{NikolausSachseWockel}
T.~Nikolaus, C.~Sachse, C.~Wockel,
{\it A smooth model for the string 2-group},
[{\tt arXiv:1104.4288}]

\bibitem{NG}
H. Nishino and S. J. Gates, Jr.,
{\it Dual versions of higher dimensional supergravities and 
anomaly cancellations in lower dimensions},
Nucl. Phys. {\bf B268} (1986) 532--542.

\bibitem{Redden}
C. Redden {\it String structures and canonical 3-forms}, 
[{\tt 	arXiv:0912.2086}] [{\tt math.DG}].



\bibitem{SS}
A. Salam and E. Sezgin,
{\it Anomaly freedom in chiral supergravities},
Phys. Scripta {\bf 32} (1985) 283.







\bibitem{S1}
H. Sati, {\it Geometric and topological structures related to M-branes}, 
Proc. Symp. Pure Math. {\bf 81} (2010), 
181--236, [{\tt arXiv:1001.5020}] [{\tt math.DG}].
 
 \bibitem{S2}
H. Sati, {\it Geometric and topological structures related to M-branes II: 
Twisted String and String${}^c$ structures}, 
J. Australian Math. Soc. {\bf 90} (2011), 93--108,
[{\tt arXiv:1007.5429}] [{\tt hep-th}].

\bibitem{S3}
H. Sati, {\it Twisted topological structures related to M-branes},
Int. J. Geom. Meth. Mod. Phys. {\bf 8} (2011), 1097--1116,
	[{\tt arXiv:1008.1755}] [{\tt hep-th}].





\bibitem{SSS1} H. Sati, U. Schreiber and J. Stasheff,
{\it $L_\infty$-connections and applications to String- and
Chern-Simons $n$-transport}, in {\it Recent Developments in QFT},
eds. B. Fauser et al., Birkh{\"a}user, Basel (2008), [{\tt
arXiv:0801.3480}] [{\tt math.DG}].

\bibitem{SSS2}
H. Sati, U. Schreiber, and J. Stasheff,
{\it Fivebrane structures}, Rev. Math. Phys.
{\bf 21} (2009) 1--44,
[{\tt arXiv:math/0805.0564}] 
[{\tt math.AT}].

\bibitem{SchommerPries}
C. Schommer-Pries {\it A finite-dimensional model for the String 2-group},
[{\tt arXiv:0911.2483}] [{\tt math.AT]}.

\bibitem{nactwist} 
 U.~Schreiber.
 {\it Differential cohomology in a cohesive topos}, 
 Habilitation, Hamburg (2011)
 \newline
{{\small\tt{http://ncatlab.org/schreiber/show/differential+cohomology+in+a+cohesive+topos}}}

\bibitem{action}
U. Schreiber, {\it On $\infty$-Lie theory}, [{\tt
http://www.math.uni-hamburg.de/home/schreiber/action.pdf}].

 \bibitem{SWI} U. Schreiber and K. Waldorf,
{\it Parallel transport and functors},  [{\tt  arXiv:0705.0452}] [{\tt math.DG}].

  \bibitem{SWII} U. Schreiber and K. Waldorf,
{\it Smooth functors vs. differential forms},
  [{\tt  arXiv:0802.0663}] [{\tt math.DG}].

\bibitem{SWIII}
U. Schreiber and K. Waldorf,
{\it Connections on nonabelian gerbes and their holonomy},
   [{\tt arXiv:0808.1923}] [{\tt math.DG}].





\bibitem{Sing}
W. M. Singer,
{\it Connective fiberings over ${\rm BU}$ and ${\rm U}$},
Topology {\bf 7} (1968) 271--303. 

%
\bibitem{Jim}
J.~Stasheff,
{\it Constrained Poisson algebras and strong homotopy representations},
 Bull. Amer. Math. Soc. (N.S.)  {\bf 19}  (1988),  no. 1, 287--290.


\bibitem{ST1}
S. Stolz and P. Teichner, {\it What is an elliptic object?}, in Topology, geometry 
and quantum field theory, 247--343, Cambridge Univ. Press, Cambridge, 2004.


\bibitem{Stong}
R. Stong,
{\it Determination of $H\sp{*} ({\rm BO}(k,\cdots,\infty ),Z\sb{2})$\ and
$H\sp{*} ({\rm BU}(k,\cdots,\infty ),Z\sb{2})$},
Trans. Amer. Math. Soc. {\bf 107} (1963) 526--544.


\bibitem{Toen2} B. To{\"e}n,
{\it Notes on non-abelian cohomology},
  Lecture in MSRI, January 2002,
    \newline
    [{\tt http://www.math.univ-toulouse.fr/$\sim$toen/msri2002.pdf}].

\bibitem{VW}
C. Vafa and E. Witten,
{\it A one-loop test of string duality},
Nucl. Phys. {\bf B447} (1995) 261-270,
[{\tt arXiv:hep-th/9505053}].


\bibitem{Waldorf}
K. Waldorf, {\it String connections and Chern-Simons theory},
[{\tt 	arXiv:0906.0117}][{\tt math.DG}].

\bibitem{Wang}
B.-L. Wang,
{\it Geometric cycles, index theory and twisted K-homology},
J. Noncommut. Geom. {\bf 2} (2008), no. 4, 497--552,
[{\tt arXiv:0710.1625v1}] [{\tt math.KT}].

\bibitem{Witt}
F. Witt,
{\it Special metric structures and closed forms},
DPhil Thesis, University of Oxford, 2004. 



 \bibitem{Flux}
E. Witten,
{\it On flux quantization in M-theory and the effective action},
J. Geom. Phys. {\bf 22} (1997) 1-13,  \newline
[{\tt arXiv:hep-th/9609122}].

\bibitem{Effective}
E. Witten,
{\it Five-brane effective action in M-theory},
J. Geom. Phys. {\bf 22} (1997) 103-133,  
[{\tt arXiv:hep-th/9610234}].

\bibitem{Wi1}
E.~Witten,
{\it D-Branes and K-Theory},
J. High Energy Phys. {\bf 12} (1998) 019,
[{\tt arXiv:hep-th/9810188}].

\bibitem{Among}
E. Witten,
{\it Duality relations among topological effects in string theory},
J. High Energy Phys. {\bf 0005} (2000) 031,
[{\tt arXiv:hep-th/9912086}].


 \end{thebibliography}
\end{document}